\documentclass[preprint,12pt]{elsarticle}

\usepackage{amssymb}
\usepackage{amsmath,blkarray,booktabs,bigstrut}
\usepackage{mathtools}
\usepackage{xcolor}
\usepackage{algorithm}
\usepackage{algpseudocode}
\usepackage{siunitx}
\usepackage{mathabx}
\def\vr{\kern-\arraycolsep & \kern-\arraycolsep}
\def\VR{\kern-\arraycolsep\strut\vrule &\kern-\arraycolsep}
\newcommand{\red}[1]{\textcolor{black}{#1}}
\newcommand{\newred}[1]{\textcolor{black}{#1}}
\newcommand{\newredbis}[1]{\textcolor{black}{#1}}

\journal{Engineering Applications of Artificial Intelligence}

\begin{document}

\begin{frontmatter}



\title{Bayesian Quality-Diversity approaches for constrained optimization problems with mixed continuous, discrete and categorical variables}


\author[inst1]{Loïc Brevault}

\affiliation[inst1]{organization={ONERA / DTIS, Université Paris -- Saclay},%Department and Organization
            city={Palaiseau},
            country={France}}

\author[inst1]{Mathieu Balesdent}

\begin{abstract}
Complex system design problems, such as those involved in aerospace engineering, require the use of numerically costly simulation codes in order to  predict the performance of the system to be designed. In this context, these codes are often embedded into an optimization process to provide the best design while satisfying the design constraints. Recently, new approaches, called Quality-Diversity, have been proposed in order to enhance the exploration of the design space and to provide a set of optimal diversified solutions with respect to some feature functions. These functions are interesting to assess trade-offs. Furthermore, complex design problems often involve mixed continuous, discrete, and categorical design variables allowing to take into account technological choices in the optimization problem. \red{Existing Bayesian Quality-Diversity approaches suited for intensive high-fidelity simulations are not adapted to mixed variables constrained optimization problems. In order to overcome these limitations, a new Quality-Diversity methodology based on mixed variables Bayesian optimization strategy is proposed in the context of limited simulation budget. Using adapted covariance models and dedicated enrichment strategy for the Gaussian processes in Bayesian optimization, \newred{this approach allows to reduce the computational cost up to two orders of magnitude}, with respect to classical Quality-Diversity approaches while dealing with discrete choices and the presence of constraints. The performance of the proposed method is assessed on a benchmark of analytical problems as well as on two aerospace system design problems highlighting its efficiency in terms of speed of convergence. The proposed approach provides valuable trade-offs for decision-markers for complex system design.} 
\end{abstract}

\begin{graphicalabstract}
\includegraphics[width = 0.85\linewidth]{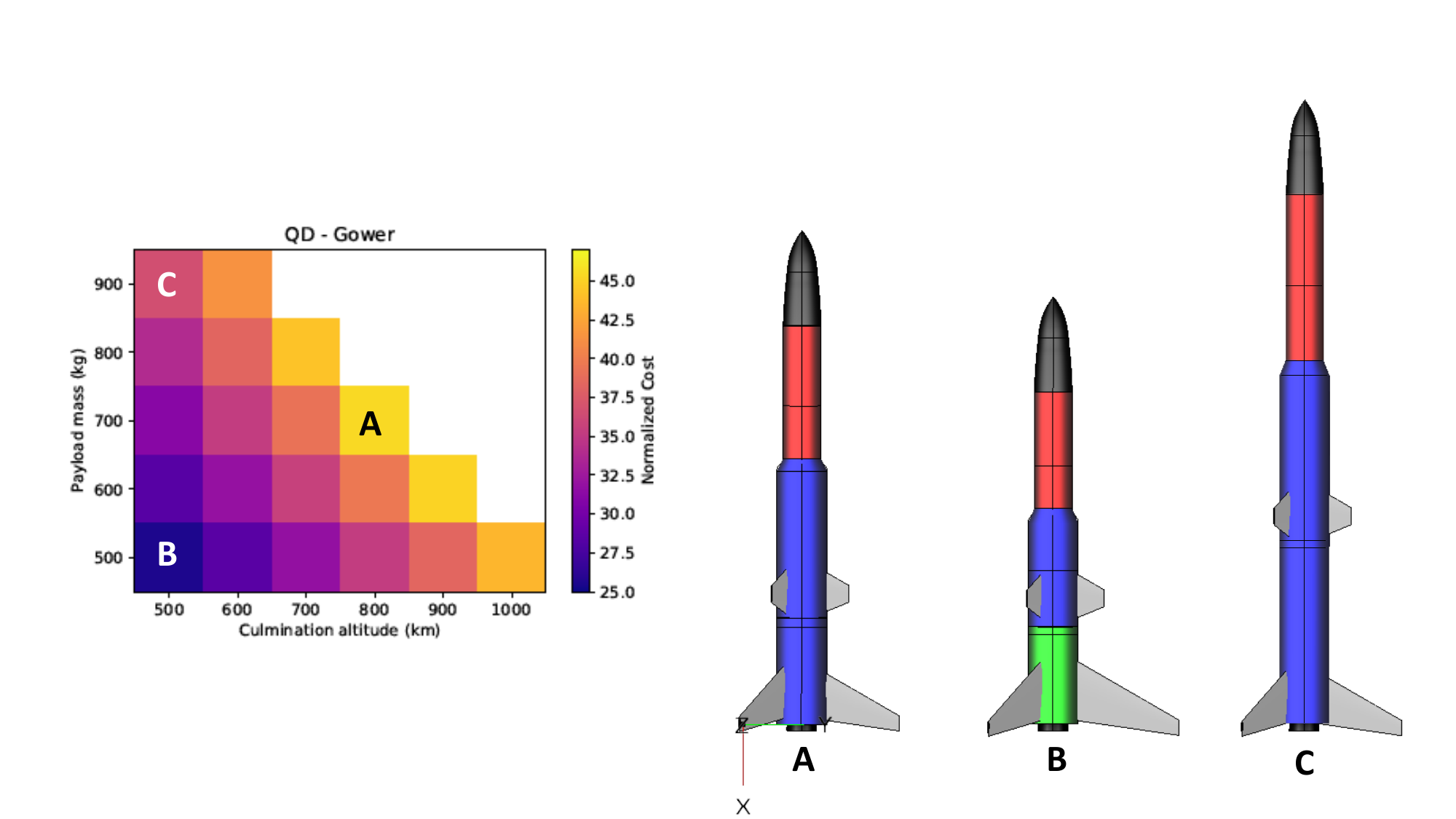}
\end{graphicalabstract}

\begin{highlights}
\item Development of a Bayesian Quality-Diversity approach to deal with constrained problems involving mixed continuous, discrete and categorical variables using adapted Gaussian processes
\item Extension of MAP-Elites algorithms to deal with mixed variables and constraints
\item Application on a benchmark of analytical problems and two aerospace design problems
\end{highlights}


\end{frontmatter}

\section{Introduction}
Complex system design problems, such as those involved in aerospace \cite{keane2005computational,adler2023hydrogen}, civil  \cite{mei2021structural}, or energy  \cite{al2022towards} engineering fields, require the use of simulation codes in order to numerically predict the behavior and performance of the system to be designed. In the different design phases of the system, these codes can be used in optimization problems in order to minimize a given objective function (\textit{e.g.}, mass, cost, consumption) with respect to some design variables characterizing the system while taking into account several specifications (\textit{e.g.}, reliability, environmental impact) as constraints. In the early phases of the design process, one aims at exploring a large design space in order to assess the best architecture for the given mission \cite{balesdent2023multidisciplinary}. This induces the handling of classical continuous design variables but also discrete variables (\textit{e.g.,} number of engines) and categorical variables that represent the different available technology options (\textit{e.g.}, for a launch vehicle, type of propulsion - solid or liquid, for a wind turbine, type of material - composite or metallic). The introduction of continuous, discrete and categorical variables presents some challenges in optimization algorithms dealing with computationally expensive engineering simulation codes. \red{For the design of complex systems, high-fidelity simulations such as Computational Fluid Dynamics for aerodynamics or Finite Element Analysis for structural assessment are necessary. These simulation models may take from few minutes up to hours of calculation despite intensive use of clusters and parallelization. In this context, in this paper, it is assumed that a limited evaluation budget of the simulation model is affordable (in the order of few hundreds of evaluations). The family of surrogate-assisted optimization algorithms is particularly suited in the context of limited evaluation budget \cite{forrester2009recent,garnett2023bayesian}.} 
This type of approaches consists in substituting the computationally costly simulation codes by surrogate models (\textit{e.g.}, neural networks \cite{du2021rapid}, Gaussian process \cite{gramacy2020surrogates}, support vector machines \cite{stromberg2020efficient}) and using them in the optimization process. One popular way to perform such optimization strategies is to use Gaussian Processes (GPs) of the objective function and constraints that are enriched all along the optimization process (via active learning) in order to converge to the optimum while exploring the design space \cite{chaiyotha2020comparative}. This is often called Bayesian optimization \cite{garnett2023bayesian,le2021revisiting}.\\

In the early design phase of innovative systems, trade-offs have often to be assessed between several quantities of interest that characterize the performance of the system (\textit{e.g.}, costs, environmental impact, availability rate). In such a context, single-objective optimization approaches are not sufficient and two different types of methodologies can be used. The first deals with multi-objective algorithms in which the different quantities of interests are gathered into a vector of objective functions that need to be optimized together. The result is a Pareto front that allows the decision makers to assess trade-offs between different solutions \cite{gunantara2018review,brevault2020multi}.
Recently, another family of optimization strategies, called Quality-Diversity \cite{chatzilygeroudis2021quality,pugh2016quality}, has been proposed. These algorithms optimize the objective function while exploring some quantities of interest, called features. \red{Features correspond to a set of system characteristics or preferences that may be useful to establish a trade-off by the decision-makers}. Such algorithms provide the decision makers with a set (an archive) of diversified solutions that aims at promoting diversity with respect to the features while optimizing a single objective function. Compared to the multi-objective optimization, Quality-Diversity (QD) algorithms allow to improve the exploration of the search space (especially in the case of non antagonistic behavior between the objective and the features) and preserve the optimization with respect to the objective function. The main QD algorithms are derived from population-based optimization techniques. The  MAP-Elites algorithm  \cite{mouret2015illuminating} \newred{(Multi-dimensional Archive of Phenotypic Elites)} is one of the most popular methods in this family. When dealing with computationally expensive simulation codes and therefore a limited simulation budget, several Bayesian QD algorithms have been proposed, such as Surrogate-Assisted ILlumination (SAIL) \cite{gaier2017aerodynamic}, Surrogate-assisted PHEnotypic Niching (SPHEN) \cite{hagg2020designing} or Bayesian Optimization of Elites (BOP-Elites) \cite{kent2020bop,kent2023bop}. \red{However, the existing algorithms within the Bayesian QD family (SAIL, SPHEN, BOP-Elites) are dedicated to continuous unconstrained optimization problems. There exist non-Bayesian QD algorithms dedicated to mixed continuous, discrete and categorical variables \cite{boige2023gradient}, however, the required number of evaluations of the computationally intensive objective, features or constraints functions is not affordable in practice for limited evaluation budget. Therefore, there is a need for Bayesian QD to solve mixed variables constrained problems with a limited simulation budget.}

\red{The contributions of this paper are the development and evaluation of a new Bayesian Quality-Diversity approach in order to handle of computationally intensive objective, constraints and feature functions as well as mixed continuous, discrete and categorical variables. In that purpose, the Gaussian processes involved in Bayesian QD optimization are adapted with specific covariance functions \cite{pelamatti2019efficient,saves2022general} in order to handle mixed variables. Furthermore, a new active learning strategy (with a suited infill criterion) is developed in order to refine the surrogate models along the optimization convergence while promoting both quality and diversity. Eventually, to carry out the optimization of the infill criterion, an adaptation of MAP-Elites approaches is performed to handle discrete and categorical variables as well as the constraints. The interest of the proposed approach lies in its ability to carry out QD optimization for constrained problems involving high fidelity simulation models in a context of limited evaluation budget. Moreover, the proposed approach is adapted for complex system design involving discrete and categorical optimization variables.}
The proposed approach is compared with classical mixed continuous / discrete version of MAP-Elites on three analytical test-cases and \red{two aerospace engineering problems. The considered performance metrics for the comparison correspond to the Quality-Diversity score \cite{pugh2016quality} and the number of discovered alternatives (defined by the features)}. \\

The paper is organized as follows. In Section \ref{sec:QD}, the QD optimization problem is described to explain the main differences with classical optimization problem. In Section \ref{sec:QD_BO}, the proposed Bayesian QD approach and the general algorithm are presented, with two versions depending on the covariance model used in the Gaussian Process. Eventually, in Section \ref{sec:numerical_expe}, the proposed approach is compared to classical mixed continuous / discrete versions of MAP-Elites on three analytical analytical test cases of increasing complexity. In the last part of the paper, \red{two engineering problems dealing with the aerodynamic design of an aircraft wing and with the multidisciplinary design of a two-stage sounding rocket are carried out to assess the performance of the proposed approach on representative industrial complexity test cases.}

\section{From classical optimization to Quality-Diversity (QD) optimization}
\label{sec:QD}

Classical continuous optimization problems are often formulated as:

\begin{eqnarray}
\min_{\mathbf{x}} & & f(\mathbf{x}) \\ 
\text{s.t.} & &  g_i(\mathbf{x}) \leq 0 \;\;\; \mbox{for} \: i=1,\ldots,n_g \\ \label{probRBDO}
& & h_j(\mathbf{x}) = 0 \;\;\; \mbox{for} \: j=1,\ldots,n_h \\
& & \mathbf{x}_\text{lb} \leq \mathbf{x} \leq \mathbf{x}_\text{ub}
\end{eqnarray}
where $\mathbf{x}\in[\mathbf{x}_\text{lb},\mathbf{x}_\text{ub}] \subset{\mathbb{R}^d}$ is a vector of continuous variables (with $\mathbf{x}_\text{lb}$ and $\mathbf{x}_\text{ub}$ the vectors of lower bounds and upper bounds), $f(\cdot)$ is a scalar objective function, $g_i(\cdot)$ is the $i^\text{th}$ inequality constraint function for $i\in\{1,\ldots,n_g\}$ and $h_j(\cdot)$ is the $j^\text{th}$ equality constraint function for $j\in\{1,\ldots,n_h\}$. The number of inequality constraints is equal to $n_g$ and the number of equality constraints is equal to $n_h$.

Moreover, in various applications such as complex engineering design problems, in addition to continuous variables, it is necessary to consider the presence of discrete and categorical variables. Categorical variables are qualitative variables that can be unordered (also known as nominal variables, \textit{e.g.}, type of propulsion, type of material) or ordered (also known as ordinal variables, \textit{e.g.}, small, medium, large). The notion of distance is not properly defined between categories (also called variable levels) and although the ordinal variables present an order, there is no distance between the different categories. The discrete variables are quantitative variables taking specific values with a notion of order and metric to estimate a distance between the different possible variable values (\textit{e.g.}, number of engines on an aircraft).\\

The optimization problem with these different variables may be formulated as:
\begin{eqnarray}
\min_{\mathbf{x}^c,\mathbf{x}^d,\mathbf{x}^q} & & f(\mathbf{x}^c,\mathbf{x}^d,\mathbf{x}^q) \\ 
\text{s.t.} & &  g_i(\mathbf{x}^c,\mathbf{x}^d,\mathbf{x}^q) \leq 0 \;\;\; \mbox{for} \: i=1,\ldots,n_g \\ \label{probRBDO}
& & h_j(\mathbf{x}^c,\mathbf{x}^d,\mathbf{x}^q) = 0 \;\;\; \mbox{for} \: j=1,\ldots,n_h \\
& & \mathbf{x}^c_\text{lb} \leq \mathbf{x}^c \leq \mathbf{x}^c_\text{ub} \\
& & \mathbf{x}^d \in \mathcal{X}^d \\
& & \mathbf{x}^q \in \mathcal{X}^q
\end{eqnarray}
with $\mathbf{x}^c,\mathbf{x}^d,\mathbf{x}^q$ respectively the continuous variables, the discrete variables and the categorical variables. $\mathcal{X}^d$ and $\mathcal{X}^q$ correspond to the definition domains for the discrete and categorical variables. The sizes of the different continuous, discrete and categorical search spaces are noted $n_{\mathbf{x}^c}$, $n_{\mathbf{x}^d}$ and $n_{\mathbf{x}^q}$.\\

Different families of optimization algorithms have been proposed to solve optimization problems: the gradient-based algorithms \cite{pedregal2004introduction}, the grid-search algorithms \cite{audet2006mesh,abramson2009mesh}, the branch-and-bound approaches \cite{clausen1999branch}, the evolutionary algorithms \cite{simon2013evolutionary}, the surrogate-based algorithms \cite{forrester2009recent}, \textit{etc}. 
In case only continuous design variables are involved in the optimization problem, gradient-based algorithms (\textit{e.g.}, Broyden – Fletcher – Goldfarb – Shanno (BFGS) \cite{fletcher2000practical}, Sequential Quadratic Programming (SQP) \cite{nocedal2006quadratic}) exploit the information of the gradient of the objective function and the constraints with respect to the design variables in order to converge to a local minimum (that can be a global minima for convex functions and search space). The gradient is used to determine a descent direction in order to improve the current knowledge about the minimum. \red{In case the optimization problem presents several local minima, different strategies such as multi-start approaches may be used to identify the global minimum. However, multi-start strategies may be time-consuming due to the repetition of optimization problem solving and may provide limited efficiency in the presence of a large number of local minima.} 

\begin{figure}[!h]
\begin{center}
\includegraphics[width=0.85\textwidth]{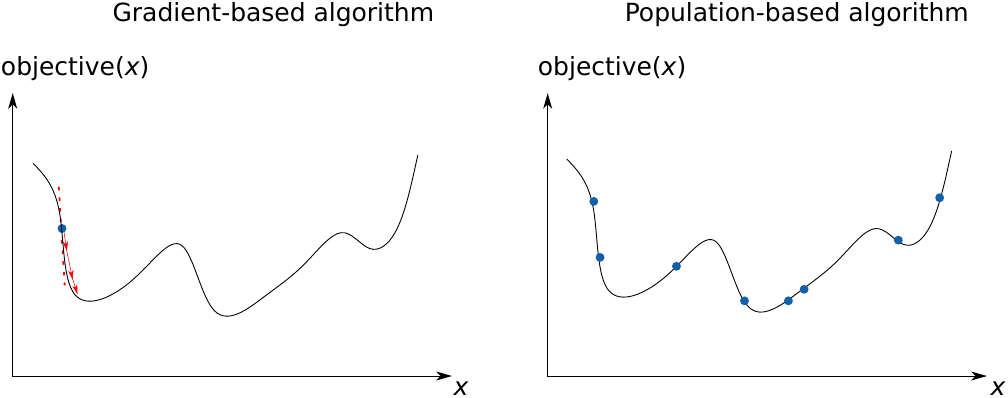}
\caption{General concepts of gradient-based (left) and population-based (right) optimization algorithms}\label{Gradient_vs_evolutionary}    
\end{center}
\end{figure}

For optimization problems with continuous, discrete and categorical variables, as the gradient is not available, when possible, adaptations (such as relaxation approaches \cite{kronqvist2019review}) are required. Alternatively, it is possible to use grid-search algorithms \cite{abramson2009mesh}, branch-and-bounds approaches \cite{clausen1999branch} or population-based algorithms \red{(\textit{e.g.}, genetic algorithm \cite{holland1975adaptation,mitchell1998introduction,sohail2023genetic}, particle swarm \cite{kennedy1995particle,gad2022particle,nayak202325}, covariance matrix adaptation - evolution strategy \cite{hansen1996adapting,hamano2022cma})}. \red{The population-based algorithms (Figure \ref{Gradient_vs_evolutionary}) are part of the family of memetic algorithms \cite{moscato2019accelerated} and rely on a set of individuals that evolve in the search space to identify the optimal regions and the global minimum. Memetic algorithms combine the population-based global search and the heuristic local search done by each of the individuals of the population \cite{moscato2019accelerated,mohammadi2020evolutionary}.} The operators of population-based algorithms differ from an algorithm to another (relying on mutation, cross-over, random generations, evolution of swarm, \textit{etc}.) and are often inspired by natural behaviors. Moreover, most of the existing population-based algorithms have been adapted to handle discrete and categorical variables through specific evolution mechanisms of the population \cite{hamano2022cma,lin2018hybrid,wang2021particle}. In addition, specific approaches have been adapted to deal with constrained optimization problems  \cite{coello2022constraint}. In order to reach convergence, population-based algorithms require in general a large number of objective function and constraint functions evaluations \red{which may be inappropriate in the context of limited evaluation budget.} 
In case the objective function  and / or the constraint functions are computationally intensive, alternative optimization strategies (\textit{e.g.}, Bayesian optimization \cite{garnett2023bayesian,jones1998efficient}) based on surrogate models have been proposed in the literature. In these approaches, each computationally intensive function is replaced by a surrogate model (\textit{e.g.}, Gaussian Process - GP, support vector machine, neural network). Starting from an initial Design of Experiments (DoE), the surrogate models for the objective function and the constraint functions are constructed. Then, an auxiliary optimization problem is solved by optimizing an infill criterion on the surrogate models in order to identify the most promising candidate solution in the search space to find the global minimum. Once the candidate solution is identified, the exact objective and constraint functions are evaluated, the DoE and the surrogate models are updated. This process continues until optimization convergence. This allows to find a global minimum while limiting the number of evaluations of the computationally intensive functions. Gaussian processes have been extensively used for that purpose \cite{le2021revisiting,jones1998efficient}.\\

All the optimization algorithms mentioned above often allow to identify a single optimal solution. \red{Even if some algorithms use a population of individuals that evolve in the search space and might explore different regions, often only a single optimal solution is retained. However, it might be possible to use the dominance concepts (local-global) and to analyze the convergence history of the population to identify different regions with local minima or a Pareto frontier.} 
A new family of optimization approaches, called Quality-Diversity (QD) algorithms \cite{chatzilygeroudis2021quality,pugh2016quality} has been proposed in order to provide a diverse set of optimal solutions characterized by various trade-offs. These techniques offer a different diversity compared to multi-objective optimization algorithms. \red{Indeed, in multi-objective optimization \cite{gunantara2018review,gonzalez2023designing}, in the presence of antagonistic objective functions, a Pareto set is obtained, corresponding to the set of the non-dominated objective solutions in the sense of Pareto dominance \cite{emmerich2018tutorial}, resulting in a trade-off between the different objective functions. This is not the purpose of QD algorithms that do not consider antagonistic objective functions.} 
Indeed, in practice, especially in the early design phases of engineering systems, the design process involves optimization algorithms in order to explore a large design space. The aim is to identify with various high quality options and examine possible trade-offs. Therefore, in early design phases, there is an interest to obtain a set of attractive solutions that can be further explored in more details in the next steps of the design process. QD approaches are based on such an idea in order to provide high-quality solutions with respect to an objective function and diversified with respect to some feature functions (that could be not antagonistic). These features represent different characteristics and preferences for the decision-markers that are not fixed in the current design phase.  For instance, in the early design phases, for the design of lifting surfaces for aerospace vehicles, it can be interesting to identify high-quality wing geometries (in terms for instance of lift-to-drag ratio) while generating diversified solutions according to features such as wing aspect ratio (quantifying how long and slender a wing is) or taper ratio (ratio between tip and root chord lengths). Then, the general idea is to identify a set of optimal solutions (with respect to lift-to-drag ratio) but with a diversity according to features (\textit{e.g.}, aspect ratio, taper ratio). The QD approaches provide a set of high-quality solutions in which the decision-makers can pick in order to further investigate depending on their preferences in terms of features. They generate also valuable information of the influence of the features on the objective function giving some aftermaths of the feature choices on the overall performance. Consequently, for decision-making, QD algorithms allow, with respect to classical optimization algorithms, to go in-depth in the analysis by providing additional information with respect to several quantities of interest modeled using feature functions. \newredbis{QD algorithms have been applied for the design of different industrial devices such as soft grippers \cite{xie2023fin}, robots \cite{mouret2015illuminating}, software product lines \cite{xiang2023automated} or shell structures \cite{sfikas2023design}.} \\

A mixed continuous-discrete-categorical QD problem may be written as:
\begin{eqnarray}\label{Eq_DQ1}
\forall \tilde{\mathbf{f}} \in \mathcal{F}_t, \;\;\;  \min_{\mathbf{x}^c,\mathbf{x}^d,\mathbf{x}^q} & & f(\mathbf{x}^c,\mathbf{x}^d,\mathbf{x}^q) \\ 
\text{s.t.} & &  g_i(\mathbf{x}^c,\mathbf{x}^d,\mathbf{x}^q) \leq 0 \;\;\; \mbox{for} \: i=1,\ldots,n_g \\ \label{probRBDO}
& & h_j(\mathbf{x}^c,\mathbf{x}^d,\mathbf{x}^q) = 0 \;\;\; \mbox{for} \: j=1,\ldots,n_h \\
& & \mathbf{f}_t(\mathbf{x}^c,\mathbf{x}^d,\mathbf{x}^q) \in  \tilde{\mathbf{f}} \\ 
& & \mathbf{x}^c_\text{lb} \leq \mathbf{x}^c \leq \mathbf{x}^c_\text{ub} \\
& & \mathbf{x}^d \in \mathcal{X}^d\\
& & \mathbf{x}^q \in \mathcal{X}^q
\label{Eq_DQ2}
\end{eqnarray}
where $\mathcal{F}_t$ is the feature space of dimension $n$, $\tilde{\mathbf{f}}$ is an element in the feature space and $\mathbf{f}_t(\cdot)$ is a vector of feature functions of size $n$. The feature functions map the design variables (continuous, discrete and categorical variables) into the feature space (a set of dimension $n$ where each coordinate corresponds to a feature of interest). 
Often in practice, the feature space is discretized into a multi-dimensional hyper-rectangular grid. For each feature index $j = 1, \dots, n$, let $s_j \in \mathbb{N}$ be the size of the discretization, representing the number of discretization nodes for the $j^\text{th}$ feature. Let $\left\{ \mathcal{F}_{t_j}^{s_j} \right\}_{j=1,\dots,n}$ be the collection of feature discretizations. The feature space $\mathcal{F}_t = \mathcal{F}_{t_1}^{s_1} \otimes \cdots \otimes \mathcal{F}_{t_n}^{s_n}$ corresponds to the tensor product of each feature coordinate $j$. 
A bin of the multi-dimensional grid corresponds to a niche (a combination of intervals, one per feature coordinate). Therefore, in this context, the elements  $\tilde{\mathbf{f}}$ of the feature space are bins (niches) of the multi-dimensional grid. 
The goal of QD algorithms is to identify the most diverse collection, which is optimal in terms of objective function with respect to the design variables and diverse with respect to the feature functions (represented by the niches). Each member of the collection is as best as possible regarding the objective function. Therefore, the QD algorithm returns a map (also called an archive or a collection) corresponding to a set of solutions that differ in terms of feature characteristics. Each niche in this map contains the solution that will be found by classical optimization algorithms associated to a specific combination of feature values (\textit{e.g.}, \newred{Cartesian} product of intervals) defining the bin.

\begin{figure}[!h]
\begin{center}
\includegraphics[width=0.95\textwidth]{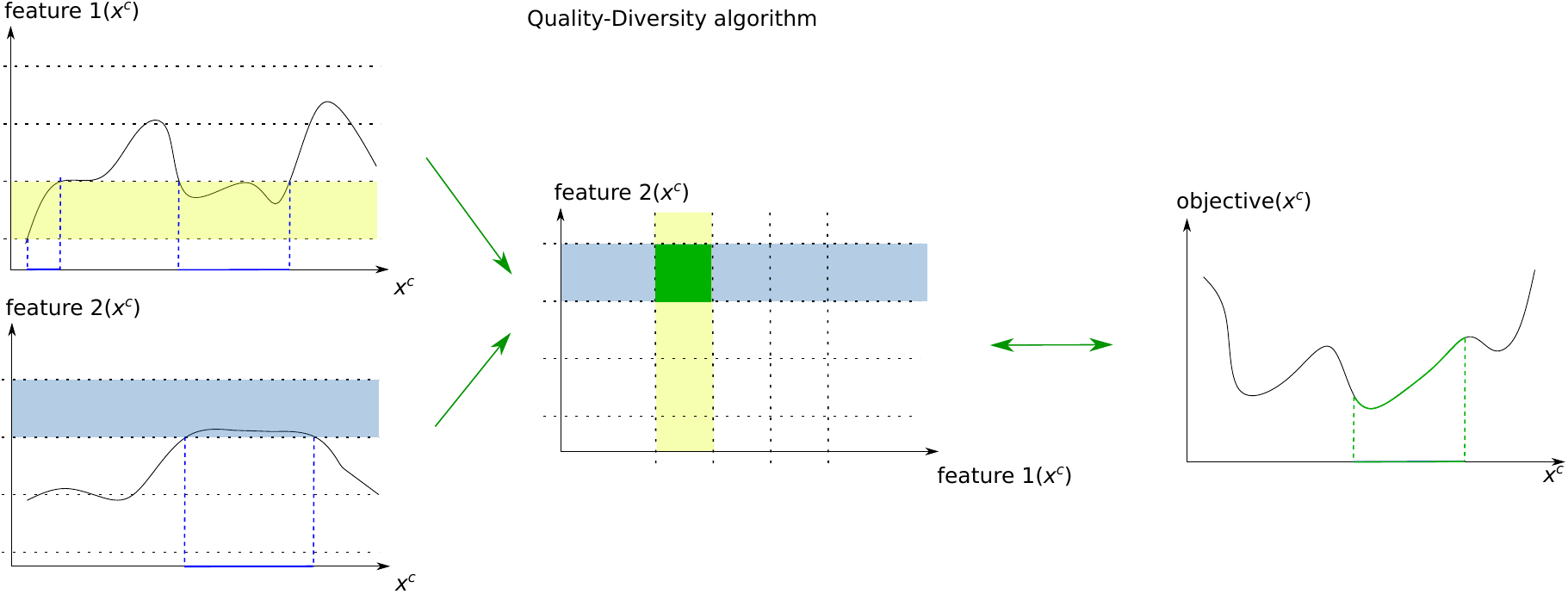}
\caption{Quality-Diversity principle for a single continuous dimensional illustration ($x^c\in\mathbb{R}$) with two feature functions. The green niche (middle of the figure) is determined by the association of the two feature functions discretization (the yellow region of the feature 1 and the blue region of the feature 2, left of the figure). This niche defines a region of the design space (that can be a union of disjoint regions) in which the minimum value of the objective function has to be found (right of the figure). }\label{QD_principle}    
\end{center}
\end{figure}

Figure \ref{QD_principle} illustrates the main characteristics of QD algorithms on a one dimensional continuous problem considering two feature functions. On the left of the figure, the mappings between the continuous design variable $x^c\in \mathbb{R}$ and the two feature functions are represented. A discretization of the two feature functions is made and two particular categories are highlighted (in yellow for the feature 1 and in blue for the feature 2). For each feature function, the feature discretization is associated with some specific regions of the input design space highlighted in blue.  
The two-dimensional grid corresponding to the tensor product of the feature function discretizations is presented in the grid in the middle of the figure. Combining the two discretizations of the two features creates different bins. The combination of the highlighted categories for each feature leads to a particular bin (outlined in green) in the two-dimensional grid. Therefore, as presented in the formulation of the QD-problem (Eqs.(\ref{Eq_DQ1}-\ref{Eq_DQ2})), the goal is to identify the minimum of the objective function (right) associated to each bin of the multi-dimensional grid. Due to the mapping with the features, for each bin, the minimum of the objective function is conditioned to a subset region of the original search space restricted to the combination of the feature associated to the bin. Therefore, the diversity is provided thanks to the combination of features and the quality is ensured thanks to the minimum value of the objective function for each bin. One approach to solve such a problem would be to repeat the solving of the optimization problem for each bin. However, such an approach is not realistic as the number of feature functions and the number of discretization of each feature increases, and the mapping between the design variables and the feature space is not known in advance.\\

Several QD algorithms have been proposed to solve such a type of problems. Firstly, algorithms derived from the population-based approaches have been developed such as: Novelty Search with Local Competition (NSLC) \cite{lehman2011evolving}, Multi-dimensional Archive of Phenotypic Elites (MAP-Elites) \cite{mouret2015illuminating} and various derived versions (\textit{e.g.}, MAP-Elites + Novelty \cite{pugh2016quality}, MAP-Elites + Passive Genetic Diversity \cite{pugh2016quality}, Covariance Matrix Adaptation MAP-Elites (CMA-ME) \cite{fontaine2020covariance}). 
More recently, QD algorithms in the family of surrogate-based approaches have been proposed. Some of the techniques are based on multi-layer perceptron such as Deep surrogate assisted MAP-Elites \cite{zhang2022deep}. However, these methods suffer from the difficulty to account for the uncertainty introduced by the use of the surrogate model in the QD process. Therefore, other approaches based on Gaussian processes have been proposed such as Surrogate-Assisted ILlumination (SAIL) \cite{gaier2017aerodynamic}, Surrogate-assisted PHEnotypic Niching (SPHEN) \cite{hagg2020designing} and Bayesian Optimization of Elites (BOP-Elites) \cite{kent2020bop,kent2023bop} to account for the uncertainty model estimation provided by the Gaussian processes.   

The existing Bayesian QD algorithms \cite{gaier2017aerodynamic,hagg2020designing,kent2020bop,kent2023bop} only handle optimization problems with continuous variables and without constraints. In the following of this paper, a new Bayesian QD approach is proposed in order to solve QD problems with mixed continuous, discrete and categorical variables and to deal with constrained optimization problems. It is assumed that the objective function, the constraints and the features are computationally intensive and replaced by surrogate models that are enriched during the QD process. The main characteristics of the proposed algorithm are introduced in the following section. 

\section{Bayesian Quality-Diversity approach for constrained optimization problems with mixed continuous, discrete and categorical variables}
\label{sec:QD_BO}

In this section, the proposed Bayesian QD algorithm is described. It allows to solve constrained QD problems with mixed continuous, discrete and categorical variables. Firstly, Gaussian process and its adaptation in order to deal with mixed variables are described in Section \ref{GP_mixed}. Then, in Section \ref{BQD_section}, the proposed Bayesian QD algorithm is presented. A focus on constrained Bayesian optimization is made along with the adaptation of mechanisms to handle the constraints in QD problems.

\subsection{Gaussian process with mixed variables}\label{GP_mixed}

Gaussian Process (GP) is a surrogate model (sometimes referred to as Kriging \cite{oliver1990kriging,santner2003design} for conditioned GP) that may be used to replace a generic computationally intensive function $f: \mathbb{R}^d \rightarrow \mathbb{R}$ with an approximation. A GP corresponds to a collection of random variables, any finite number of which has a multivariate joint Gaussian distribution. It may be seen as a generalization of the Gaussian probability distribution by encoding a distribution over a set of functions \cite{rasmussen2003gaussian}. 

\subsubsection{Gaussian process with continuous variables}
In this section, only continuous variables are considered. A GP is fully determined by its mean function $\mu(\cdot)$ and its covariance function $\text{Cov}(\cdot,\cdot)$. If the function of interest $f(\cdot)$ follows a GP, it can be expressed as  $f(\cdot) \sim \text{GP}(\mu(\cdot),\text{Cov}(\cdot,\cdot))$. 
To build a GP, it is necessary to solve a supervised learning problem.
A GP is trained on a Design of Experiments constituted of an input training set of size $M$, $\mathcal{X}_M=\left\{ \mathbf{x}_{1}^c, \ldots , \mathbf{x}_{M}^c \right\}$ and the corresponding set of computationally expensive function responses  $\mathcal{Y}_M= \left\{y_{1}=f\left(\mathbf{x}_{1}^c\right),\ldots, \right.$ $\left. y_{M}=f\left(\mathbf{x}_{M}^c\right) \right\}$. The function responses form a vector noted $\mathbf{y}_M=\left[y_1,\cdots,y_M \right]^T$. In a regression context, a GP prior is assumed on the mean function and on the covariance function. 
Regarding the mean function, as the tendency of the exact function is often unknown, a constant function $\mu$ is generally assumed as GP prior, resulting in ordinary Kriging. Depending on available \textit{a-priori} knowledge, other types of mean function may be assumed such as quadratic or more general basis functions. 
The covariance function $\text{Cov}(\cdot,\cdot)$ is usually defined through the use of a parameterized covariance function called a \textit{kernel} $k^{\boldsymbol{\Theta}}(\cdot,\cdot)$. The kernel is used to model the covariance between two elements (\textit{e.g.}, $\mathbf{x}^c $ and $ \mathbf{x}^{c'}$) as a symmetric positive definite function of the values of the coordinates of the elements $k^{\boldsymbol{\Theta}}( \mathbf{x}^c , \mathbf{x}^{c'})$. This covariance function depends on some hyper-parameters $\boldsymbol{\Theta}$ that need to be determined during the training process with the DoE. The most known kernels \cite{rasmussen2003gaussian} are the Squared Exponential kernel (also known as Radial Basis Function), the Rational Quadratic kernel, the Matérn kernel, \textit{etc.} The covariance function is a key element in GP. The covariance function allows to encode some assumptions on the behavior of the exact function (\textit{e.g.}, smoothness, periodicity, stationarity, separability). Multidimensional kernels may be obtained by combining single dimensional kernels through for instance a product operator following the formalism of Reproducing Kernel Hilbert Space (RKHS) \cite{alvarez2012kernels}. 
The prior mean and prior covariance are updated by relying on the information on the modeled function through the data set $\{\mathcal{X}_M,\mathcal{Y}_M \}$, which enables to provide a more insightful model of the considered function.

If a constant mean function is assumed, the GP is defined such that $f(\cdot) \sim \mathcal{N}\left(\mu,k^{\boldsymbol{\Theta}}\left(\cdot,\cdot \right) \right)$ with $\mathcal{N}(\cdot,\cdot)$ the Gaussian distribution. Considering the DoE $\{\mathcal{X}_M,\mathcal{Y}_M \}$, the GP has a multivariate Gaussian distribution with a covariance matrix $\mathbf{K}_{MM}$ built from the parameterized covariance function $k^\Theta(\cdot,\cdot)$ on the input dataset $\mathcal{X}^M$ (the dependence on $\boldsymbol{\Theta}$ is dropped to simplify the notations).
In the presence of experimental or numerical noisy data, the relationship between the latent function values $\mathbf{f}_M = f\left(\mathcal{X}_M\right)$  and the observed responses $\mathcal{Y}_M$ is given by: $ p\left(\mathbf{y}| \mathbf{f}_M \right) = \mathcal{N}\left(\mathbf{y}|\mathbf{f}_M, \sigma^2 \mathbf{I} \right)$ with $\sigma^2$ an assumed Gaussian noise variance.\\

Then, from these input and output training sets and the prior on the GP, it is possible to train the surrogate model using the marginal likelihood. It is obtained by integrating out the latent function giving $
p \left(\mathbf{y}_M|\mathcal{X}_M, \boldsymbol{\Theta},\mu,\sigma \right) = \mathcal{N}\left(\mathbf{y}_M| \mu, \mathbf{K}_{MM} + \right.$ $\left. \sigma^2 \mathbf{I}_{MM} \right)$ with $\mathbf{I}_{MM}$ the identity matrix of size $M$. To simplify the notations, let define $\mathbf{\hat{K}}_{MM} = \mathbf{K}_{MM}+\sigma^2 \mathbf{I}$. In practice, the GP training requires to minimize the negative Log-Marginal Likelihood (LML) with respect to the hyperparameters $\boldsymbol{\Theta}$, $\mu$ and  $\sigma$. The LML $L(\cdot)$ is given by:
\begin{eqnarray}
L \left(\boldsymbol{\Theta}, \mu, \sigma|\mathcal{X}_M,\mathcal{Y}_M \right) &=& \log \left(p \left(\mathbf{y}_M|\mathcal{X}_M, \boldsymbol{\Theta}, \mu, \sigma \right)\right) \\
&\propto &  \log \left(|\mathbf{\hat{K}}_{MM}| \right) - \mathbf{y}_M^T\mathbf{\hat{K}}^{-1}_{MM} \mathbf{y}_M 
\end{eqnarray}
where all the kernel matrices implicitly depend on the hyperparameters $\boldsymbol{\Theta}$. To solve the optimization problem, any optimizer may be used (\textit{e.g.}, gradient-based, population-based algorithms). Moreover, a closed form of the constant mean function may be sometimes determined \cite{rasmussen2003gaussian}. 

%\begin{figure}[h!]
%\begin{center}
%\includegraphics[width=0.7\linewidth]{GP_principle.pdf}
%\caption{Example of Gaussian process prediction and associated confidence interval}
%\label{Kriging_example}
%\end{center}
%\end{figure}

Once the GP has been trained (the optimal hyperparameters have been determined, noted $\boldsymbol{\hat{\Theta}}$, $\hat{\mu}$ and  $\hat{\sigma}$), the prediction $y^*$ at a new unknown location $\mathbf{x}^{c*}\in \mathbb{R}^d$ is done by using the conditional properties of a multivariate normal distribution:
\begin{equation}
p\left( y^*|\mathbf{x}^{c*},\mathcal{X}_M,\mathcal{Y}_M,\boldsymbol{\hat{\Theta}}, \hat{\mu},\hat{\sigma} \right) = \mathcal{N}\left(y^* | \hat{f}^*,\hat{s}^{*2} \right)
\end{equation}
where $\hat{f}^*,\hat{s}^{*2}$ are the mean prediction and the associated variance. These terms are defined by:
\begin{eqnarray}
\hat{f}^* &=& \hat{f}(\mathbf{x}^{c*}) = \hat{\mu} + \mathbf{k}^T_{\mathbf{x}^*}\left({\mathbf{K}}_{MM} + \hat{\sigma}^2 \mathbf{I} \right)^{-1} \left(\mathbf{y}_M - \mathbf{1}\hat{\mu} \right) \\
\hat{s}^{*2}&=& \hat{s}^{2}(\mathbf{x}^{c*}) =  k_{\mathbf{x}^{c*},\mathbf{x}^{c*}} - \mathbf{k}^T_{\mathbf{x}^{c*}}\left({\mathbf{K}_{MM}} + \hat{\sigma}^2 \mathbf{I} \right)^{-1}\mathbf{k}_{\mathbf{x}^{c*}}
\end{eqnarray}
\noindent where $k_{\mathbf{x}^{c*},\mathbf{x}^{c*}} = k(\mathbf{x}^{c*},\mathbf{x}^{c*})$ and $\mathbf{k}_{\mathbf{x}^{c*}} = \left[k \left( \mathbf{x}_{i}^{c}, \mathbf{x}^{c*} \right)\right]_{i=1,\ldots,M}$.

\noindent $\hat{f}(\mathbf{x}^{c*})$ and $\hat{s}(\mathbf{x}^{c*})$ correspond to the mean and variance of the posterior of the GP. A key element about GP is the possibility to have access to a prediction and a confidence level associated to the prediction that may be used for instance in an active learning strategy with a refinement process \cite{garnett2023bayesian}.  

\subsubsection{GP with mixed continuous, discrete and categorical variables}

In order to adapt GP to the presence of mixed continuous, discrete and categorical variables, it is necessary to define an adapted covariance function that can deal with such a type of variables. 

%In this section, to ease the notation and without loss of generality, two types of variables are considered: the continuous variables $\mathbf{x}=[x_{(1)},\cdots,x_{(n_x)}]^T\in \mathbb{R}^{n_x}$ and the categorical variables $\mathbf{z}=[z_{(1)},\cdots,z_{(n_z)}]^T\in \mathcal{F}_z$ of dimension $n_z$. In the following, the discrete variables are included in the family of categorical variables. Categorical are characterized as non-relaxable variables defined within a finite set of choices, or levels, and characterized by an absence of relation of order between the possible choices. A discrete variable may be seen as a special case of categorical variables for which a relation of order exists between the variable levels, categorical variables are therefore more general and used in the following. 

The kernel to deal with mixed continuous, discrete and categorical variables may be defined as a product of one-dimensional kernels following the RKHS formalism. The resulting mixed-variable kernel can then be defined as \cite{pelamatti2021mixed} for a couple of continuous ($\mathbf{x}^c,\mathbf{x}^{c'}$), discrete ($\mathbf{x}^d,\mathbf{x}^{d'}$) and categorical variables ($\mathbf{x}^q,\mathbf{x}^{q'}$):
\begin{equation}
\label{KernelProd2}
\begin{split}
k(\{ \mathbf{x}^c, \mathbf{x}^d,\mathbf{x}^q\},\{ \mathbf{x}^{c'}, \mathbf{x}^{d'},\mathbf{x}^{q'} \})  = \prod_{i = 1}^{n_{x^c}}  k_{x^c_{(i)}}\left(x^c_{(i)},x^{c'}_{(i)}\right) \times \prod_{j = 1}^{n_{x^d}}  k_{x^d_{(j)}}\left(x^d_{(j)},x^{d'}_{(j)}\right) \times\\
\prod_{k = 1}^{n_{x^q}}  k_{x^q_{(k)}}\left(x^q_{(k)},x^{q'}_{(k)}\right)
\end{split}
\end{equation}
with $k_{x^c_{(i)}}(\cdot,\cdot)$ the kernel associated to the $i^\text{th}$ coordinate of the continuous variable vector $\mathbf{x}^c$, and $k_{x^d_{(j)}}(\cdot,\cdot)$ and  $k_{x^q_{(k)}}(\cdot,\cdot)$ the kernels associated to the $j^\text{th}$ coordinate discrete variable vector and the $k^\text{th}$ coordinate categorical variable vector.  For the kernels associated to the continuous variables, any kernel discussed in the previous section may be used. 

Different kernels have been proposed in the literature to deal with discrete and categorical variables \cite{saves2022general,pelamatti2021mixed}. Due to the absence of classical distance measures between the values of these variables (especially categorical variables), adapted kernels have to be used. In the following, the kernels are discussed for single dimensional case and the extension to multidimensional problem is done using Eq.(\ref{KernelProd2}).  Most often, no distinction is made between discrete and categorical variables and generic kernel for these variables will be noted $k_z$. In the following of this section, discrete or categorical variables $x^d$ or $x^q$ are noted using the generic notation $z$.

Considering a scalar categorical (or discrete) variable $z$ with $L$ possible levels $\{z^1,\dots, z^L\}$. 
To ease the notation, as these variables are characterized by a finite number of levels, the kernel function returns a finite number of covariance values that can be organized into a $L \times L$ matrix $\mathbf{K}$ defined such that:
\begin{equation}
    \mathbf{K}_{m,d} = k_z(z=z^m,z'=z^d)
\end{equation}
with $(m,d)\in \{1,\dots, L\}^2$.
In order to be a valid covariance matrix, the covariance kernel $\mathbf{K}$ has to be symmetric and positive semi-definite. In the following, two different kernels \cite{saves2022general,pelamatti2021mixed} adapted for discrete and categorical variables are presented. 

\paragraph{Compound Symmetry kernel}

The Compound Symmetry (CS) kernel \cite{saves2022general,halstrup2016black,pelamatti2020overview} is also known as the Gower kernel as the distance between the possible categorical levels for a variable $z$, is expressed based on the Gower distance \cite{gower1971general}.

CS is characterized by a single covariance value for any-non identical pair of inputs. Relying on a continuous squared exponential kernel, CS kernel is defined as (for two discrete scalar variable $z$ and $z'$ with $L$ possible levels):

\begin{equation}
k_{z}(z,z') = \sigma_{z}^2 \exp \left( - \theta d_\text{gow}(z,z') \right)
\end{equation}
with
\begin{equation}
d_\text{gow}(z,z') =  \left\{
    \begin{array}{ll}
        0 & \mbox{if } z=z' \\
        1 & \mbox{if } z\neq z'
    \end{array}
\right.\end{equation}
$\sigma_{z}^2$ and $\theta \geq 0$ are respectively the variance and hyperparameter associated to the CS kernel. \red{In case the discrete or categorical variables $z$ and $z'$ have the same level then, the distance is null, otherwise independently of the values of the levels for the two variables, the kernel returns the same value. In the present derivation of Compound Symmetry, the Gower distance is always a distance as $z$ and $z'$ can take the same possible discrete levels, therefore there is no possible missing value.}

The CS kernel is very simple to account for the effect of a given  variable since it relies on a single hyperparameter. However, the covariance between any pair of non identical levels of a given  variable $z$ is the same, regardless of the level values. This assumption may be too simplistic, especially when dealing with discrete and categorical variables which present a large number of levels. To avoid such a phenomenon, alternative kernels may be considered as discussed in the following paragraph.

\paragraph{Hypersphere decomposition kernel}

An alternative kernel to deal with discrete or categorical variables is based on the hypersphere decomposition \cite{saves2022general,pelamatti2020overview,zhou2011simple}. The general idea is to use a mapping between each level of the variable $z$ and a point on the surface of a $L$-dimensional hypersphere. This decomposition allows to ensure the symmetric and positive semi-definite nature of the corresponding covariance matrix \cite{zhou2011simple,rebonato2011most}. In practice, a mapping $F$ is defined based on a polyspherical change of coordinates $F: [-\pi,\pi]^{L-1} \rightarrow \mathbb{S}^{L-1}$ with $\mathbb{S}^{L-1}$ the unit $(L-1)$-sphere such that: $\mathbb{S}^{L-1}=\{\mathbf{x}\in \mathbb{R}^{L} : \parallel \mathbf{x} \parallel = 1\}$. This corresponds to a change of coordinates between the spherical coordinates and the Cartesian coordinates. The matrix $\mathbf{K}$ associated to the kernel of the hypersphere decomposition is given by $\mathbf{K} = \sigma_z^2 \mathbf{L}^T\mathbf{L}$ with $\mathbf{L}$ the lower triangular matrix associated to the Cholesky decomposition of $\mathbf{K}$. $\mathbf{L}$ is defined such that:

\scriptsize
\begin{eqnarray}
   \mathbf{L} = \left[
    \begin{array}{lllll}
        1 & 0 & \cdots & \cdots & 0 \\
        \cos(\theta_{2,1}) & \sin(\theta_{2,1}) & 0 & \cdots & 0 \\
        \vdots & \vdots & \vdots & \vdots & \vdots \\
        \cos(\theta_{L,1}) & sin(\theta_{L,1})\cos(\theta_{L,2}) & \cdots & \cos(\theta_{L,L-1})\prod_{k=1}^{L-2}\sin(\theta_{L,k}) & \prod_{k=1}^{L-1}\sin(\theta_{L,k})
    \end{array}
\right]
\end{eqnarray}
\normalsize
The hypersphere kernel is parametrized by the $\theta_{i,j}$ hyperparameters involved in the matrix $\mathbf{L}$, corresponding to $L\times(L-1)/2$ parameters. Compared to the CS kernel, the hypersphere decomposition kernel is able to give a different covariance value for each pair of levels characterizing the variable $z$. Moreover, this covariance function can take negative values, as each covariance value is computed as the product between a number of sine and cosine functions. Therefore, it is possible to model positive and negative correlations between the  levels. 
It is important to notice that a part of the hyperparameters characterizing this kernel influences several covariance values simultaneously (due to the repeated presence of the same $\theta_{ij}$ for different covariance values) that can lead to difficulties to estimate the optimal values of the hyperparameters.

In \cite{saves2022general}, the authors proposed a formulation that unifies the definition of the CS and hypersphere decomposition kernels in the context of Gaussian kernel. This formulation is an extension of the continuous Gaussian kernel to deal with mixed continuous-categorical variables. The Gaussian hypersphere decomposition kernel consists of the composition of a Gaussian kernel with an hypersphere decomposition kernel: $\mathbf{K} = \sigma_z^2 \exp(-\mathbf{L}^T\mathbf{L})$. This composition allows to present the CS kernel and hypersphere decomposition kernel under the same formalism. However, the Gaussian hypersphere decomposition kernel is only able to represent positive correlations between two levels. An equivalence between hypersphere and Gaussian hypersphere if the correlations are positive has been presented in \cite{saves2022general} (with a restriction for the hypersphere angles to be in $[0,\pi/2]$).

There exist alternative kernels to deal with discrete and categorical variables such as the latent variable kernel \cite{zhang2020latent} or the coregionalization \cite{pelamatti2021mixed,alvarez2012kernels}. For more information, please refer to \cite{pelamatti2021mixed,pelamatti2020overview}.

%\subsection{Bayesian optimization with constraints}

\subsection{Bayesian Quality-Diversity algorithm}\label{BQD_section}
\subsubsection{Overview of the proposed algorithm}

The proposed Bayesian QD algorithm is derived from SPHEN \cite{hagg2020designing} and SAIL \cite{gaier2017aerodynamic} approaches to extend them in order to account for mixed continuous, discrete and categorical variables and to handle constrained optimization problems. 

\begin{figure}[!h]
\begin{center}
\includegraphics[width=0.85\textwidth]{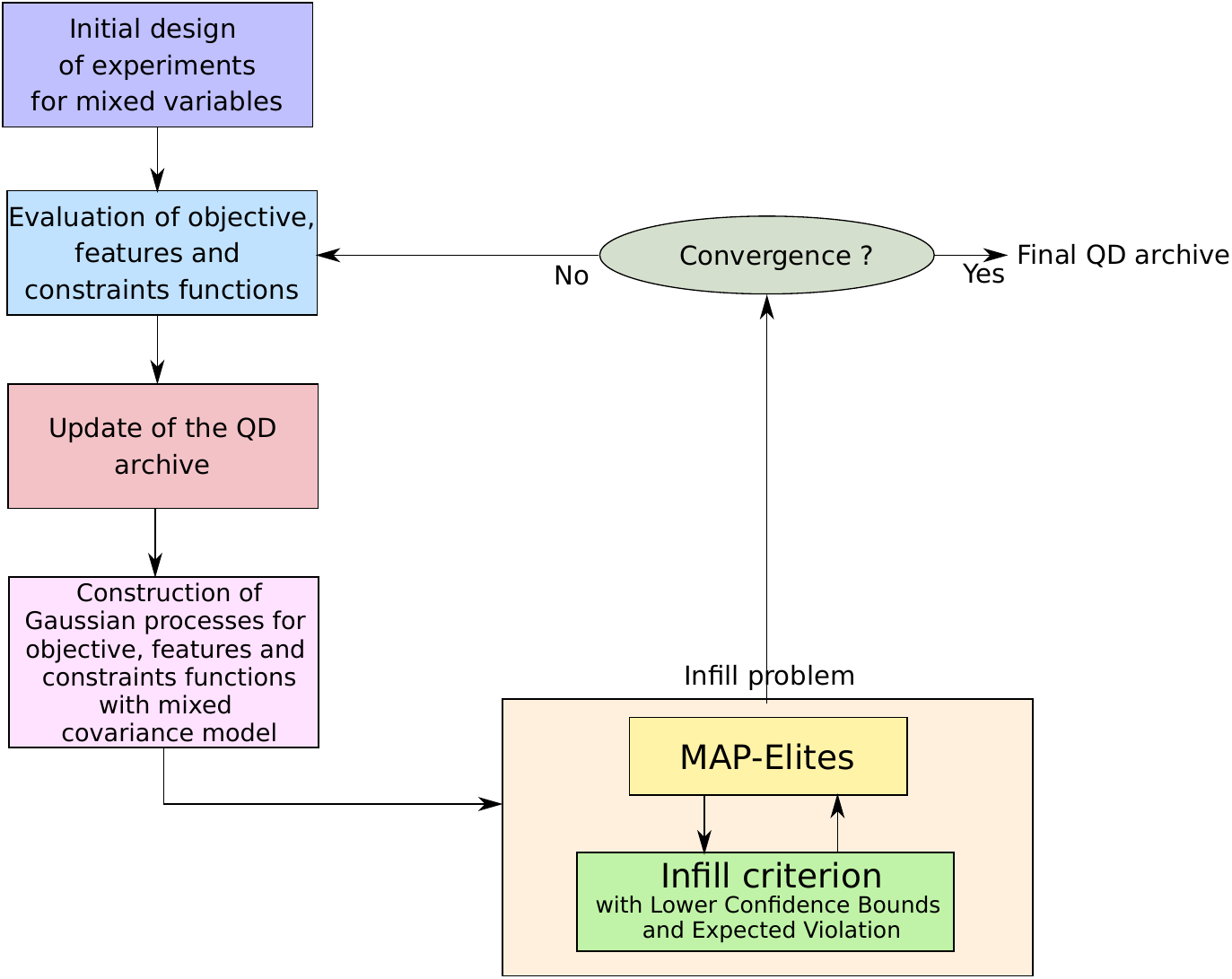}
\caption{Proposed Bayesian QD algorithm to deal with mixed continuous, discrete and categorical variables and the presence of constraints}\label{QD_algo_blocks}    
\end{center}
\end{figure}

In the proposed approach, the objective function, the feature functions and the constraint functions are replaced by Gaussian processes with adapted covariance models (either based on CS kernel or hypersphere decomposition kernel depending on the regularity of the functions to be modeled). 
The proposed approach (Figure \ref{QD_algo_blocks}, Algorithm \ref{alg:QS}) starts with a Design of Experiments of size $M$ in the joint space of continuous, discrete and categorical spaces. The continuous variables $\mathbf{x}^c$ may be sampled according to classical DoE strategies such as Latin Hypercube Sampling (LHS) \cite{mckay2000comparison}, Sobol' sampling \cite{sobol1967distribution}, random sampling, \textit{etc}. The discrete $\mathbf{x}^d$  and categorical $\mathbf{x}^q$  variables are sampled randomly based on the possible values for each variable. Joint sampling in the continuous, discrete and categorical spaces may be possible with adapted space filling \cite{deng2015design}. Based on the input DoE $\mathcal{X}_M = \{ \{\mathbf{x}_1^c,\mathbf{x}_1^d,\mathbf{x}_1^q\}, \dots, \{\mathbf{x}_M^c,\mathbf{x}_M^d,\mathbf{x}_M^q\} \}$, the objective function, the features and the constraints are evaluated to get the corresponding outputs: for the objective function $\mathcal{Y}_M = \{y_1 = f(\mathbf{x}_1^c,\mathbf{x}_1^d,\mathbf{x}_1^q), \dots, y_M = f(\mathbf{x}_M^c,\mathbf{x}_M^d,\mathbf{x}_M^q)  \}$, for the features $f_{t_j}(\cdot)$,  $\mathcal{F}_{{t_j}_M} = \{f_{t_j}(\mathbf{x}_1^c,\mathbf{x}_1^d,\mathbf{x}_1^q), \dots, f_{t_j}(\mathbf{x}_M^c,\mathbf{x}_M^d,\mathbf{x}_M^q)  \}$ and for the constraints (here inequality) $g_i(\cdot)$, $\mathcal{G}_{{i}_M} = \{g_{i}(\mathbf{x}_1^c,\mathbf{x}_1^d,\mathbf{x}_1^q), \dots,$ $ g_{i}(\mathbf{x}_M^c,\mathbf{x}_M^d,\mathbf{x}_M^q)  \}$.
Based on these DoEs, a GP with an adapted covariance model is constructed for each function to be modeled. Moreover, a QD archive is created from the input DoE and the corresponding evaluations of the exact functions. Based on the DoE, only feasible solutions with respect to the exact constraints are added in the archive. Moreover, for each discovered bin (defined with the exact features), the best solution with respect to the exact objective function is added in the niche.

Using these GPs, an auxiliary optimization problem is solved in order to identify the most promising candidates to be added to the DoEs in order to improve the current archive and the GPs. 

\subsubsection{Infill problem optimization}

The auxiliary optimization problem consists in minimizing an infill criterion under some constraints. The infill criterion has to account for the constraint functions in order to identify interesting solutions with respect to the objective function while being feasible with respect to the constraints. Different infill criteria can be used \cite{jones1998efficient}. In this paper, the selected infill criterion is the Lower Confidence Bound combined with the Expected Violation for the constraints. Consequently, the auxiliary optimization problem is formulated as follows:
\begin{eqnarray}
\label{probaux_QD_1}
\forall \tilde{\mathbf{f}} \in \mathcal{F}_t, \;\;\;  \min_{\mathbf{x}^c,\mathbf{x}^d,\mathbf{x}^q} & & \hat{f}(\mathbf{x}^c,\mathbf{x}^d,\mathbf{x}^q) - k \times \hat{s}(\mathbf{x}^c,\mathbf{x}^d,\mathbf{x}^q) \\ 
\text{s.t.} & &  EV_{\hat{g}_i}(\mathbf{x}^c,\mathbf{x}^d,\mathbf{x}^q) \leq t_i \;\;\; \mbox{for} \: i=1,\ldots,n_g \\ 
& & \mathbf{\hat{f}}_t(\mathbf{x}^c,\mathbf{x}^d,\mathbf{x}^q) \in  \tilde{\mathbf{f}} \\ 
& & \mathbf{x}^c_\text{lb} \leq \mathbf{x}^c \leq \mathbf{x}^c_\text{ub} \\
& & \mathbf{x}^d \in \mathcal{X}^d, \;\; \mathbf{x}^q \in \mathcal{X}^q  
\label{probaux_QD_2}
\end{eqnarray}
with $EV_{\hat{g}_i}(\cdot) = \hat{g}_i(\cdot)\times \Phi\left(\frac{\hat{g}_i(\cdot)}{\hat{s}_{g_i}(\cdot)}\right)+\hat{s}_{g_i}(\cdot)\times \phi \left(\frac{\hat{g}_i(\cdot)}{\hat{s}_{g_i}(\cdot)} \right) $, $\Phi(\cdot)$ and $\phi(\cdot)$ the Cumulative Distribution Function and the Probability Density Function of the standard Normal distribution, $\hat{g}_i(\cdot)$ the prediction of the posterior GP associated to the constraint $g_i(\cdot)$ and $\hat{s}_{g_i}(\cdot)$ the standard deviation associated to the posterior GP. $EV_{\hat{g}_i}(\cdot)$ corresponds to the Expected Violation (EV) \cite{audet2000surrogate} associated to the constraint $g_i(\cdot)$ and accounts for both the GP prediction and its associated uncertainty to evaluate the risk of violation of the exact constraint functions. $k$ is a positive scalar parameter specifying the exploitation / exploration balance in the Bayesian optimization process. $t_i$ is a threshold corresponding to the maximum accepted constraint violation considering uncertainty associated to the GP of the constraint $g_i(\cdot)$. In this optimization problem, the objective function, the features and the constraints are replaced by their respective GPs, reducing the computational cost associated to the infill problem optimization. To solve this auxiliary constrained optimization problem, an adaptation from original MAP-Elites algorithm \cite{mouret2015illuminating} is used to handle continuous, discrete and categorical variables and also to manage the presence of constraints. This algorithm is described in the following section.

\begin{algorithm}[h!]
\caption{Bayesian QD algorithm for constrained optimization problem with mixed continuous, discrete and categorical variables } \label{alg:QS}
\begin{algorithmic}
\State 1) Initialization: $M$ (initial DoE size), $n$ (number of features), $n_g$ (number of constraints)
\State 2) Initial DoE:
\State \hspace{0.5cm}  $\mathcal{X}_M \gets \{ \{\mathbf{x}_1^c,\mathbf{x}_1^d,\mathbf{x}_1^q\}, \dots, \{\mathbf{x}_M^c,\mathbf{x}_M^d,\mathbf{x}_M^q\} \}$
\State 3) Evaluation of the objective function: 
\State \hspace{0.5cm}  $\mathcal{Y}_M \gets \{y_1 = f(\mathbf{x}_1^c,\mathbf{x}_1^d,\mathbf{x}_1^q), \dots, y_M = f(\mathbf{x}_M^c,\mathbf{x}_M^d,\mathbf{x}_M^q)  \}$ %\Comment{Evaluation of the objective function}
\State 4) Evaluation of the feature functions: 
\State \hspace{0.5cm}  $\mathcal{F}_{{t_j}_M} \gets \{f_{t_j}(\mathbf{x}_1^c,\mathbf{x}_1^d,\mathbf{x}_1^q), \dots, f_{t_j}(\mathbf{x}_M^c,\mathbf{x}_M^d,\mathbf{x}_M^q)  \} \;\; \text{for } j\in\{1,\dots,n\}$
\State 5) Evaluation of the constraint functions: 
\State \hspace{0.5cm}  $\mathcal{G}_{{i}_M} = \{g_{i}(\mathbf{x}_1^c,\mathbf{x}_1^d,\mathbf{x}_1^q), \dots,$ $ g_{i}(\mathbf{x}_M^c,\mathbf{x}_M^d,\mathbf{x}_M^q)  \} \;\; \text{for } i\in\{1,\dots,n_g\}$
\State 6) Creation of an exact archive 
\While{(stopping criterion not reached)}
    \State 7) Build GPs for the objective function, the feature functions and the constraint functions
    \State 8) Run MAP-Elites to solve auxiliary optimization problem Eqs.(\ref{probaux_QD_1}-\ref{probaux_QD_2})
        \State \hspace{0.5cm}   8a) Creation of a MAP-Elite archive
        \State \hspace{0.5cm} \textbf{while} (stopping criterion not reached) \textbf{do} \Comment{See Algorithm \ref{alg:MAP_Elites}}
            \State \hspace{1.cm}  8b) Generation of the MAP-Elites population
            \State \hspace{1.cm}  8c) Evaluation of the auxiliary functions 
            \State \hspace{1.cm}  8d) Random selection of elites and mutation operation
            \State \hspace{1.cm}  8e) Evaluation of the children on the auxiliary functions
            \State \hspace{1.cm}  8f) Update of the MAP-Elite archive
            \State \hspace{0.5cm} \textbf{end while} 
            %\EndWhile
    \State 9) Use of Sobol' sequence to select $p$ elites 
    \State \hspace{0.5cm}   $\mathcal{X}_p \gets \{ \{\mathbf{x}_1^{c*},\mathbf{x}_1^{d*},\mathbf{x}_1^{q*}, \dots, \{\mathbf{x}_p^{c*},\mathbf{x}_p^{d*},\mathbf{x}_p^{q*}\} \}$
    \State 10) Evaluations of the exact functions based on $\mathcal{X}_p$
    \State 11) Update of the DoEs 
    \State \hspace{0.5cm}   $\mathcal{X}_M \gets \mathcal{X}_M \cup \mathcal{X}_p$, $\mathcal{Y}_M \gets \mathcal{Y}_M \cup \mathcal{Y}_p$, $\mathcal{F}_{{t_j}_M} \gets \mathcal{F}_{{t_j}_M}\cup \mathcal{F}_{{t_j}_p} \text{ for } j\in\{1,\dots,n\}$ and  $\mathcal{G}_{{i}_M} \gets \mathcal{G}_{{i}_M}\cup \mathcal{G}_{{i}_p} \text{ for } i\in\{1,\dots,n_g\}$
    \State 12) Update the exact archive 
\EndWhile
\end{algorithmic}
\end{algorithm}

\subsubsection{Adaptation of MAP-Elites algorithm for infill problem optimization \label{discreteMAPElites}}
The adapted MAP-Elites algorithm \red{(see Algorithm \ref{alg:MAP_Elites} in \ref{Appendix_MAP_Elites_algo}}) starts with the generation of an initial population of individuals randomly distributed in the search space \red{(Step 2 of Algorithm \ref{alg:MAP_Elites})}. For each discrete and categorical variable, each possible level has an equal probability to be chosen using a multinomial distribution. Based on this initial population, the auxiliary objective function, the auxiliary feature functions and the auxiliary constraints are evaluated using the surrogate models (Eqs.(\ref{probaux_QD_1})-(\ref{probaux_QD_2})) \red{(Step 3 of Algorithm \ref{alg:MAP_Elites})}. An initial archive is defined based on the responses of the functions and their association to the corresponding bins \red{(Step 4 of Algorithm \ref{alg:MAP_Elites})}. 
Then, among the individuals of the population, some individuals are randomly selected, called the elites \red{(Step 5 of Algorithm \ref{alg:MAP_Elites})}. For each elite, each coordinate of the vector $[\mathbf{x}^c,\mathbf{x}^d,\mathbf{x}^q]^T$ has a chance of being mutated (depending on a mutation probability). The continuous coordinates are mutated according to a Gaussian distribution (characterized by a variance parameter) and the discrete and categorical variables are mutated according to a multinomial distribution on the existing levels for the coordinate. For the discrete variables, the  probability mass functions  of the multinomial distribution used for the mutation can depend on the distance between the different levels for the coordinate \red{(Step 6 of Algorithm \ref{alg:MAP_Elites})}. The elites are then evaluated on the auxiliary functions based on GPs (objective, features and constraints, in \red{Step 7 of Algorithm \ref{alg:MAP_Elites}}). The new solutions are added to the current archive based on different rules \red{(Step 8 of Algorithm \ref{alg:MAP_Elites})}. First, only feasible solutions with respect to the auxiliary optimization problem are added in the archive (constraint dominance approach \cite{coello2002constraint}). Then, if a new solution discovers an unoccupied bin, the elite is added to the archive bin. If a new solution is better in terms of objective function than the current solution in the archive bin, the elite replaces the existing solution. This process is repeated for a certain number of generations leading to an increasing number of discovered bins and better solutions in each bin.\\
Once the maximal number of generations has been reached in the MAP-Elites process on the surrogate models, among the elites of the final archive, a certain number of elites are selected to be evaluated on the exact objective, features and constraints. As the archive may have a large number of solutions (for instance due to a large number of bins), a Sobol' sequence  \cite{sobol1967distribution} is used to select a limited number of optimal elites that uniformly covers the feature space.  These selected elites are evaluated on the exact objective, features and constraints and added to the DoE $\mathcal{X}_M$. The different GPs of the objective function, features and constraints are updated and a new iteration of the Bayesian QD algorithm is performed. Moreover, the exact archive is updated based on the exact objective function, features and constraints evaluations added to the DoE.

A new iteration of the Bayesian QD algorithm is therefore carried out.  Due to the highly computational cost context, the stopping criterion is based on the maximum number of objective, feature and constraints evaluations. This number is based on the affordable computational budget defined by the user. Moreover, the algorithm may stop if no new elite in the exact archive are generated for a given number of iterations defined by the user (stagnation of the algorithm). 

\section{Numerical experiments}
\label{sec:numerical_expe}

In order to evaluate the efficiency of the proposed Bayesian QD algorithm, three different analytical problems of increasing complexity and \red{two engineering test problems} are proposed. 
\newred{Three different algorithms are compared: the MAP-Elites algorithm, chosen as it corresponds to a reference algorithm in the family of evolutionary QD algorithms ; the proposed Bayesian QD algorithm with a Gower kernel and the Bayesian QD algorithm with the hypersphere decomposition kernel.} This MAP-Elites algorithm has been modified in order to deal with mixed continuous, discrete and categorical variables and the presence of constraints. The modifications derived in the MAP-Elites are the same as those involved in the proposed Bayesian-QD algorithm to optimize the auxiliary infill problem (Section \ref{discreteMAPElites}). \red{In the considered test problems, because of the presence of mixed continuous, discrete and categorical variables, the high computational cost of feature functions and the presence of optimization constraints, it is not possible to compare with existing Bayesian QD algorithms (SPHEN \cite{hagg2020designing}, SAIL \cite{gaier2017aerodynamic}, BOP-Elites \cite{kent2023bop}) that are not compatible with such problems.}\\

For each test problem, in order to account for the stochastic nature of MAP-Elites algorithm and the initial random DoE for the Bayesian QD algorithm, 10 repetitions are carried out from random initializations. To be able to compare the obtained results, for each repetition, the same initial samples are considered for all the algorithms either under the form of an initial DoE (Bayesian QD algorithms) or an initial population (MAP-Elites). LHS (with a given seed associated to each repetition) is used to generate these initial samples. The number of initial samples corresponds to ten times the dimension of the QD problem \cite{jones1998efficient}.  All the numerical settings for the different algorithms are provided in \ref{Appendix_num_settings}.

Two main indicators are used to compare the algorithms efficiency, in a context of limited number of evaluations of the exact functions. The QD-score, corresponding to the sum of the objective function values for all the illuminated niches, is used to evaluate the overall performance of the algorithm. Moreover, a second indicator corresponding the number of illuminated niches is used to evaluate the ability of the algorithms to create diversity.  The three analytical problems correspond to  modified versions of classical optimization problems with the Rosenbrock function \cite{rosenbrock1960automatic} (Section \ref{Rosen_section}), the Trid function \cite{neumaier1999some} (Section \ref{Trid_section}) and the Styblinski-Tang function \cite{styblinski1990experiments} (Section \ref{Styblinski_section}). The first engineering problem consists of the aerodynamic design of an aircraft wing and is described in Section \ref{Wing_design_pb_section}. \red{The second aerospace problem consists of the design of a two-stage sounding rocket and it is presented in Section \ref{Rocket_section}}. In these problems, as the discrete and categorical variables are handled in the same way in the different algorithms, only categorical variables are involved without loss of generality.

\red{Regarding the hyperparameter tunings and the optimization algorithm settings, numerical details are provided in \ref{Appendix_num_settings}. In terms of computational cost, the training of the Gaussian processes and the solving of the auxiliary infill optimization problem is negligible (in the order of few seconds) compared to the evaluation of the exact simulation model (in the order of few minutes on a cluster of 12 cores for the aerospace problems, see Section \ref{Wing_design_pb_section} for more details).}

\subsection{Rosenbrock problem}\label{Rosen_section}

The Rosenbrock problem is derived from the classical Rosenbrock optimization problem \cite{rosenbrock1960automatic} which has been modified in order to incorporate mixed continuous and categorical variables, a constraint function and two features. This QD problem is in dimension four: two continuous variables ($d_c=2$) and two categorical variables ($d_q=2$). The QD problem is defined as:
\begin{eqnarray}
\forall \tilde{\mathbf{f}} \in \mathcal{F}_t, \;\;\;  \min_{\mathbf{x}^c,\mathbf{x}^q} & & f(\mathbf{x}^c,\mathbf{x}^q) \\ 
\text{s.t.} & &  g_1(\mathbf{x}^c,\mathbf{x}^q) \leq 0 \\ 
& & \mathbf{f}_t(\mathbf{x}^c,\mathbf{x}^q) \in  \tilde{\mathbf{f}} \\ 
& & \mathbf{x}^c_\text{lb} \leq \mathbf{x}^c \leq \mathbf{x}^c_\text{ub} \\
& & \mathbf{x}^q=[x_1^q,x_2^q]^T \in \{0,1,2,3,4,5\}\times\{0,1\} 
\end{eqnarray}

The objective function is defined as:
\begin{equation}
    f(\mathbf{x}^c, \mathbf{x}^q) = - \frac{\displaystyle\sum_{i=1}^{d_c-1} a_q(\mathbf{x}^q)\times \left(x^c_{i+1}-({x_{i}^c})^2\right)^2+b_q(\mathbf{x}^q)\times\left(e_q(\mathbf{x}^q)-x_{i+1}^c\right)^2}{f_q(\mathbf{x}^q)}
\end{equation}
with $\mathbf{x}^c = [x^c_{1},x^c_{2}]^T \in [-5,5]^2$ the vector of the continuous design variables. $a_q, b_q, e_q$ and $f_q$ are variables whose values depend on the value taken by the categorical variables $\mathbf{x}^q= [x^q_{1},x^q_{2}]^T$. The two categorical variables can respectively take six and two levels such that: $x^q_{1} \in \{0,1,2,3,4,5\}$ and $x^q_{2} \in \{0,1\}$. The two features are defined such that $\mathbf{f}_t(\cdot,\cdot) =[f_{t_1}(\cdot,\cdot),f_{t_2}(\cdot,\cdot)]^T$ with:
 \begin{eqnarray}
     f_{t_1}(\mathbf{x}^c, \mathbf{x}^q) & = &  j_q(\mathbf{x}^q)*(x^c_{1} - k_q(\mathbf{x}^q))^{r_q(\mathbf{x}^q)}+s_q(\mathbf{x}^q) \\
     f_{t_2}(\mathbf{x}^c, \mathbf{x}^q) & = &  v_q(\mathbf{x}^q)*(x^c_{2}-t_q(\mathbf{x}^q) )^2 +u_q(\mathbf{x}^q)    
 \end{eqnarray}

The correspondence between the values of the categorical variables and the values of  $a_q, b_q, e_q, f_q, j_q, k_q, r_q, s_q, t_q, u_q$ and $v_q$ is  given in \ref{Rosen_pb_matrix}. 

%by the following matrix:
%\small 
%\begin{equation}
%     \bordermatrix{ & x^q_{1} & x^q_{2}  & \vr a_q & b_q & e_q & f_q & j_q & k_q & r_q & s_q & t_q & u_q & v_q\cr
%      & 0 & 0 & \VR 100 &  1 & 0.7 & 2000 &1 & 0 & 1 & -1.2 & 0 & 0 & -1\cr
%      & 0 & 1 & \VR 103 & 1.6 & 0.2 & 1950 &-1  & 0 & 1 & -0.2 & 0 & 0 & 0.97\cr
%      & 1 & 0 & \VR 98 & 2 & 0.3 & 2100  &1& 0 & 1 & -0.7 &  0 & 0& 0.95 \cr
%      & 1 & 1 & \VR 100 & 1.7 & 0.5 & 2020 &  1 & 0 & 1 & 0.15 & 0 & 0 & 1.1 \cr
%      & 2 & 0 & \VR 95 & 4.7 & 1.5 & 1970 & 1 & 0.15 & 2 & 0 & 0.5 & 0 & -0.8 \cr
%      & 2 & 1 & \VR 97 & 2.4 & 1.2 & 2100 & 1 & -0.55 & 2 & 0.4 & 0 & -0.8 & 0.7 \cr
%      & 3 & 0 & \VR 103 & 1.7 & 2.5 & 2070 & -1 & -1.15 & 2 & 0 & -1.5 & 0 & 1.8 \cr
%      & 3 & 1 & \VR 100 & 0.2 & 1 & 1890 & 1 & -1.3 & 2 & 1.4 & 0 & 0.8& -1.7 \cr
%      & 4 & 0 & \VR 96 & 1.1 & 0.5 & 2140 & -1 & 0.5 & 2 & 0 & -2.3 & 0 & -0.8 \cr
%      & 4 & 1 & \VR 104 & 1.5 & 2 & 1930 & -1 & 1.4 & 2 & -2.4 & 0 & 1.8 & -0.8 \cr
%      & 5 & 0 & \VR 99 & 1.1 & 0.5 & 2140 & 1 & -1.5 & 2 & 0& 2 & 0 & -0.9 \cr
%      & 5 & 1 & \VR 104 & 1.5 & 2 & 2030 & 1 & 1.8 & 2 & 0.4 & 0 & 1 & -0.3 } \qquad
%\end{equation}
%\normalsize
The inequality constraint corresponds to : $$g_1(\mathbf{x}^c, \mathbf{x}^q) = \left((x^c_{1}-0.5)^2+x^c_{2}-5.6\right)/10.\leq0$$

In order to define the archive, the feature functions are discretized in a two-dimensional grid with for each feature axis: $\tilde{f}_{t_1} =[-50,-40,-30,$ $-20,-10,0,10,20,30,40,50]$ and $\tilde{f}_{t_2} =[-50,-40,-30,-20,-10, 0,\;$$10,20,$ $30,40,50,60,70,80]$.

\begin{figure}[!h]
\begin{center}
\includegraphics[width=0.75\textwidth]{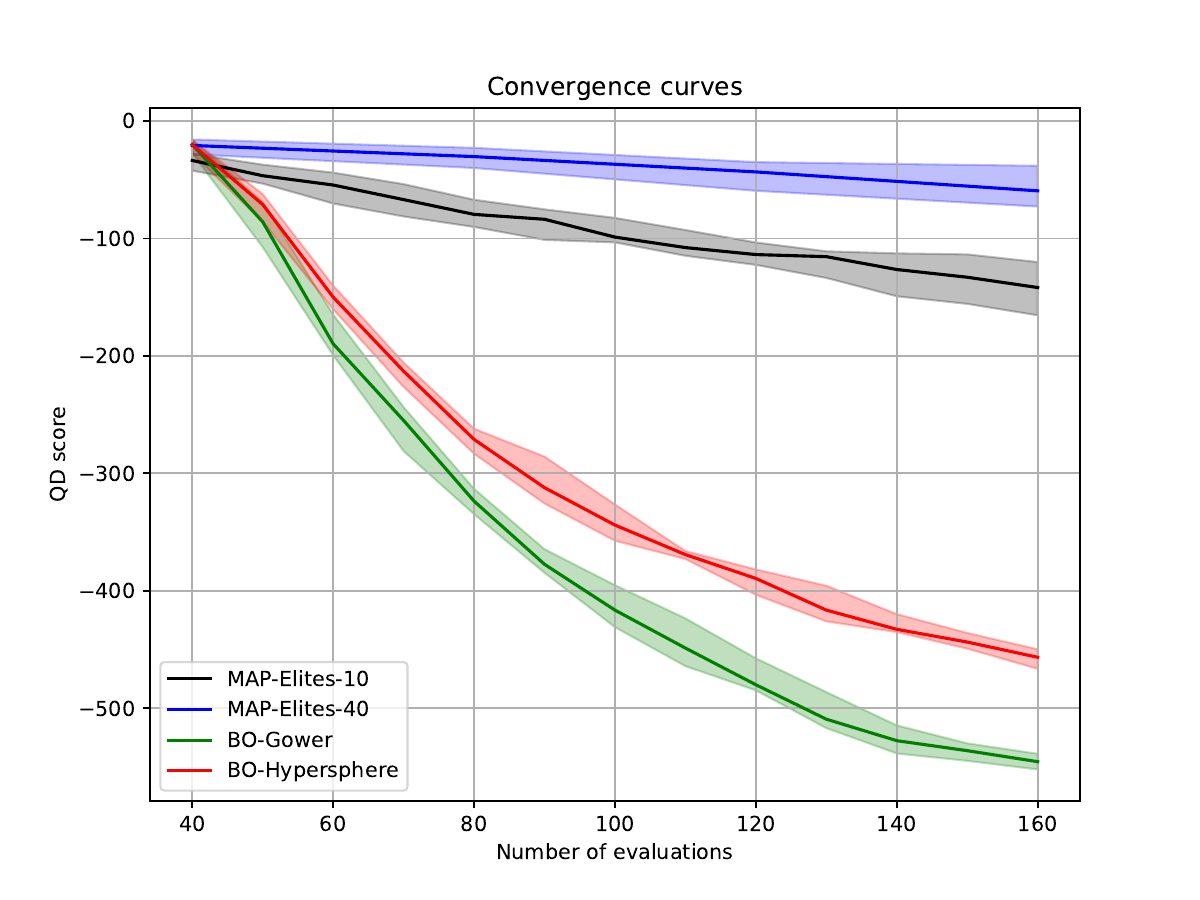}
\caption{Convergence curves (normalized QD score, the lower, the better) for the Rosenbrock problem with MAP-Elites and Bayesian QD algorithm with Gower and hypersphere kernels. For the ten repetitions, the curves correspond to the median whereas the upper and lower limits of the shade areas correspond to the $75^\text{th}$ and $25^\text{th}$ quantiles.}\label{Rosen_Convergence}    
\end{center}
\end{figure}

Figure \ref{Rosen_Convergence} illustrates the QD-score convergence for the MAP-Elites algorithm and the proposed Bayesian QD with the Gower and hypersphere kernels. For the ten repetitions, the curves correspond to the median whereas the upper and lower limits of the shade correspond to the $75^\text{th}$ and $25^\text{th}$ quantiles. In order to compare with Bayesian QD, two sizes of the population are tested with MAP-Elites (population of sizes 10 and 40). A computational budget of 160 evaluations of the exact functions (objective, features and constraint) is allowed. It corresponds to 40 initial samples for the initial DoE plus 120 new added samples during BO. As it can be seen, the convergence of the Bayesian QD algorithms is better in terms of QD-score (the lower, the better) compared to MAP-Elites algorithms. \newred{Indeed, the QD score at the end of allocated budget for BO-Gower is lower than -500 and lower than -450 for BO-Hypersphere while the best MAP-Elites configuration gives a QD-score of around -150.} The MAP-Elites algorithm seems to have a linear convergence with the number of evaluations of the exact functions with the slope influenced by the number of individuals in the population (MAP-Elites with 10 individuals seems to converge better than MAP-Elites with 40 individuals). 

Bayesian QD algorithms have a faster decrease in terms of convergence. Moreover, comparing both Bayesian QD algorithms, in this test case, Bayesian QD with Gower kernel seems to be more efficient than with Bayesian QD the hypersphere kernel. The ability of hypersphere kernel to give a different covariance value for each pair of levels characterizing the categorical variables should give an advantage to the associated GPs to better model the exact functions. However, it is counterbalanced by the larger number of hyperparameters to be determined in the training of the GPs. Indeed, with the Gower kernel, 5 hyperparameters have to be optimized (two lengthscales for the continuous variables, two hyperparameters for the categorical variables and one amplitude parameter) whereas 19 hyperparameters have to be optimized for the hypersphere kernel (two lengthscales for the continuous variables, sixteen hyperparameters for the categorical variables and one amplitude parameter). Therefore, the optimization problem involved with GPs with the hypersphere kernel is more difficult to solve. 

Eventually, in Figure \ref{Rosen_Convergence}, the repetitions with different initial populations for the MAP-Elites and different initial DoEs for the Bayesian QD algorithms illustrate the robustness of the algorithms (the dispersion represented by the shade area around the median is limited). 

\begin{figure}[!h]
\begin{center}
\includegraphics[width=0.75\textwidth]{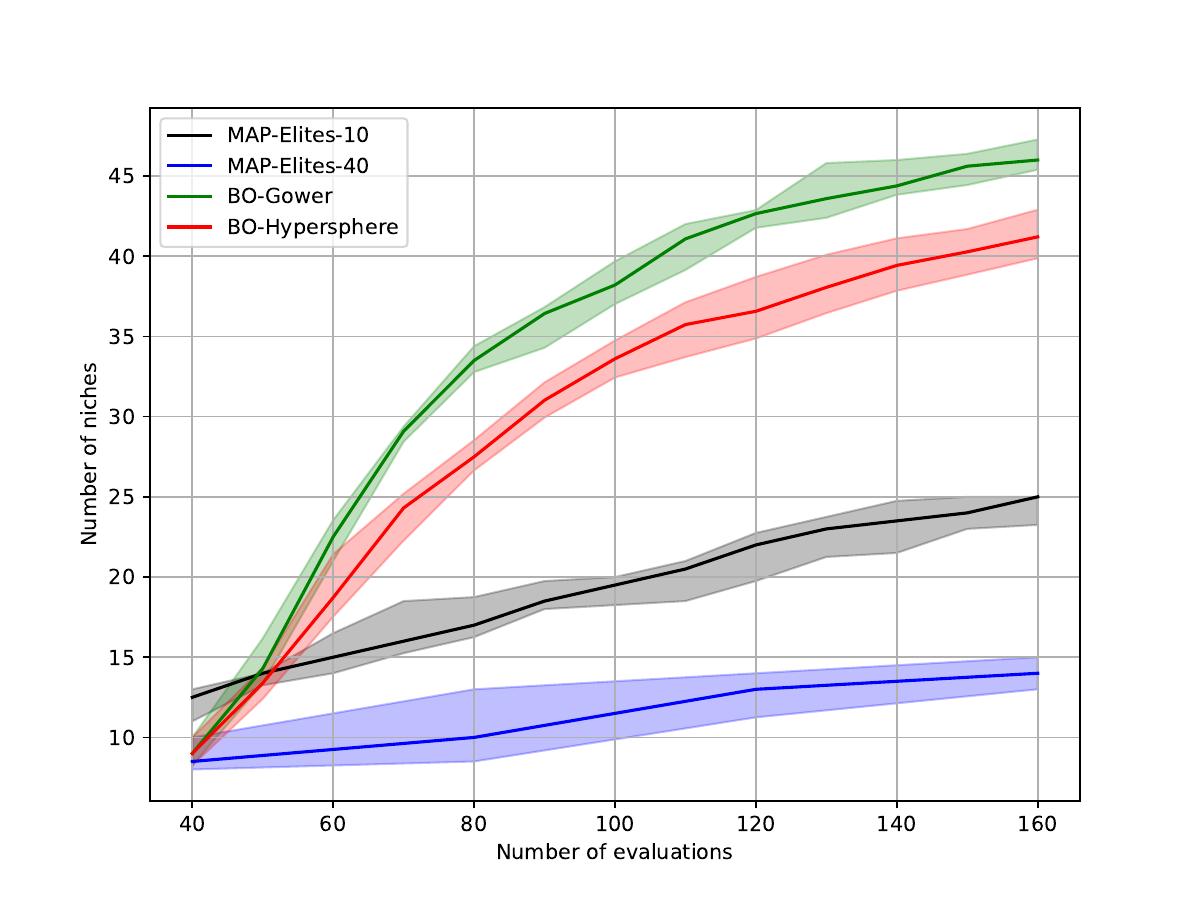}
\caption{Number of discovered niches (the higher, the better) for the Rosenbrock problem with MAP-Elites and Bayesian QD algorithm with Gower and hypersphere kernels. For the ten repetitions, the curves correspond to the median whereas the upper and lower limits of the shade corresponds to the $75^\text{th}$ and $25^\text{th}$ quantiles.}\label{Rosen_Convergence_niche}    
\end{center}
\end{figure}

Regarding the number of discovered niches for the Rosenbrock problem (Figure \ref{Rosen_Convergence_niche}), the Bayesian QD algorithms provide better results than MAP-Elites algorithm. Indeed, QD algorithms discover almost twice as much niches compared to MAP-Elites on this test problem. Bayesian QD with the Gower kernel discovers the largest number of niches for all the repetitions. 

\begin{figure}[!h]
\begin{center}
\includegraphics[width=0.49\textwidth]{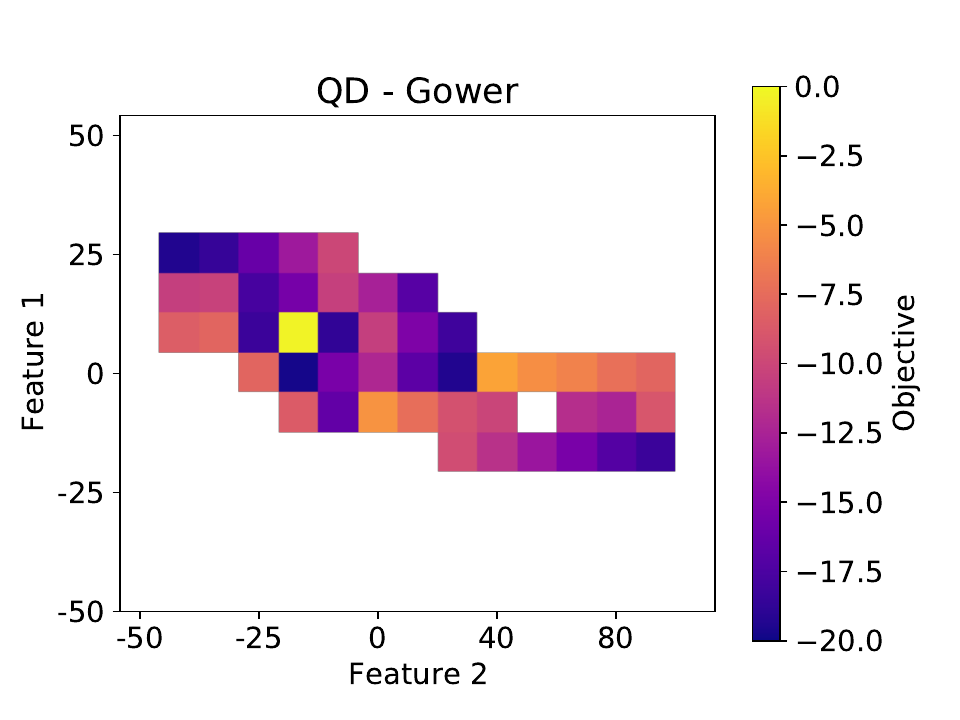}
\includegraphics[width=0.49\textwidth]{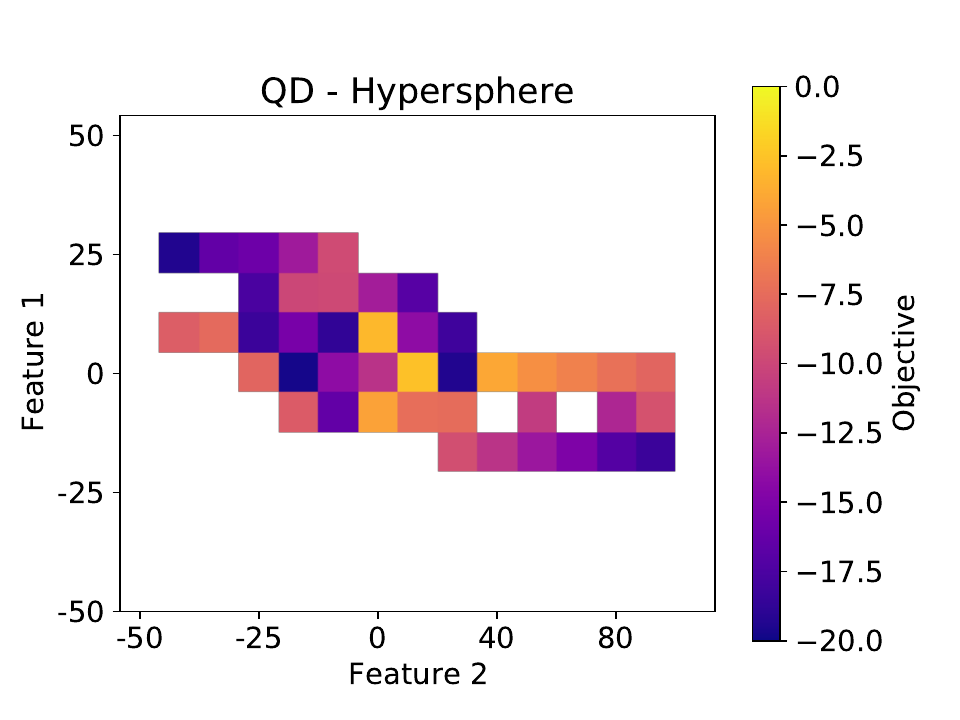}
\includegraphics[width=0.49\textwidth]{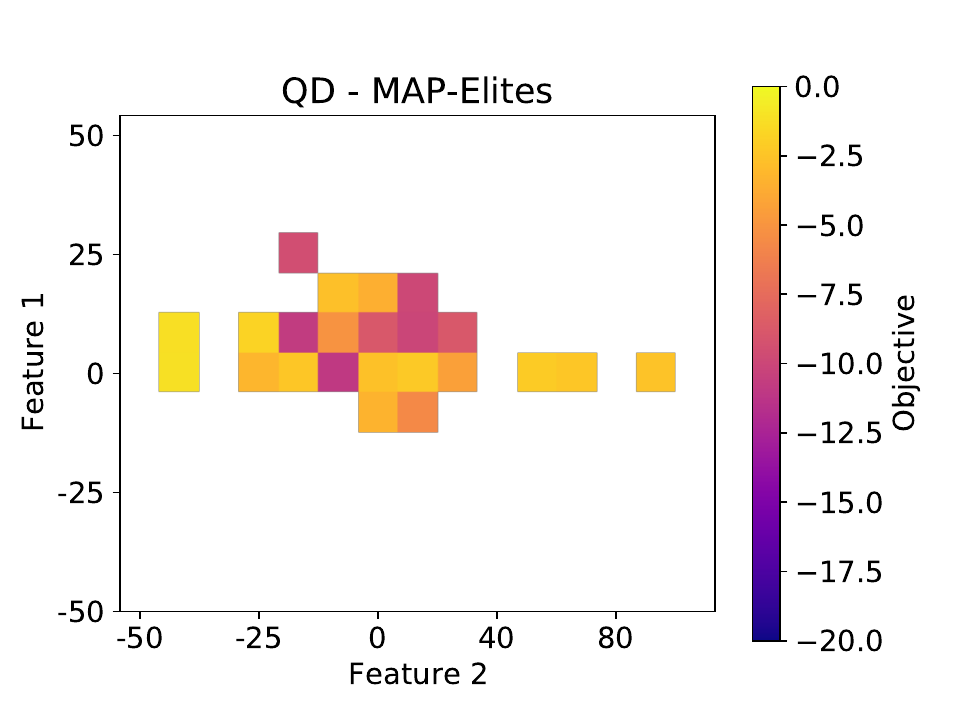}
\includegraphics[width=0.49\textwidth]{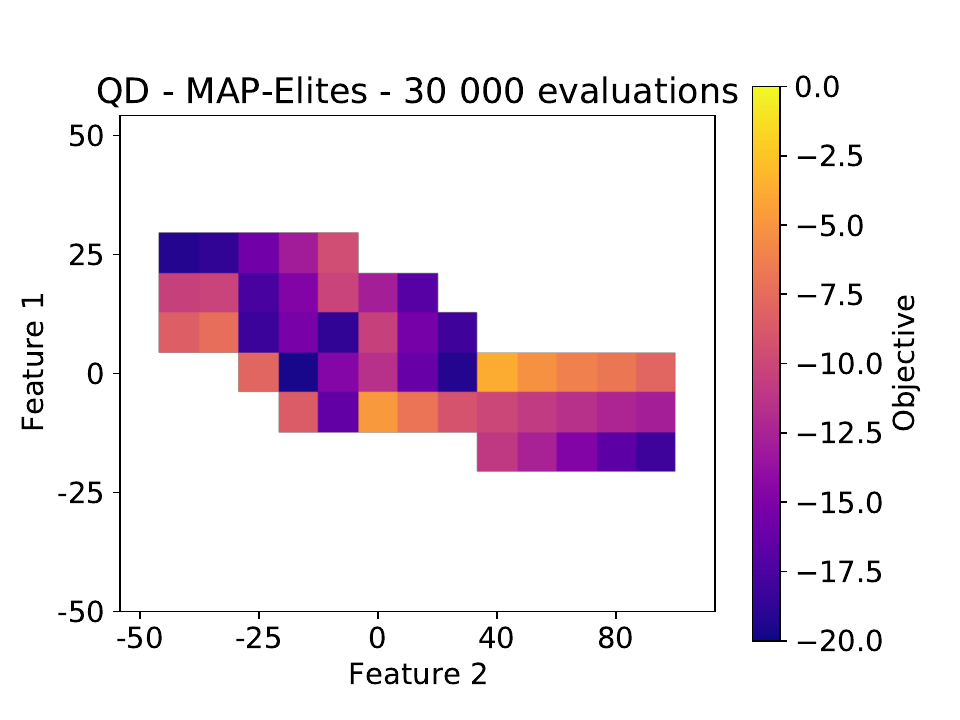}
\caption{Final archive for the Rosenbrock problem obtained by Bayesian QD with the Gower kernel with 160 evaluations (top left), by Bayesian QD with the hypersphere kernel with 160 evaluations (top right), by QD MAP-Elites with 160 evaluations (bottom left) and with QD MAP-Elites with 30 000 evaluations (bottom right) }\label{Rosen_final_archive}    
\end{center}
\end{figure}

\begin{figure}[!h]
\begin{center}
\includegraphics[width=1.\textwidth]{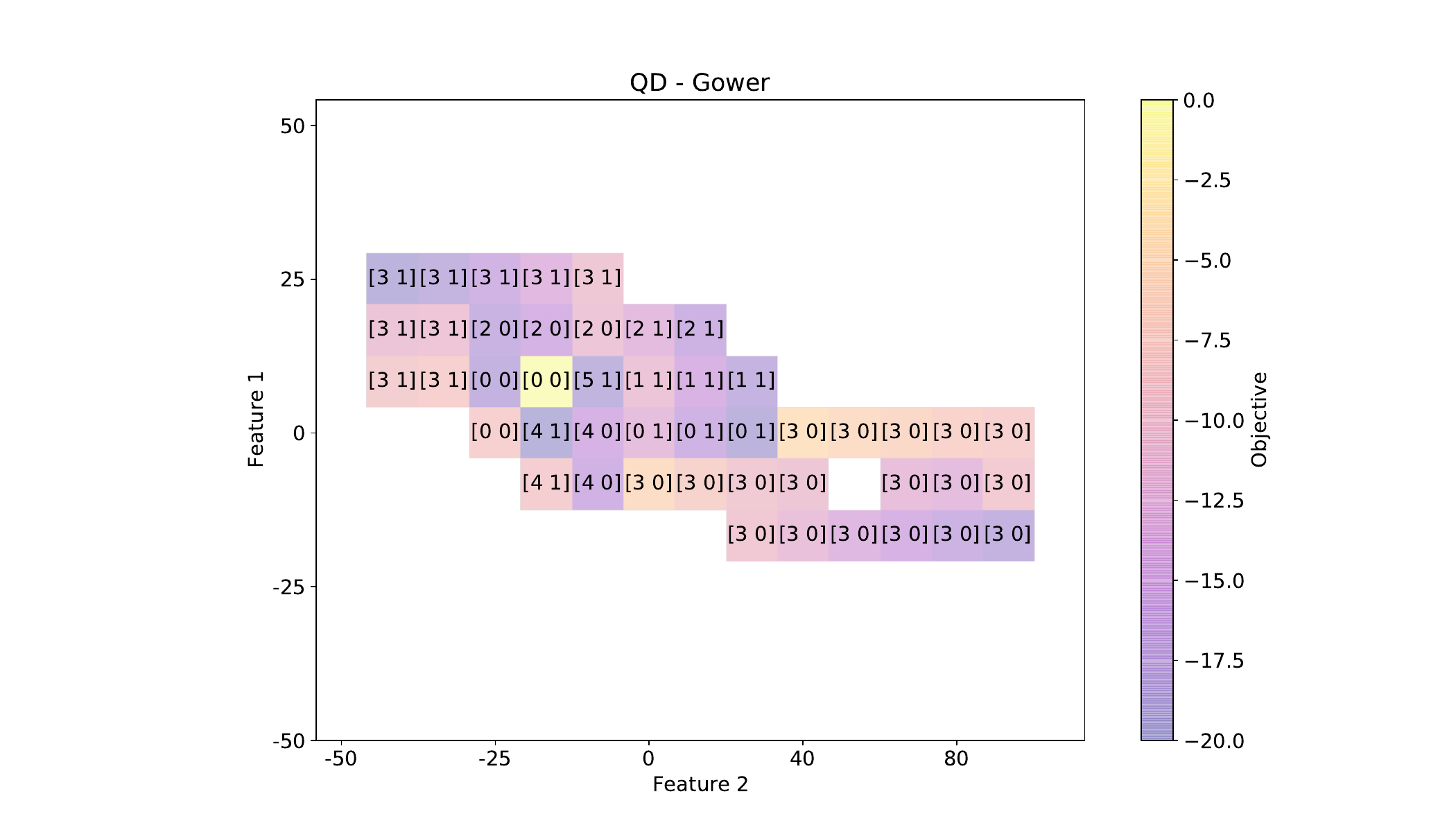}
\caption{Final archive for the Rosenbrock problem obtained by Bayesian QD with Gower kernel with 160 evaluations. The numbers correspond to value the best values of categorical variables $\mathbf{x}^q=[x_1^q,x_2^q]^T$ in each niche}\label{Rosen_Gower_final_map_cate}    
\end{center}
\end{figure}

\begin{figure}[!h]
\begin{center}
\includegraphics[width=1.\textwidth]{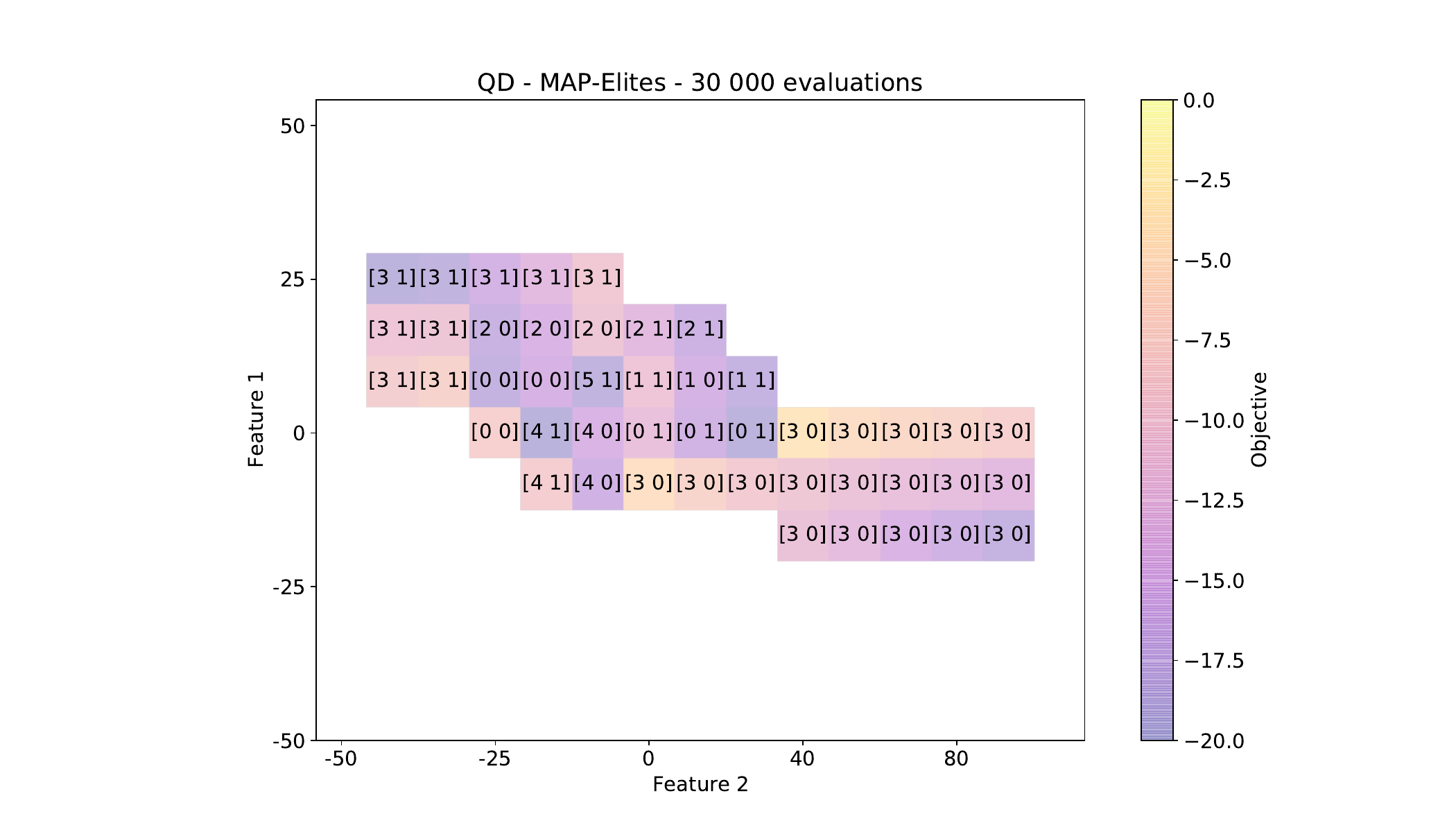}
\caption{Final archive for the Rosenbrock problem obtained by MAP-Elites with 30 000 evaluations. The numbers correspond to value the best values of categorical variables $\mathbf{x}^q=[x_1^q,x_2^q]^T$ in each niche}\label{Rosen_MAPElites_final_map_cate}    
\end{center}
\end{figure}

In Figure \ref{Rosen_final_archive}, the final QD archive obtained (for one representative repetition) with the 160 evaluations of the exact functions (corresponding to the initial DoEs in addition with the chosen candidates during the enrichment process) are represented for the Bayesian QD algorithms (QD-Gower and QD-Hypersphere) and for the MAP-Elites. In order to compare, the final archive obtained with a MAP-Elites algorithm involving 30 000 evaluations of the exact functions is also represented as the reference map. It can be seen that with only 160 evaluations, the final QD archive provided by Bayesian QD algorithms are close to the archive obtained by MAP-Elites with 30 000 evaluations, that illustrates the efficiency of the proposed algorithm. However, the map obtained with the QD MAP-Elites with 160 evaluations is far from the one obtained with Bayesian QD algorithms both in terms of niches discovered and quality of the best individual in each discovered niche.

The extension of the Bayesian QD algorithms to deal with mixed continuous, discrete and categorical variables offers the possibility to converge into different categorical optimal solutions for different niches. Indeed, for the Rosenbrock problem, in Figures \ref{Rosen_Gower_final_map_cate} and \ref{Rosen_MAPElites_final_map_cate}, the optimal solutions in terms of categorical variables $\mathbf{x}^q=[x_1^q,x_2^q]^T$ are displayed for each niche for respectively the Bayesian QD Gower algorithm with 160 evaluations and for the MAP-Elites algorithm with 30 000 evaluations. It can be seen that depending on the niches, optimal solutions may belong to different categories. Moreover, the Bayesian QD algorithm provides optimal solutions in terms of categorical variable values (with only 160 exact function evaluations) that are similar to the ones found by MAP-Elites algorithm (with 30 000 evaluations). The ability to deal with categorical variables is interesting as this allows to converge to different categorical solutions into different niches and therefore to offer diversity with respect to the categorical variables. This is of particular interest in the field of engineering design with categorical variables representing for instance technological choices or architectural choices. This aspect is further explored in the aerospace design problems in Sections \ref{Wing_design_pb_section} and \ref{Rocket_section}.\\

\begin{figure}[!h]
\begin{center}
\includegraphics[width=0.75\textwidth]{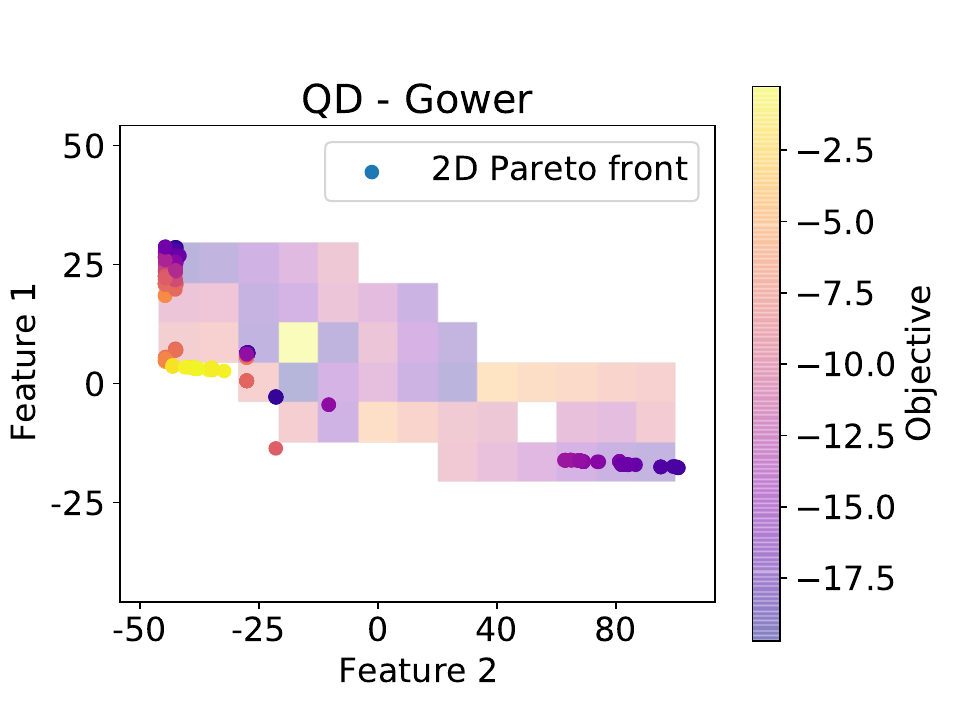}
\caption{\newred{Pareto front defined by the dots in the feature space, and colored by the objective function values.}}\label{Pareto}    
\end{center}
\end{figure}

\newred{One could be interested by comparing the results provided by QD algorithm with results from obtained multi-objective optimizations considering the features as additional objectives. For this test-case, a constrained three-objective problem is solved using a mixed continuous-discrete NSGAIII algorithm with a full convergence \cite{deb2013evolutionary}. The obtained results is under the form of a three-dimensional Pareto front, representing the non dominated solutions between the initial objective function and the two features. Figure \ref{Pareto} illustrates the two-dimensional Pareto front in the features space to compare with the feature map obtained with QD algorithm. It can be noted that in this test case, because the feature and objective functions are antagonistic, the Pareto front corresponds to the lower niches in the feature map. It could be possible to obtain other niches by considering higher rank Pareto fronts (corresponding to dominated solutions). However, this comparison with multi-objective approach is only valuable when the objective and features functions are antagonistic, which is not always the case it can be seen in the physical engineering problems in Sections \ref{Wing_design_pb_section} and \ref{Rocket_section}.}

\newpage
\subsection{Trid problem}\label{Trid_section}

The Trid problem is derived from the classical Trid optimization problem \cite{neumaier1999some} (also called Neumaier number 3 function) which has been modified in order to involve mixed continuous and categorical variables, a constraint and two features. This QD problem is of dimension six: four continuous variables ($d_c=4)$ and two categorical variables. The QD problem is defined by:
\begin{eqnarray}
\forall \tilde{\mathbf{f}} \in \mathcal{F}_t, \;\;\;  \min_{\mathbf{x}^c,\mathbf{x}^q} & & f(\mathbf{x}^c,\mathbf{x}^q) \\ 
\text{s.t.} & &  g_1(\mathbf{x}^c,\mathbf{x}^q) \leq 0 \\ 
& & \mathbf{f}_t(\mathbf{x}^c,\mathbf{x}^q) \in  \tilde{\mathbf{f}} \\ 
& & \mathbf{x}^c_\text{lb} \leq \mathbf{x}^c \leq \mathbf{x}^c_\text{ub} \\
& & \mathbf{x}^q=[x_1^q,x_2^q]^T \in \{0,1,2\}\times\{0,1\}
\end{eqnarray}

The objective function is defined such that:
\begin{equation}
    f(\mathbf{x}^c, \mathbf{x}^q) = \sum_{i=1}^{d_c} a_q(\mathbf{x}^q) \times \left(x^c_{i}-b_q(\mathbf{x}^q)\right)^2 -\sum_{i=2}^{d_c} c_q(\mathbf{x}^q) \times x^c_{i}x^c_{i-1}
\end{equation}
with $\mathbf{x}^c = [x^c_{1},x^c_{2},x^c_{2},x^c_{4}]^T \in [0,1]^4$ the vector of the continuous design variables. $a_q, b_q, c_q$ are variables whose values depend on the value taken by the categorical variables $\mathbf{x}^q= [x^q_{1},x^q_{2}]^T$. The two categorical variables can respectively take three and two levels such that: $x^q_{1} \in \{0,1,2\}$ and $x^q_{2} \in \{0,1\}$. 

The two feature functions are defined such that $\mathbf{f}_t(\cdot,\cdot) =[f_{t_1}(\cdot,\cdot),f_{t_2}(\cdot,\cdot)]^T$ with:
 \begin{eqnarray}
     f_{t_1}(\mathbf{x}^c, \mathbf{x}^q) & = &  e_q(\mathbf{x}^q)*x^c_{3}+(f_q(\mathbf{x}^q)\times x^c_{1}-j_q(\mathbf{x}^q))^2+k_q(\mathbf{x}^q)\times x^c_{2} \\
     f_{t_2}(\mathbf{x}^c, \mathbf{x}^q) & = &  r_q(\mathbf{x}^q)\times x^c_{2}-s_q(\mathbf{x}^q)+(t_q(\mathbf{x}^q)\times x^c_{4}\times x^c_{3}-u_q(\mathbf{x}^q))^2     
 \end{eqnarray}

The correspondence between the values of the categorical variables and the values of  $a_q, b_q, c_q, e_q, f_q, j_q, k_q, r_q, s_q, t_q$ and $u_q$ is given in \ref{Trid_pb_matrix}.

%by the following matrix:
%\small 
%\begin{equation}
%     \bordermatrix{ & x^q_{1} & x^q_{2}  & \vr a_q & b_q & c_q & e_q & f_q & j_q & k_q & r_q & s_q & t_q & u_q \cr
%      & 0 & 0 & \VR 1. &  1 & 1 & 1 & 1 & 0.7 & 1 & 1 & 1.5 & 1 & 0.4\cr
%      & 1 & 0 & \VR 0.95 & 1 & 1.1 & 0.8  &1 & 0.4 & 1.1 & 1 & 1.9 & 1 & 0.1 \cr
%      & 2 & 0 & \VR 1 & 1.3 & 0.97 & 1.1 & 0.8  & 0.1 & 1 & 0.9 &1.5. & 1.1 & 0.4\cr
%      & 0 & 1 & \VR 1.1 & 0.7 & 1 & 1 &  1 & 0.7 & 1 & 1 & 0.7 & 1 & 1.4 \cr
%      & 1 & 1 & \VR 0.7 & 0.5 & 0.4 & 1.5 & 1 & 1.7 & 0.7& 0.7 & 0.5 & 1 & 0.9 \cr
%      & 2 & 1 & \VR 0.7 & 1 & 1.5 & 1 & 1.3 & 0.91 & 1 & 1 & 1.5 & 0.7 & 0.1 } \qquad
%\end{equation}
The inequality constraint is defined as: $$g_1(\mathbf{x}^c, \mathbf{x}^q) = (x^c_{1}-0.4)^2+1.5\times x^c_{3}-1.3\leq0$$

In order to define the archive, the feature functions are discretized in a two-dimensional grid with for each feature axis: $\tilde{f}_{t_1} =[-1.5,-0.5,0.5,1.5,2.5,$ $3.5,4.5]$ and $\tilde{f}_{t_2} =[-2.5,-1.5,-0.5,0.5,1.5,2.5]$.

\normalsize
\begin{figure}[!h]
\begin{center}
\includegraphics[width=0.75\textwidth]{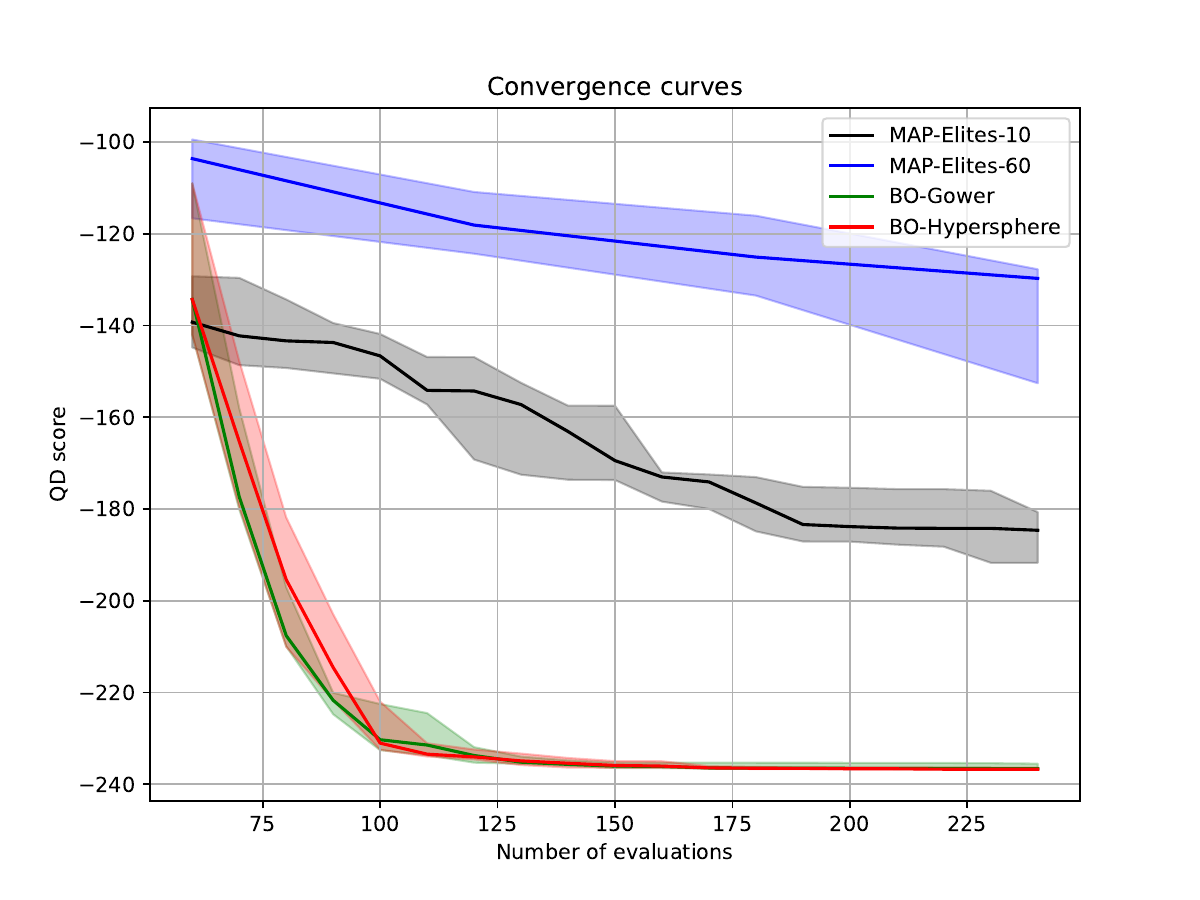}
\caption{Convergence curves (normalized QD-score, the lower, the better) for the Trid problem with MAP-Elites and Bayesian QD algorithm with Gower and hypersphere kernels. For the ten repetitions, the curves correspond to the median whereas the upper and lower limits of the shade area corresponds to the $75^\text{th}$ and $25^\text{th}$ quantiles.}\label{Trid_Convergence}    
\end{center}
\end{figure}

\begin{figure}[!h]
\begin{center}
\includegraphics[width=0.75\textwidth]{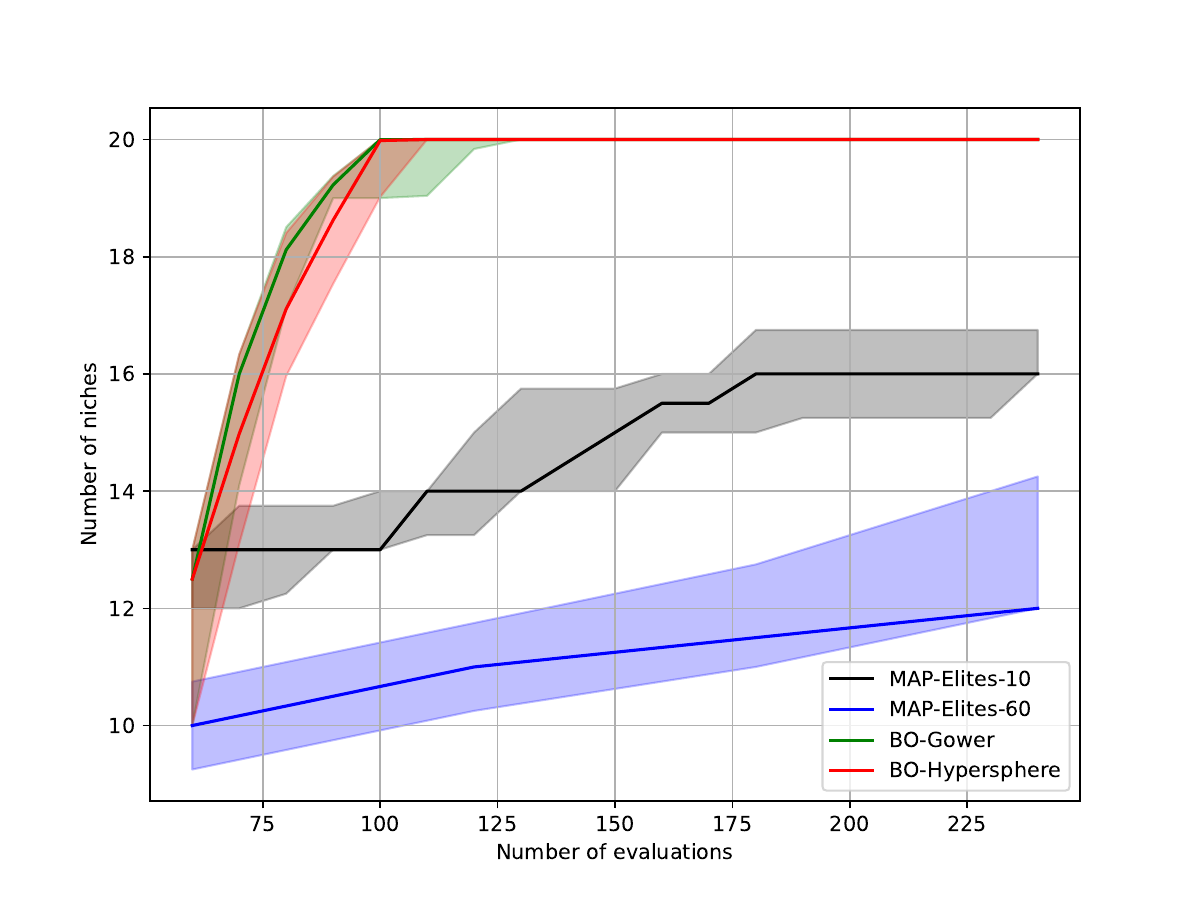}
\caption{Number of discovered niches (the higher, the better) for the Trid problem with MAP-Elites and Bayesian QD algorithm with Gower and hypersphere kernels. For the ten repetitions, the curves correspond to the median whereas the upper and lower limits of the shade corresponds to the $75^\text{th}$ and $25^\text{th}$ quantiles.}\label{Trid_Convergence_niche}    
\end{center}
\end{figure}

Both Figures \ref{Trid_Convergence} and \ref{Trid_Convergence_niche} representing the convergence for the Trid problem respectively of the QD-score and the number of discover niches illustrate the better efficiency of the Bayesian QD algorithms compared to MAP-Elites. Similar trends as for the Rosenbrock problem may be described for this problem in dimension 6. The MAP-Elites algorithms seem to converge with a linear trend (determined by the population size) whereas the Bayesian QD algorithms converge faster with respect to the number of exact function evaluations. Moreover, the Bayesian QD algorithms discover a larger number of niches in the feature space. \newred{The full convergence of BO algorithms is reached after a hundred evaluations whereas the best MAP-Elites provides at the end of simulation budget (240 evaluations) a number of illuminated niches of 16 instead of 20. The dispersion of the results of BO-algorithms for the different repetitions is very low compared to MAP-Elites algorithms. Indeed, for all the repetitions, BO-algorithms illuminate all the niches (20).}

\begin{figure}[!h]
\begin{center}
\includegraphics[width=0.49\textwidth]{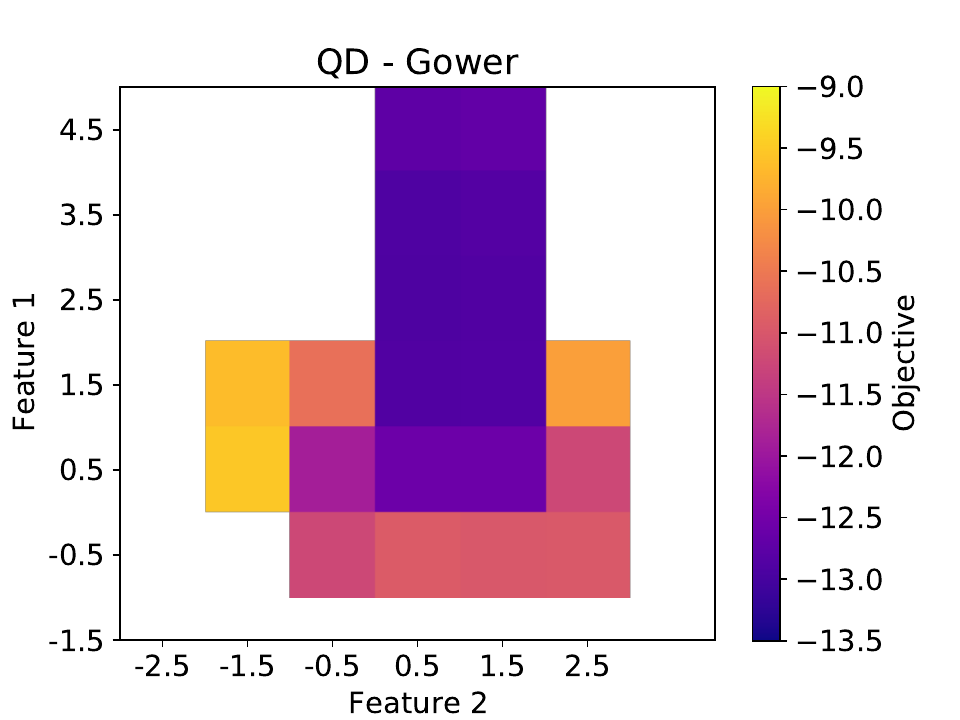}
\includegraphics[width=0.49\textwidth]{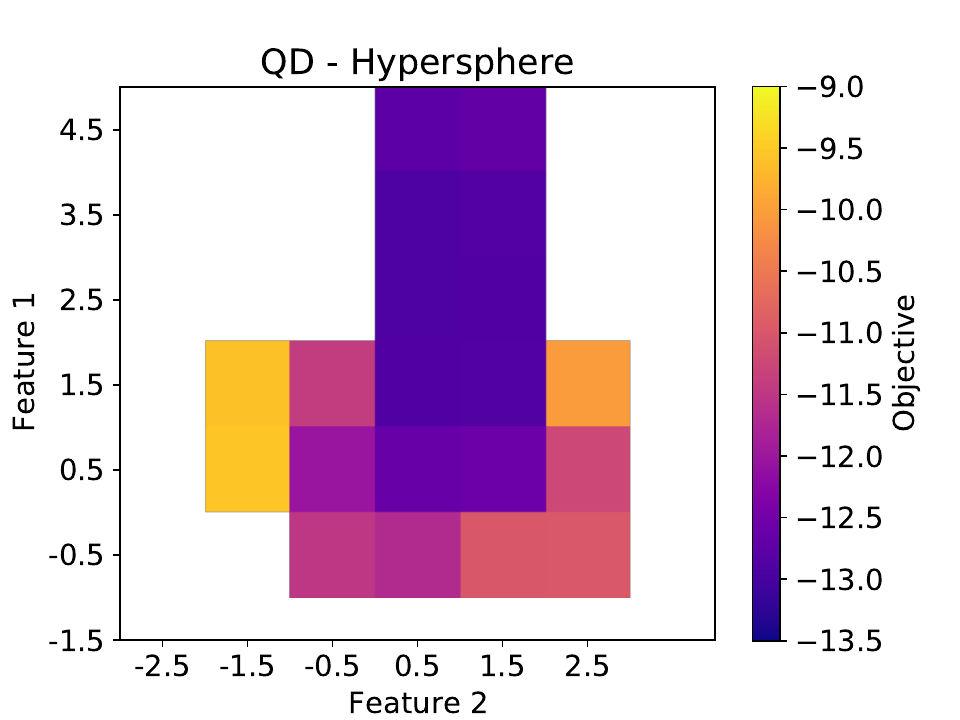}
\includegraphics[width=0.49\textwidth]{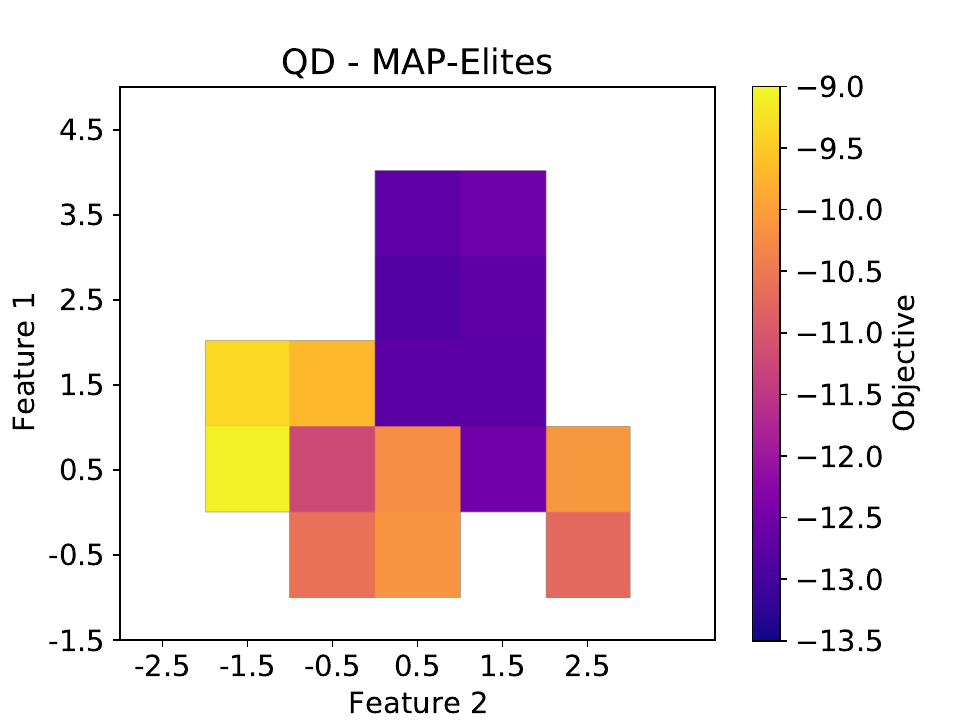}
\includegraphics[width=0.49\textwidth]{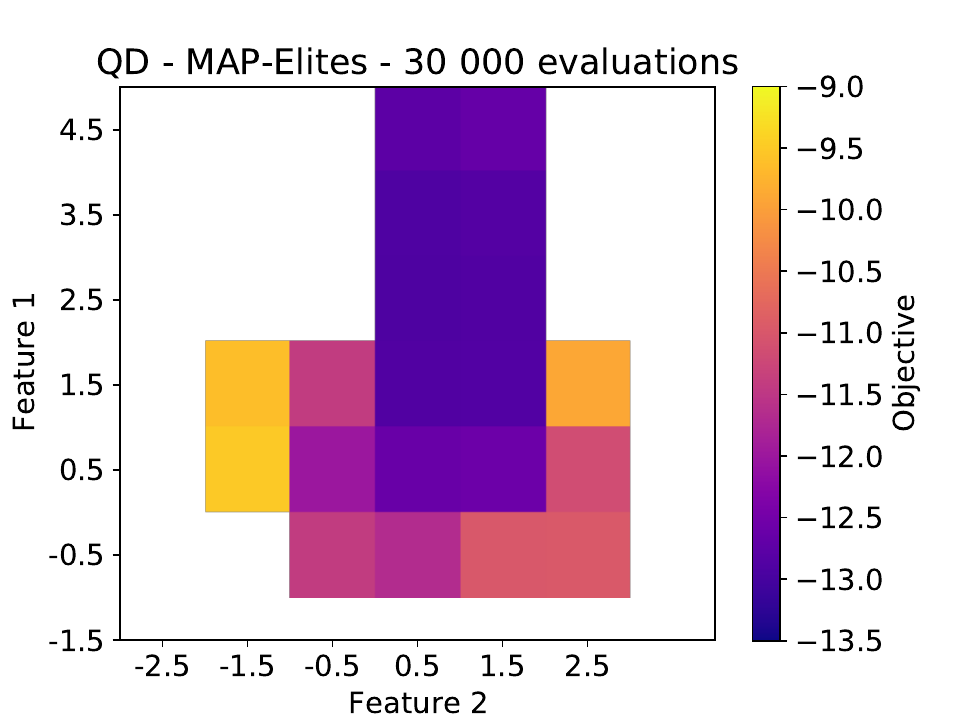}
\caption{Final archive for the Trid problem obtained by Bayesian QD with the Gower kernel with 240 evaluations (top left), by Bayesian QD with the hypersphere kernel with 240 evaluations (top right), by QD MAP-Elites with 240 evaluations (bottom left) and with QD MAP-Elites with 30 000 evaluations (bottom right) }\label{Trid_final_archive}    
\end{center}
\end{figure}

Eventually, in Figure \ref{Trid_final_archive}, the final archives obtained with 240 evaluations of the exact functions for the Bayesian QD-algorithms are very similar to the one obtained with MAP-Elites algorithm after 30 000 evaluations. However, MAP-Elites with only 240 evaluations of the exact functions is still far from the converged archive illustrating the interest of QD algorithms based on GPs with adapted covariance models to handle mixed continuous and categorical variables.  

\newpage
\subsection{Styblinski-Tang problem}\label{Styblinski_section}

The Styblinski-Tang problem is derived from the classical Styblinski-Tang optimization problem \cite{styblinski1990experiments} which has been modified in order to incorporate mixed continuous and categorical variables, two constraints and two features. This QD problem is in dimension nine: six continuous variables ($d_c=6)$ and three categorical variables. The QD problem is defined as:
\begin{eqnarray}
\forall \tilde{\mathbf{f}} \in \mathcal{F}_t, \;\;\;  \min_{\mathbf{x}^c,\mathbf{x}^q} & & f(\mathbf{x}^c,\mathbf{x}^q) \\ 
\text{s.t.} & &  g_1(\mathbf{x}^c,\mathbf{x}^q) \leq 0 \\ 
& &  g_2(\mathbf{x}^c,\mathbf{x}^q) \leq 0 \\ 
& & \mathbf{f}_t(\mathbf{x}^c,\mathbf{x}^q) \in  \tilde{\mathbf{f}} \\ 
& & \mathbf{x}^c_\text{lb} \leq \mathbf{x}^c \leq \mathbf{x}^c_\text{ub} \\
& & \mathbf{x}^q=[x_1^q,x_2^q,x_3^q] \in \{0,1\}^3
\end{eqnarray}

The objective function is defined such that:
\begin{equation}
    f(\mathbf{x}^c, \mathbf{x}^q) = \sum_{i=1}^{d_c} \left(a_q(\mathbf{x}^q) \times \left(x^c_{i}\right)^4-b_q(\mathbf{x}^q) \times \left(x^c_{i}\right)^2 +c_q(\mathbf{x}^q)\times x^c_{i}\right)
\end{equation}
with $\mathbf{x}^c = [x^c_{1},x^c_{2},x^c_{2},x^c_{4},x^c_{5},x^c_{6}]^T \in [0,1]^6$ the vector of the continuous variables. $a_q, b_q, c_q$ are variables whose values depend on the value taken by the categorical variables $\mathbf{x}^q= [x^q_{1},x^q_{2},,x^q_{3}]^T$. The three categorical variables can respectively take two levels such that: $x^q_{1} \in \{0,1\}$,  $x^q_{2} \in \{0,1\}$ and $x^q_{3} \in \{0,1\}$. The two features are defined such that $\mathbf{f}_t(\cdot,\cdot) =[f_{t_1}(\cdot,\cdot),f_{t_2}(\cdot,\cdot)]^T$ with:
 \begin{eqnarray}
     f_{t_1}(\mathbf{x}^c, \mathbf{x}^q) & = & (x^c_{3}- e_q(\mathbf{x}^q))^2+(x^c_{5}- f_q(\mathbf{x}^q))^2 \\
     f_{t_2}(\mathbf{x}^c, \mathbf{x}^q) & = & x^c_{2}+j_q(\mathbf{x}^q)+(x^c_{4}-k_q(\mathbf{x}^q))^2
 \end{eqnarray}

The correspondence between the values of the categorical variables and the values of  $a_q, b_q, c_q, e_q, f_q, j_q, $ and $k_q$ are given in \ref{Sty_pb_matrix}. 

%by the following matrix:
%\small 
%\begin{equation}
%     \bordermatrix{ & x^q_{1} & x^q_{2} & x^q_{3} & \vr a_q & b_q & c_q & e_q & f_q & j_q & k_q  \cr
%      & 0 & 0 & 0 & \VR 1 &  16 & 5 & 1.2 & 0.7 & 3.5 & 0.7 \cr
%      & 1 & 0 & 0 & \VR 1.1 & 18 & 6.1 & 1.4  & 0.9 & 3.8 & 0.2 \cr
%      & 1 & 1 & 0 & \VR 0.95 & 17 & 4.9 & 1.7 & 1.3  & 2.8 & 0.7 \cr
%      & 0 & 1 & 0 & \VR 0.94 & 12 & 6.9 & 1.4 &  0.2 & 1.4 & 0.2  \cr
%      & 0 & 0 & 1 & \VR 0.75 & 10 & 7 & 2.2 & 1.7 & 1.5 & 0.5 \cr
%      & 1 & 0 & 1 & \VR 1.2 & 19 & 4.2 & 1.5  & 2.9 & 1.4 & 1.2 \cr
%      & 1 & 1 & 1 & \VR 0.97 & 12 & 1.9 & 0.7 & 2.3  & 3.8 & 0.4 \cr
%      & 0 & 1 & 1 & \VR 1.1 & 18 & 4.2 & 1.9 &  0.7 & 2.7 & 0.4 } \qquad
%\end{equation}
The two inequality constraints are defined as: $$g_1(\mathbf{x}^c, \mathbf{x}^q) = x^c_{1} + x^c_{2} -1\leq0$$
$$g_2(\mathbf{x}^c, \mathbf{x}^q) = x^c_{4} + x^c_{6} -2\leq0$$

In order to define the archive, the features are discretized in a two-dimensional grid with for each feature axis: $\tilde{f}_{t_1} =[0,2,4,6,8,10,12]$ and $\tilde{f}_{t_2} =[-5,-3,-1,1,3,5]$.

\begin{figure}[!h]
\begin{center}
\includegraphics[width=0.75\textwidth]{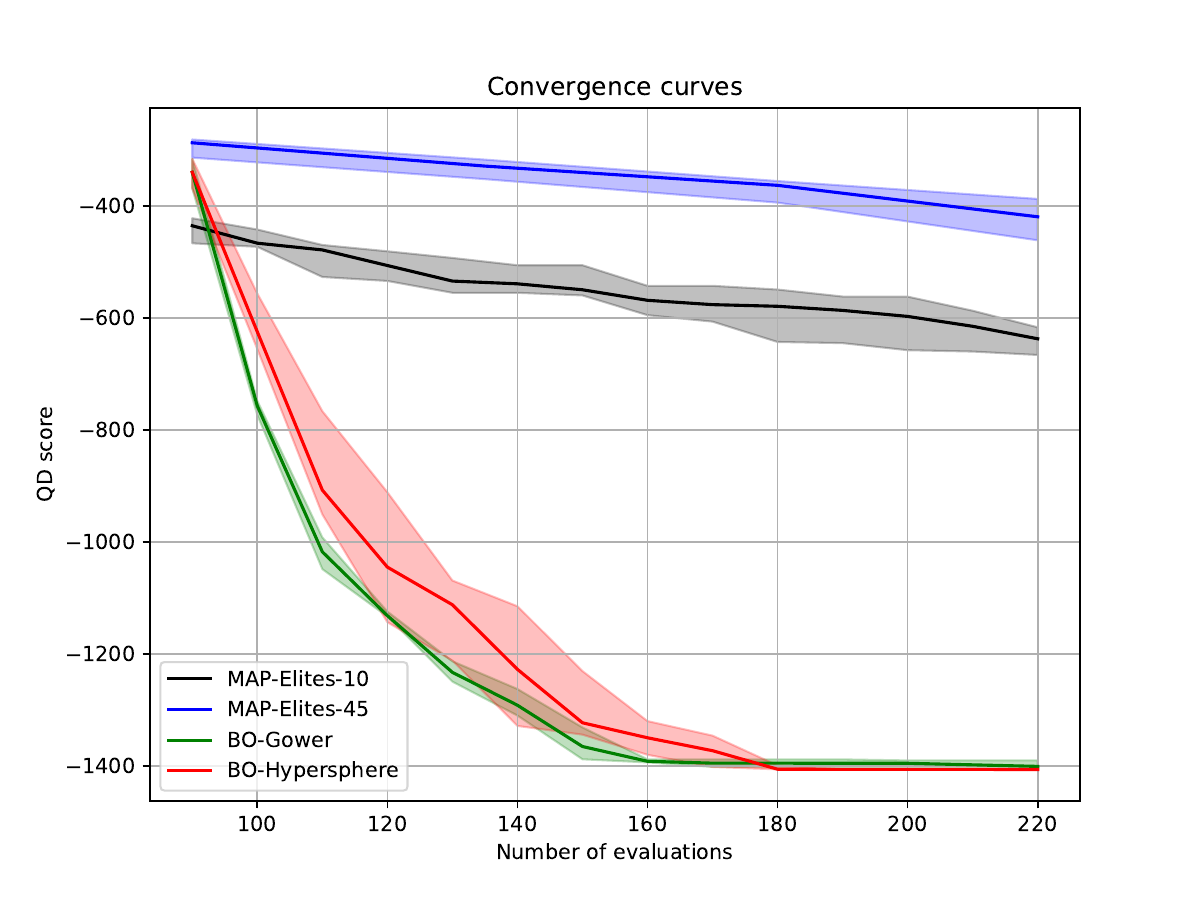}
\caption{Convergence curves (normalized QD score, the lower, the better) for the Styblinski-Tang problem with MAP-Elites and Bayesian QD algorithm with Gower and hypersphere kernels. For the ten repetitions, the curves correspond to the median whereas the upper and lower limits of the shade area corresponds to the $75^\text{th}$ and $25^\text{th}$ quantiles.}\label{Tang_Convergence}    
\end{center}
\end{figure}

\begin{figure}[!h]
\begin{center}
\includegraphics[width=0.75\textwidth]{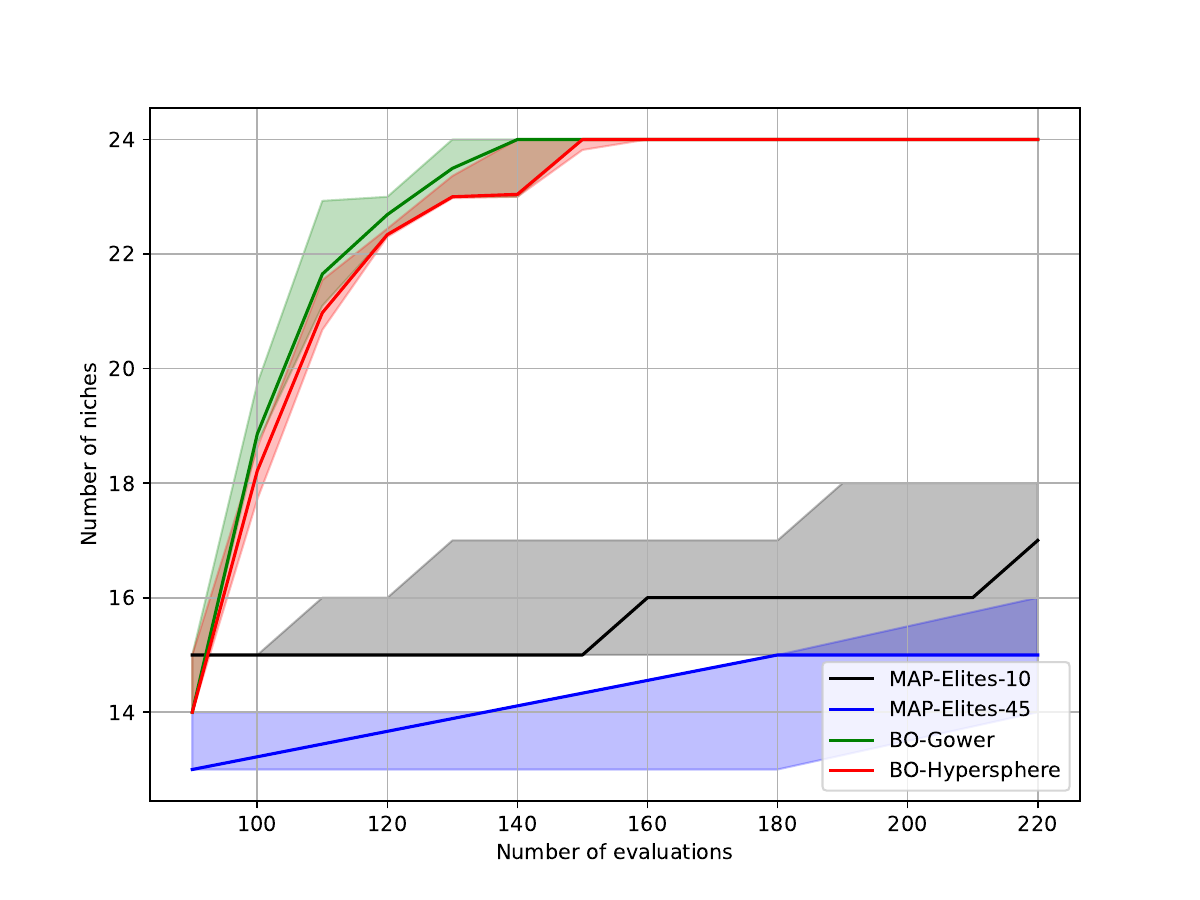}
\caption{Number of discovered niches (the higher, the better) for the Styblinski-Tang problem with MAP-Elites and Bayesian QD algorithm with Gower and hypersphere kernels. For the ten repetitions, the curves correspond to the median whereas the upper and lower limits of the shade corresponds to the $75^\text{th}$ and $25^\text{th}$ quantiles.}\label{Tang_Convergence_niche}    
\end{center}
\end{figure}

This QD problem is of higher complexity in terms of number of continuous variables, number of categorical variables and number of constraints. The convergence curves for the QD-score (Figure \ref{Tang_Convergence}) and the number of discovered niches (Figure \ref{Tang_Convergence_niche}) provide the same tendencies as for the two previous analytical problems. \newred{ Indeed, the dispersion of the illuminated niches for BO-algorithms is null at the end of the simulation budget whereas MAP-Elites algorithms present a large dispersion in the number of illuminated niches.}

The Bayesian QD algorithms converge faster toward a better solution compared to MAP-Elites algorithms that displays a linear convergence rate in terms of number of exact function evaluations. Both Bayesian QD algorithms converge for all the repetitions to the same solution illustrating the robustness to the initial DoE. Moreover, the number of illuminated niches corresponds (in median) to 17 for MAP-Elites with 10 individuals whereas it increases up to 24 for the Bayesian QD algorithms. 

\begin{figure}[!h]
\begin{center}
\includegraphics[width=0.49\textwidth]{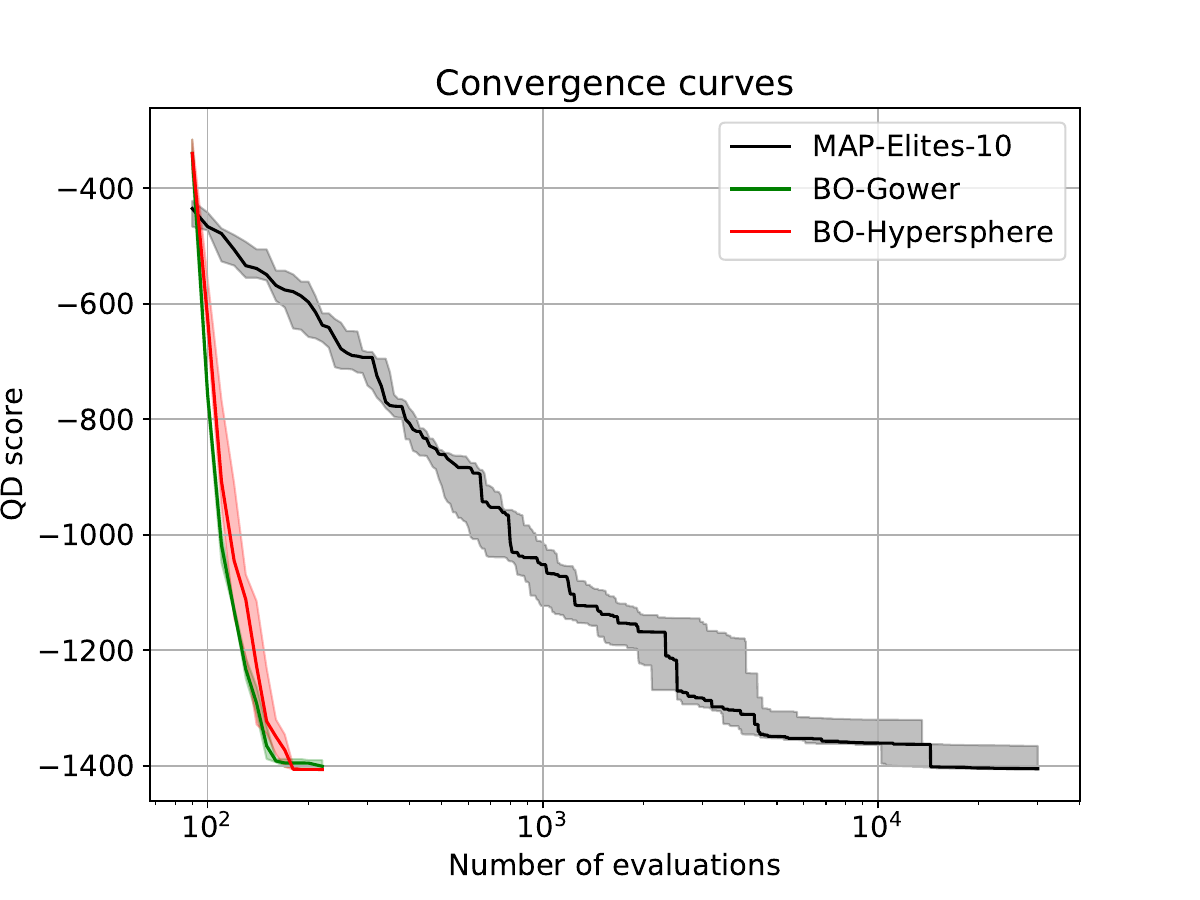}
\includegraphics[width=0.49\textwidth]{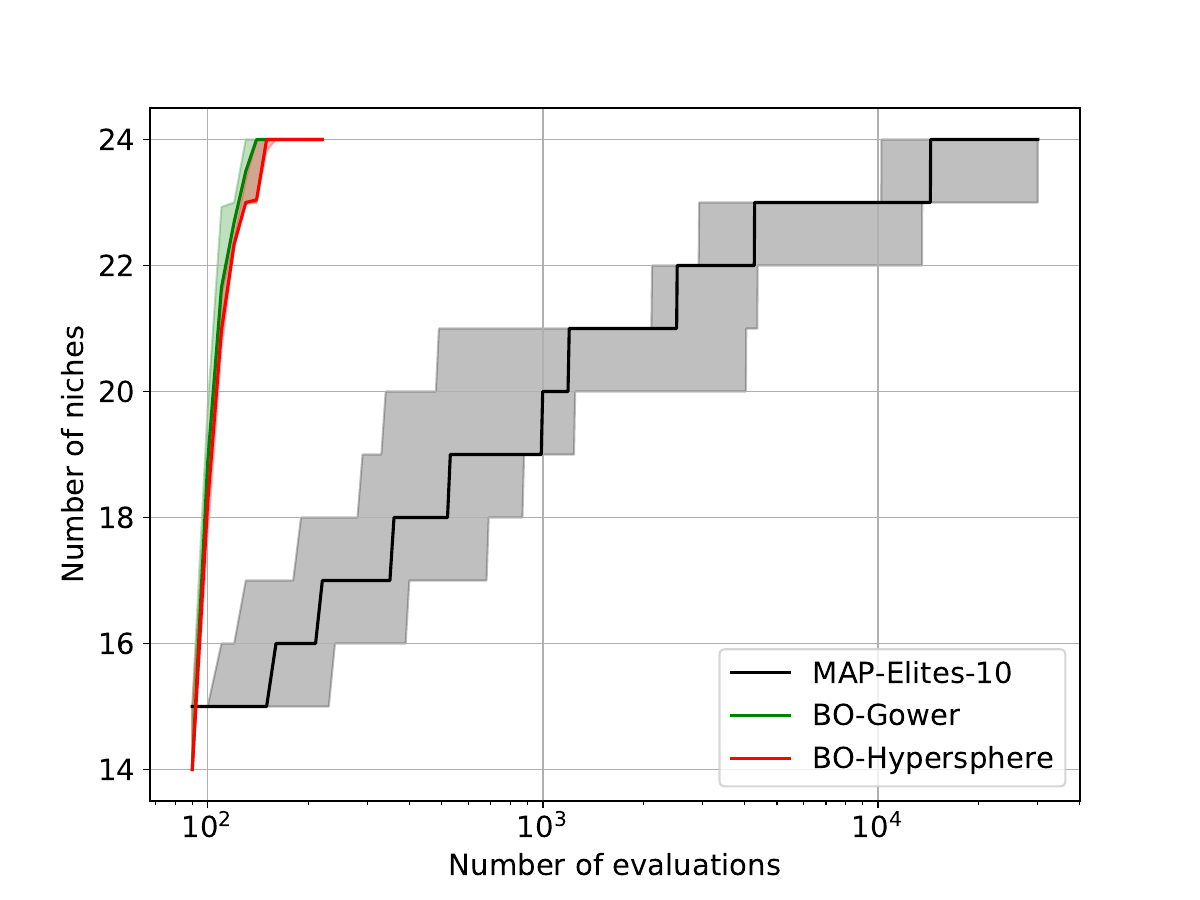}
\caption{Full convergence curves (with log-scale of the abscissa) for both Bayesian QD algorithms and MAP-Elites algorithm (with 10 individuals) for the QD-score (left) and the number of discovered niches (right)}\label{Tang_full_convergence_curve}    
\end{center}
\end{figure}

For this test problem, a study of full convergence of MAP-Elites with 10 individuals is illustrated in Figure \ref{Tang_full_convergence_curve}. It can be seen the difference with respect to the orders of magnitude in terms of the number of exact function evaluations to reach convergence with MAP-Elites ($\sim 30 000$ evaluations) compared to Bayesian QD algorithms ($\sim 220$ evaluations). Moreover, through the repetitions, it can be seen that the Bayesian QD algorithms are more robust to the initialization (with respect to the initial DoE) compared to the  MAP-Elites (with respect to the initial population).

\begin{figure}[!h]
\begin{center}
\includegraphics[width=0.49\textwidth]{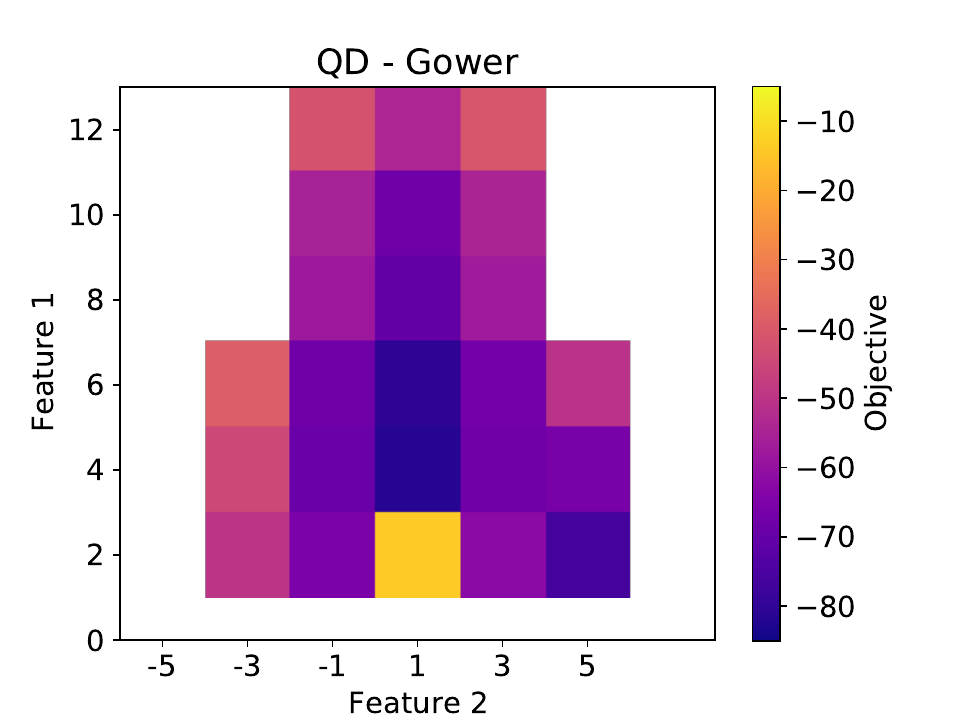}
\includegraphics[width=0.49\textwidth]{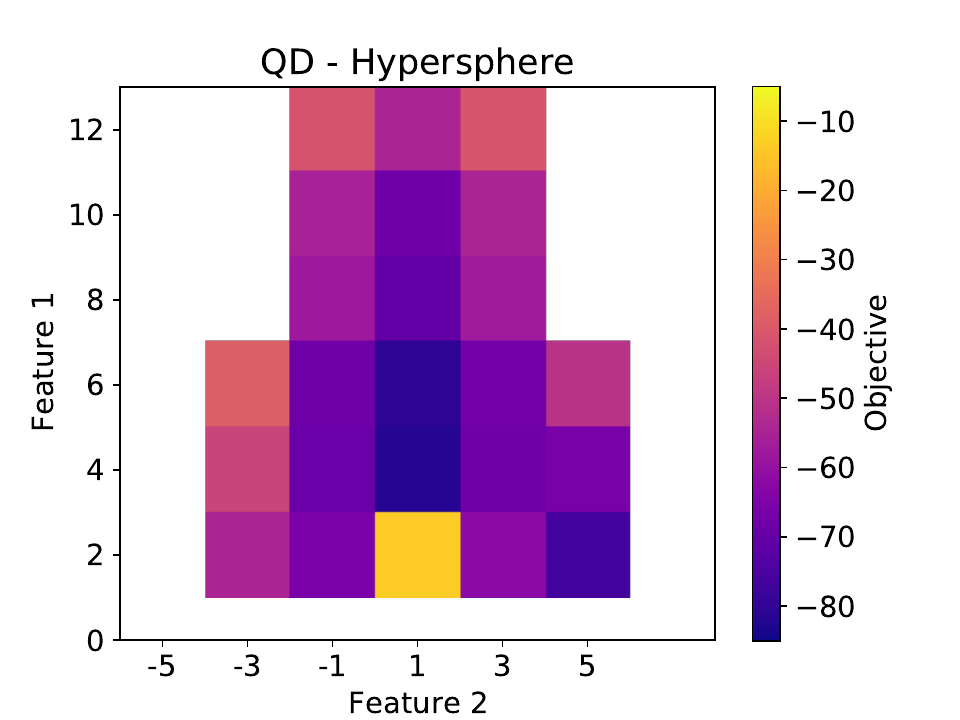}
\includegraphics[width=0.49\textwidth]{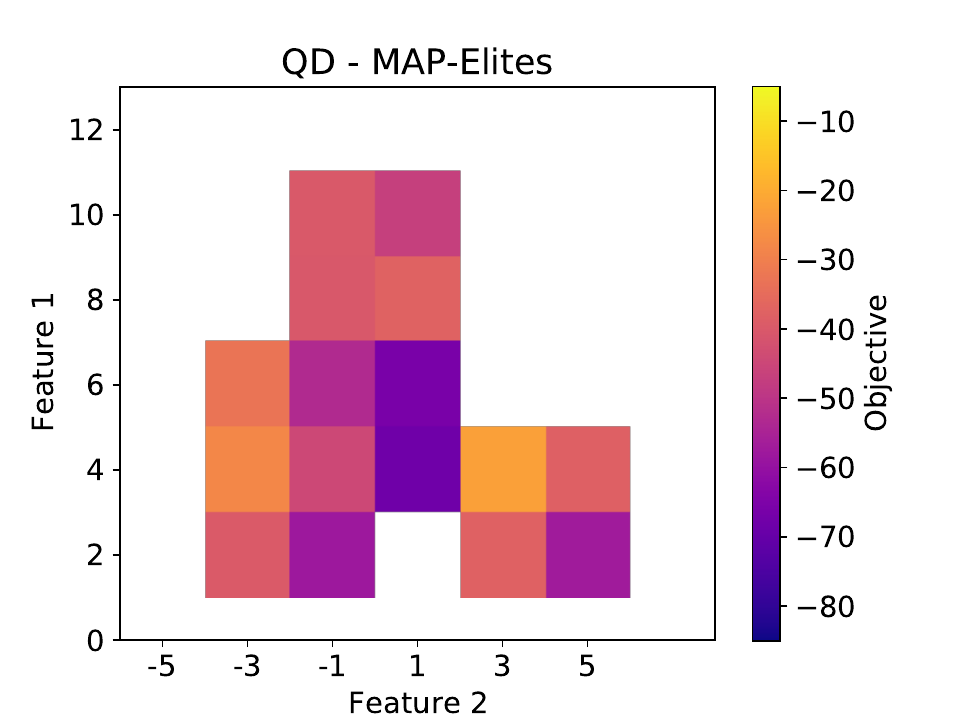}
\includegraphics[width=0.49\textwidth]{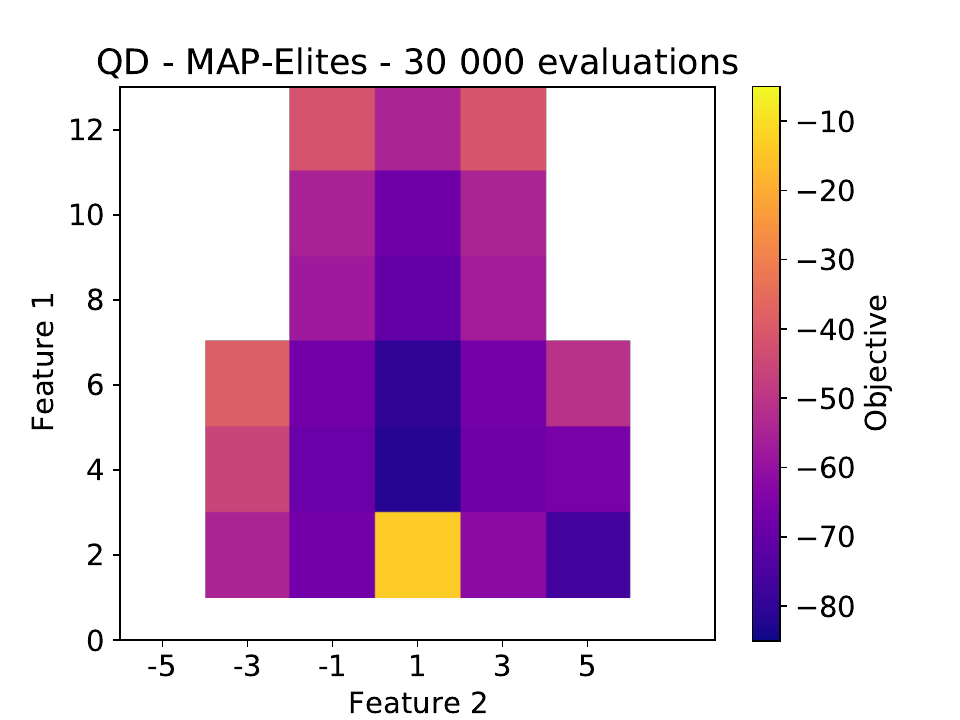}
\caption{Final archive for the Styblinski-Tang problem obtained by Bayesian QD with Gower kernel with 220 evaluations (top left), by Bayesian QD with hypersphere kernel with 220 evaluations (top right), by QD MAP-Elites with 220 evaluations (bottom left) and with QD MAP-Elites with 30 000 evaluations (bottom right) }\label{Tang_final_archive}    
\end{center}
\end{figure}

The final archive obtained by both Bayesian QD algorithms  with 220 evaluations of the exact function are identical to the final archive obtained with MAP-Elites algorithm with 30 000 evaluations to the exact function. \\

\subsection{Summary of numerical experiments on analytical test cases}
The main lessons learnt on the analytical problems are first that, in a context of computationally intensive problems, Bayesian QD algorithms converged faster than MAP-Elites approaches both in terms of QD-score and in terms of number of discovered niches. Moreover, Bayesian QD algorithms present more robust performance with respect to the initial Design of Experiments (illustrated with the problem repetitions). In addition, the choice of the kernel type of the categorical variables (between Gower and Hypersphere) is not critical with respect to the speed of convergence or the robustness to the initial DoE. However, when the number of categorical levels is large (as in the Rosenbrock problem), the hypersphere kernel might present some difficulties compared to the Gower kernel due to the number of hyperparameters to be optimized in the training of Gaussian processes. In addition, the extension of Bayesian QD algorithms to handle mixed continuous, discrete and categorical variables is interesting as it offers the possibility to have different optimal discrete/categorical values in each niche resulting in different technological or architectural choices available for the decision-makers. Eventually, independently of the problem dimension, for MAP-Elites algorithm, a smaller population size with a higher number of generations seems to be more efficient than a larger population size. 
Based on these conclusions, \red{two aerospace engineering design problems} are solved in the next sections. 

\subsection{Three dimensional wing design problem}\label{Wing_design_pb_section}

The first engineering design problem consists in optimizing the aerodynamic shape of a three dimensional wing. In the design of an aircraft, the optimization of the lifting surfaces is an essential aspect as an aerodynamically efficient wing may allow to reduce the aircraft fuel consumption. In the early design phases, the designer might be interested in a diversity of efficient wing shapes in terms of aerodynamics but that could offer different behaviors with respect to the aircraft aerodynamic center and the flight qualities for the next design steps. Indeed, in the preliminary design phases, as the main characteristics of the aircraft are not frozen, the designer might be interested in determining a collection of efficient wings depending on some trade-offs in terms of disciplinary requirements (for instance between structural, flying qualities or aerodynamics). Classical features such as the aspect ratio (the ratio between the wing span to its mean chord) or the taper ratio (the ratio between the root to the tip chord lengths of a wing) help the designer to learn about the overall aerodynamic characteristics of a wing. In such a context, Quality-Diversity algorithms may provide the designer with different high performance wing configurations in early design stages based on these two geometric features. 
 
\begin{figure}[!h]
\begin{center}
\includegraphics[width=1.\textwidth]{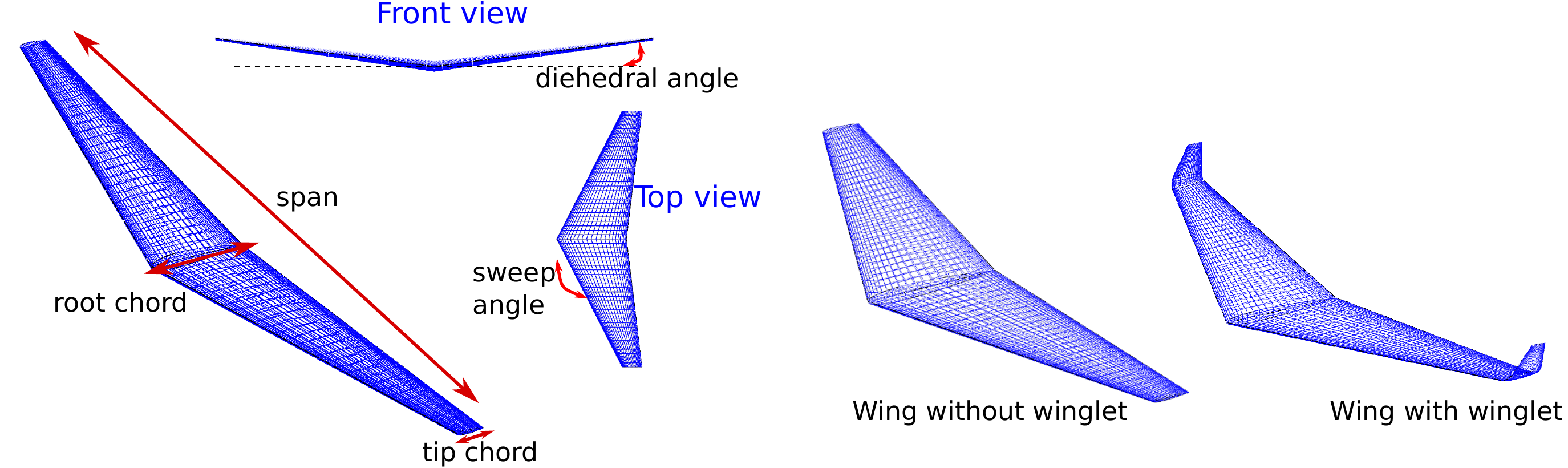}
\caption{Wing parameterization and the categorical choice corresponding to the presence or absence of winglets}\label{wing_parameters}    
\end{center}
\end{figure}

\begin{figure}[!h]
\begin{center}
\includegraphics[width=1.0\textwidth]{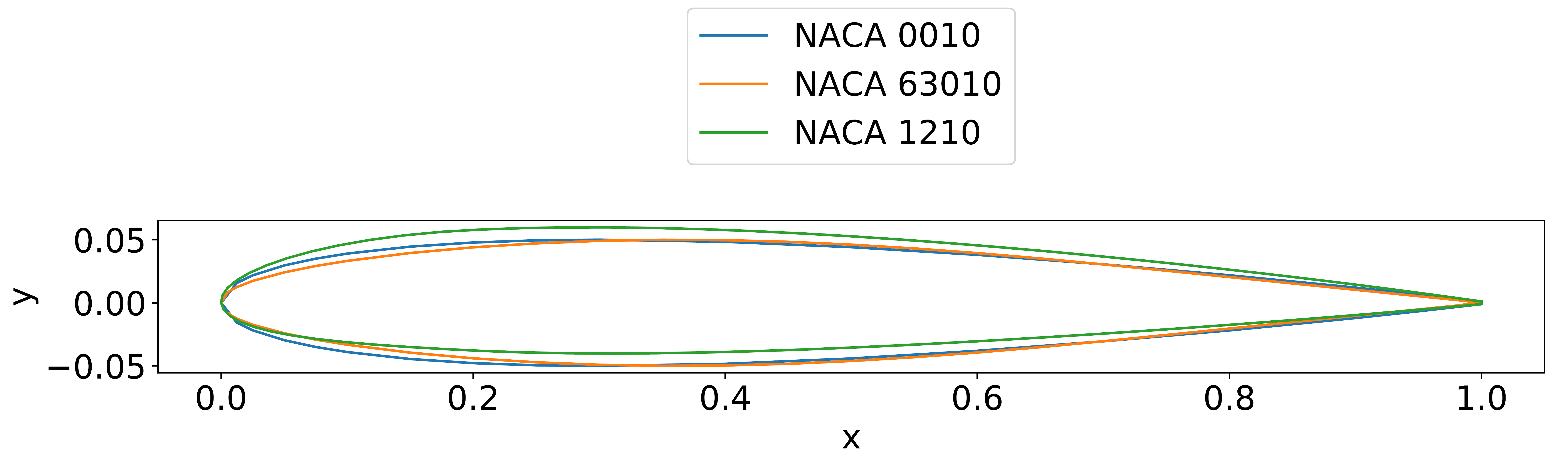}
\caption{Different choices of airfoil profiles: NACA-0010, NACA-1210 and NACA-63010}\label{Wing_airfoils}    
\end{center}
\end{figure}

In the following, the proposed Bayesian QD algorithm is applied to a wing design problem with two features: the aspect ratio and the taper ratio. The airfoil is constituted of a single section parameterized by a root chord, a tip chord, a span, a dihedral angle and a sweep angle (Figure \ref{wing_parameters}). These parameters correspond to the continuous variables of the QD problem. Two categorical variables are also involved: the presence of winglets (with two possible choices: on and off, see Figure \ref{wing_parameters}) and the airfoil type (with three possible airfoil profiles, NACA-0010 \newred{(National Advisory Committee for Aeronautics)}, NACA-1210, NACA-63010, see Figure \ref{Wing_airfoils}). The considered flight conditions are a Mach number of 0.5 and an angle-of-attack of 5 degrees.  \\

\begin{figure}[!h]
\begin{center}
\includegraphics[width=0.65\textwidth]{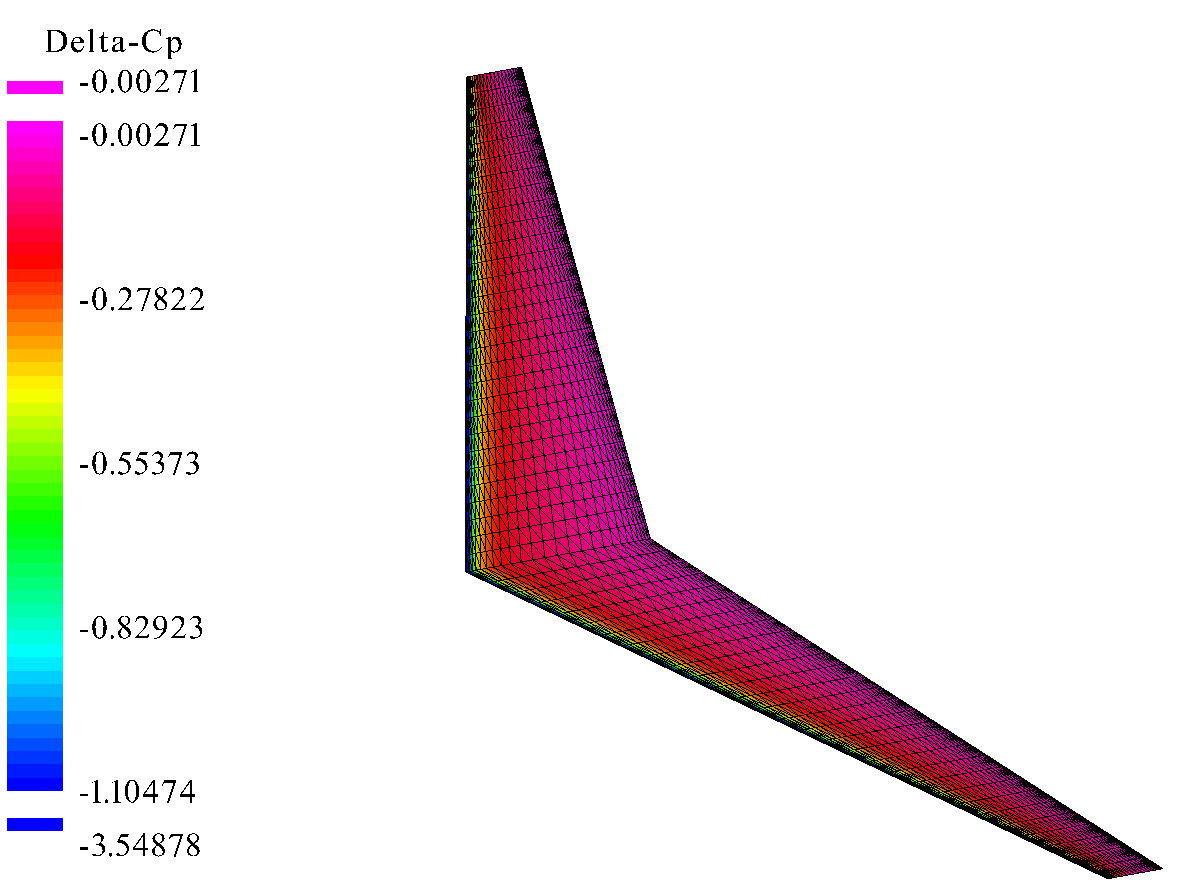}
\caption{Pressure distribution obtained by VLM calculation for one wing geometry, without winglet and with the airfoil NACA-63010}\label{Wing_QD_Gower_Cp}    
\end{center}
\end{figure}

To estimate the aerodynamic performances of the wing, a Vortex Lattice Method (VLM) is used with OpenVSP and VSPAero \cite{mcdonald2016advanced}. This consists of a simplified computational fluid dynamics model used in preliminary design phases to calculate aerodynamic forces and moments acting on the wing. The method represents lifting surfaces as an infinitely thin sheet of discrete vortices. By solving the governing equations using the Biot-Savart Law and Kutta-Joukovsky theorem \cite{anderson2011ebook}, VLM solvers can determine lift and pressure distributions, induced drag, and pitching moments endured by the wing. For the wing design problem, VLM is used to estimate the drag coefficient and the lift coefficient of a parametric wing.

The QD problem is formalized as follows:
\begin{eqnarray}\label{Eq_DQ_wing1}
\forall \tilde{\mathbf{f}} \in \mathcal{F}_t, \;\;\;  \min_{\mathbf{x}^c,\mathbf{x}^q} & & C_D(\mathbf{x}^c,\mathbf{x}^q) \\ 
\text{s.t.} & &  -C_L(\mathbf{x}^c,\mathbf{x}^q) + C_{L_T} \leq 0  \\ 
& & [f_{AR}(\mathbf{x}^c,\mathbf{x}^q),f_{TR}(\mathbf{x}^c,\mathbf{x}^q)]^T \in  \tilde{\mathbf{f}} \\ 
& & \mathbf{x}^c_\text{lb} \leq \mathbf{x}^c \leq \mathbf{x}^c_\text{ub} \\
& & \mathbf{x}^q \in \mathcal{X}^q
\label{Eq_DQ_wing2}
\end{eqnarray}
where $\mathbf{x}^c=[Rc,Tc,Sp,Di,Sw]^T$ with $Rc$ the root chord, $Tc$ the tip chord, $Sp$ the span, $Di$ the dihedral angle and $Sw$ the sweep angle. Table \ref{Cont_vb_wing} gives the domain of definition of the continuous design variables. Moreover, Figure \ref{wing_parameters} presents the parameterization of the wing with the different continuous design variables.   $C_D(\mathbf{x}^c,\mathbf{x}^q)$ is the drag coefficient of the wing, $C_L(\mathbf{x}^c,\mathbf{x}^q)$ is the lift coefficient and $C_{L_T}=0.2$ is the target lift coefficient of the wing.

\begin{table}[!h]
\begin{center}
\caption{Continuous design variables and their definition domains}\label{Cont_vb_wing}
\small
\begin{tabular}{||l c ||} 
 \hline
 Variable & Definition domain  \\ [0.5ex] 
 \hline\hline
\small
Root chord $Rc$ (m) & [3.5, 6]   \\ 
 \hline
\small
Tip chord $Tc$ (m) & [1, 1.5]   \\
 \hline
\small
Span $Sp$ (m) & [7, 11]  \\
  \hline
\small
Diehedral angle $Di$ (deg) & [0, 5] \\
  \hline
\small
Sweep angle $Sw$ (deg) & [20, 30]  \\ [1ex] 
 \hline
\end{tabular}
\end{center}
\end{table}

\normalsize
Moreover, two categorical variables are considered $\mathbf{x}^q=[Wg,Ai]^T$ with $Wg$ the presence or not of winglet and $Ai$ the choice of the airfoil profile $\{$NACA-0010, NACA-1210, NACA-63010$\}$. An example of VLM results is illustrated in Figure \ref{Wing_QD_Gower_Cp} for a wing without winglet and the airfoil NACA-63010. The results present the pressure distribution on the extrados of the wing and the associated mesh. $f_{AR}(\cdot,\cdot)$ and $f_{TR}(\cdot,\cdot)$ are the two considered feature functions corresponding respectively to the aspect ratio and the taper ratio. 

For this design problem, the use of surrogate models is relevant as one evaluation of the exact function on a cluster of 12 cores of Skylake Intel® Xeon® Gold 6152 CPU represents about 5 min. Therefore, for this QD optimization problem in dimension 7, the number of exact evaluations with algorithms such as MAP-Elites to reach convergence (in the order of several thousand evaluations) is hardly affordable. 

Considering the computational cost for this problem and the results obtained with the three analytical QD problems, only the MAP-Elites with a population of 10 individuals is considered in this test case. Moreover, for both Bayesian QD algorithms and for the MAP-Elites algorithm, five repetitions with different initial DoEs and populations are considered. 

\begin{figure}[!h]
\begin{center}
\includegraphics[width=0.75\textwidth]{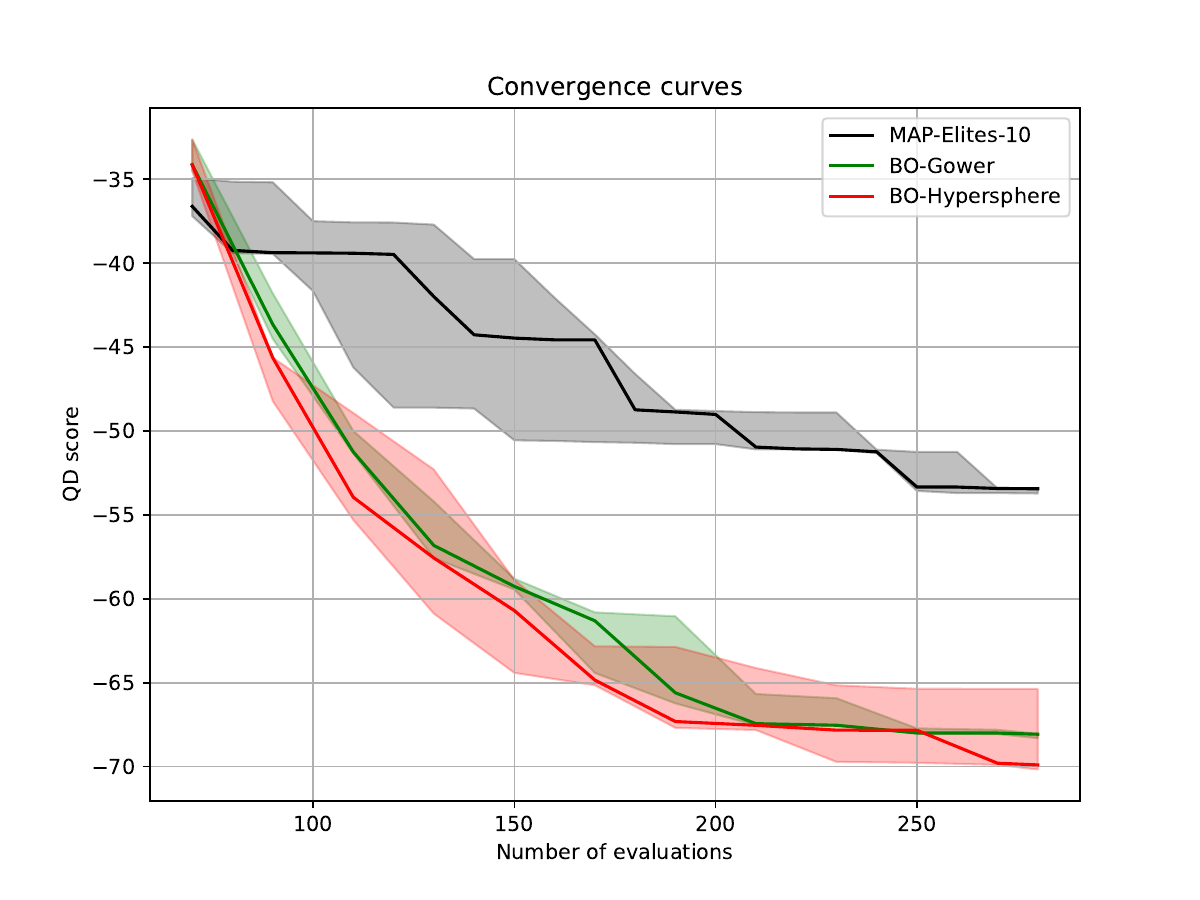}
\caption{Convergence curves (normalized QD score, the lower, the better) for the Wing design problem with MAP-Elites and Bayesian QD algorithm with Gower and hypersphere kernels. For the five repetitions, the curves correspond to the median whereas the upper and lower limits of the shade area corresponds to the $75^\text{th}$ and $25^\text{th}$ quantiles.}\label{Wing_Convergence}    
\end{center}
\end{figure}

\begin{figure}[!h]
\begin{center}
\includegraphics[width=0.75\textwidth]{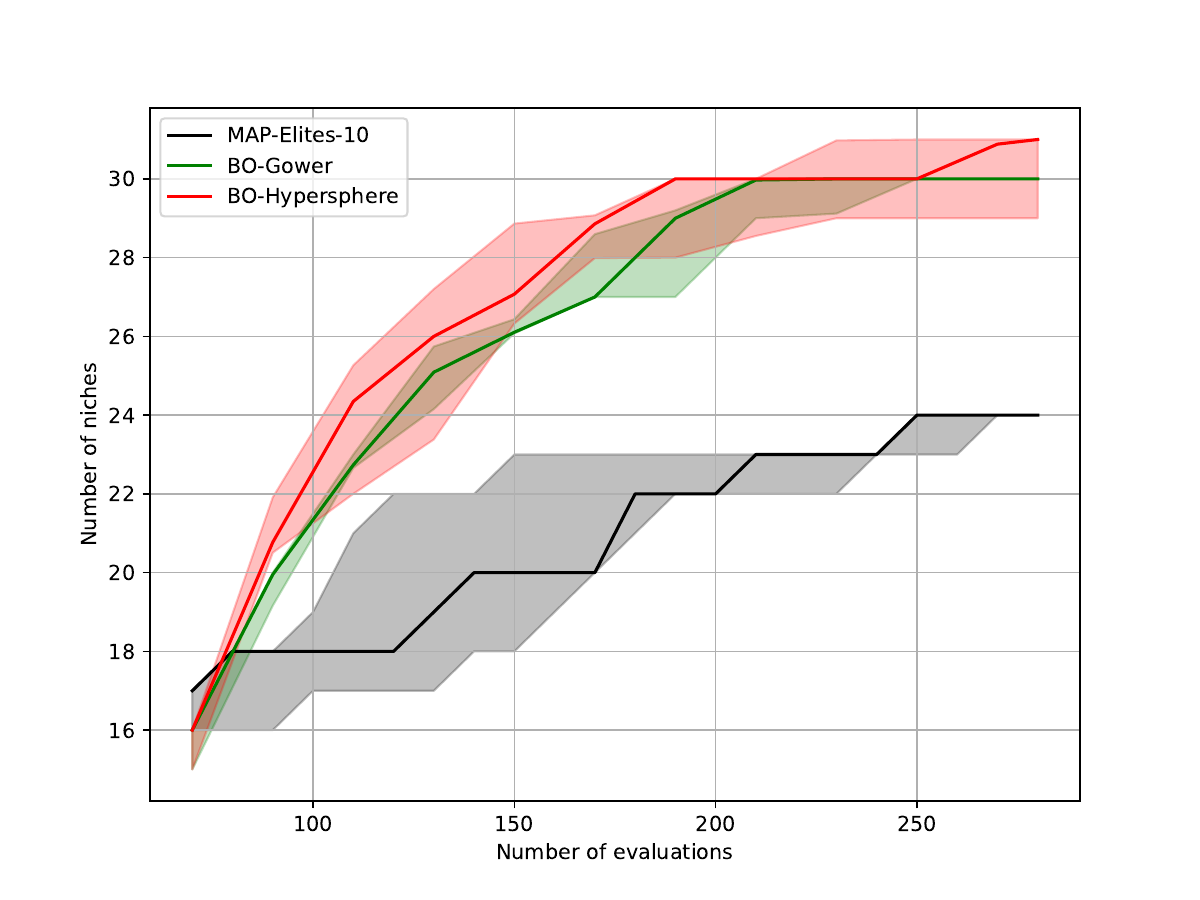}
\caption{Number of discovered niches (the higher, the better) for the Wing design problem with MAP-Elites and Bayesian QD algorithm with Gower and hypersphere kernels. For the five repetitions, the curves correspond to the median whereas the upper and lower limits of the shade corresponds to the $75^\text{th}$ and $25^\text{th}$ quantiles. }\label{Wing_Convergence_niche}  
\end{center}
\end{figure}

Figures \ref{Wing_Convergence} and \ref{Wing_Convergence_niche} present the convergence curves for the QD-score and the number of discovered niches in the wing design problem.  As for the analytical test problems, for the wing design problem, the Bayesian QD algorithms converge faster to a better solution in terms of QD-score and illuminate a larger number of niches compared to the MAP-Elites algorithm. Moreover, for this engineering design problem, the choice of the covariance model to deal with the discrete and categorical variables has a limited influence, both the Gower and the hypersphere kernels provide similar results. 

\begin{figure}[!h]
\begin{center}
\includegraphics[width=0.49\textwidth]{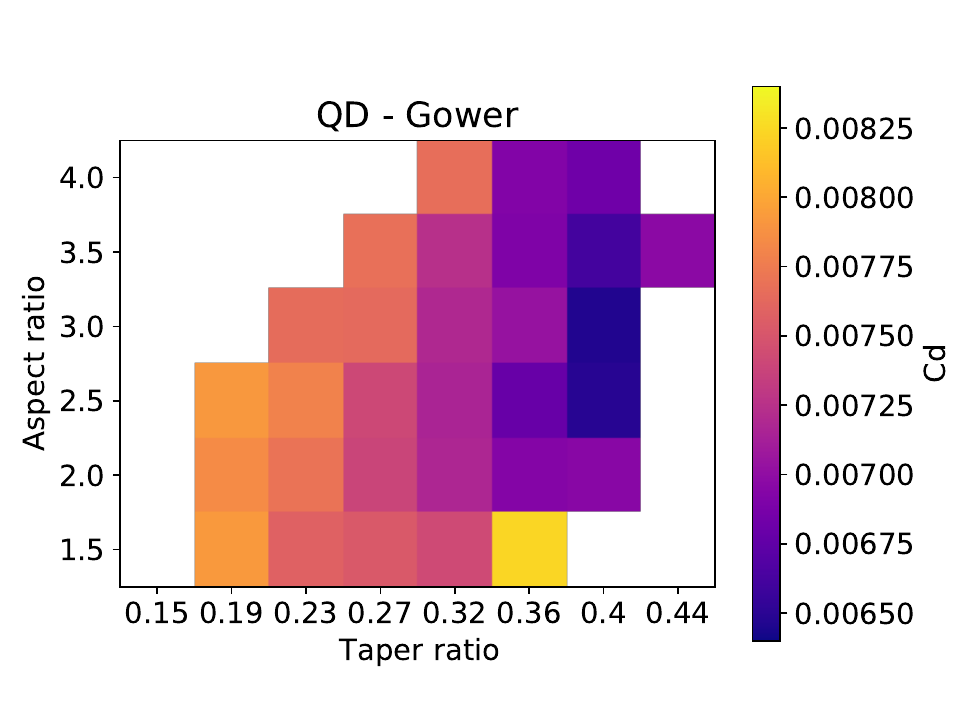}
\includegraphics[width=0.49\textwidth]{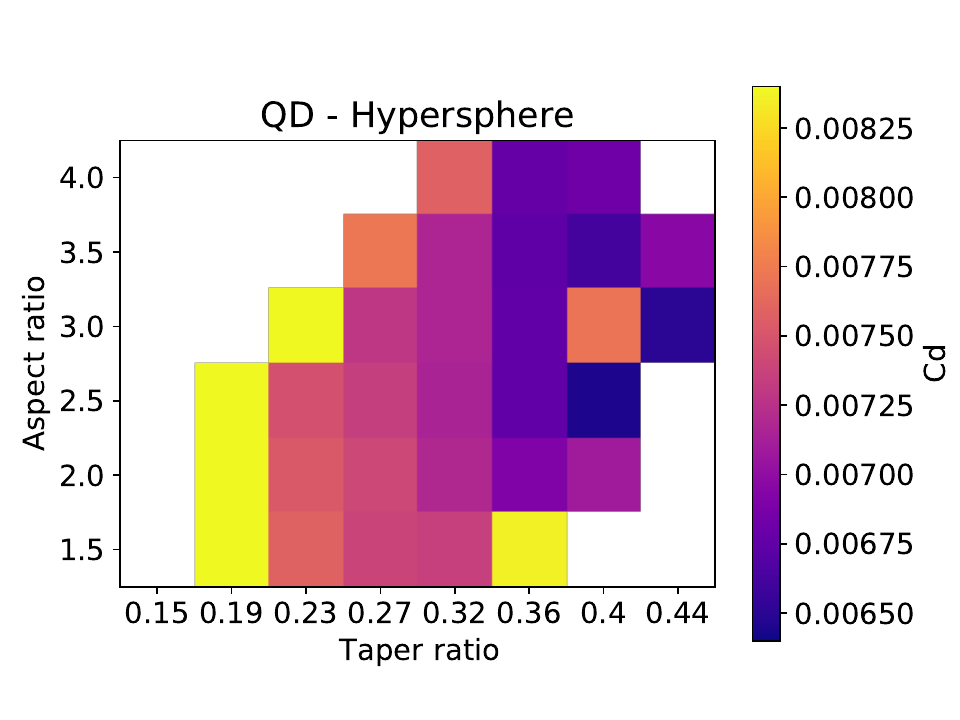}
\includegraphics[width=0.49\textwidth]{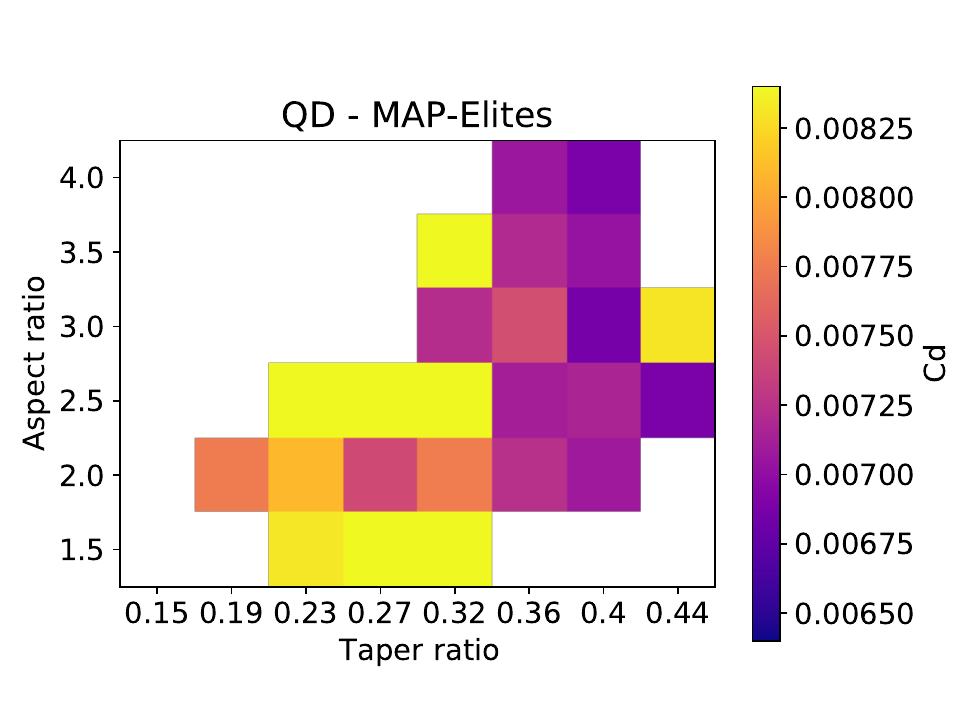}
\caption{Final archive for the Wing design problem obtained by Bayesian QD with Gower kernel (top left), by Bayesian QD with hypersphere kernel with (top right), by MAP-Elites (bottom left)}\label{Wing_final_archive}    
\end{center}
\end{figure}

\red{The final archive obtained (with one of the repetitions) for the wing design problem is presented in Figure \ref{Wing_final_archive}, for the Bayesian QD with the Gower and the hypersphere kernels and with MAP-Elites. The archives for both Bayesian QD algorithms are similar with a small advantage to the Gower kernel with lower objective values in some niches. It highlights the robustness of the proposed approach to different choices of kernels to take into account mixed continuous, discrete and categorical variables in Gaussian process.} However, as for the analytical test problems, the archive obtained by MAP-Elites does not illuminate as many niches as the Bayesian QD algorithms and the quality in each discovered niche is less optimal.

\begin{figure}[!h]
\begin{center}
\includegraphics[width=0.49\textwidth]{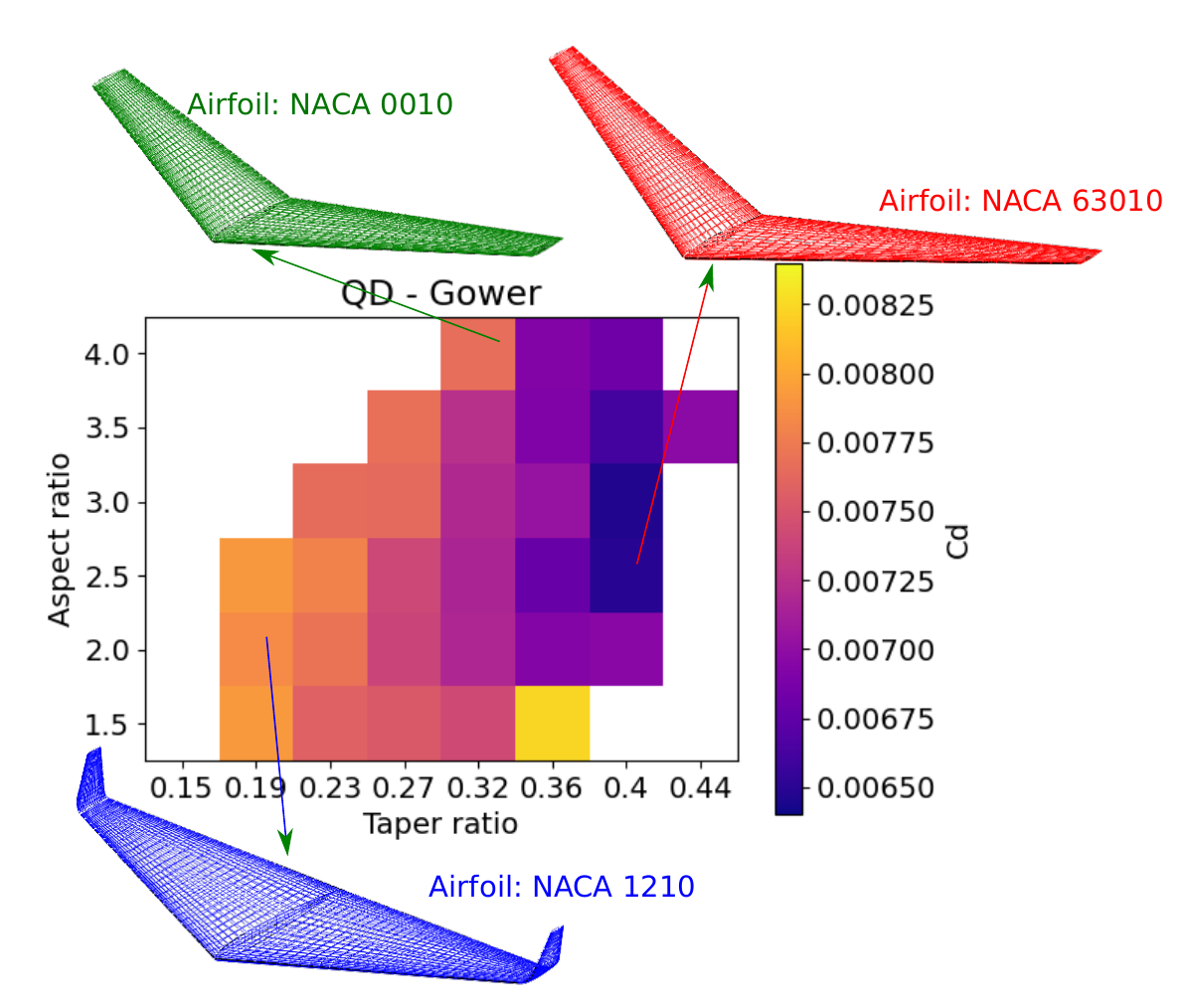}
\includegraphics[width=0.49\textwidth]{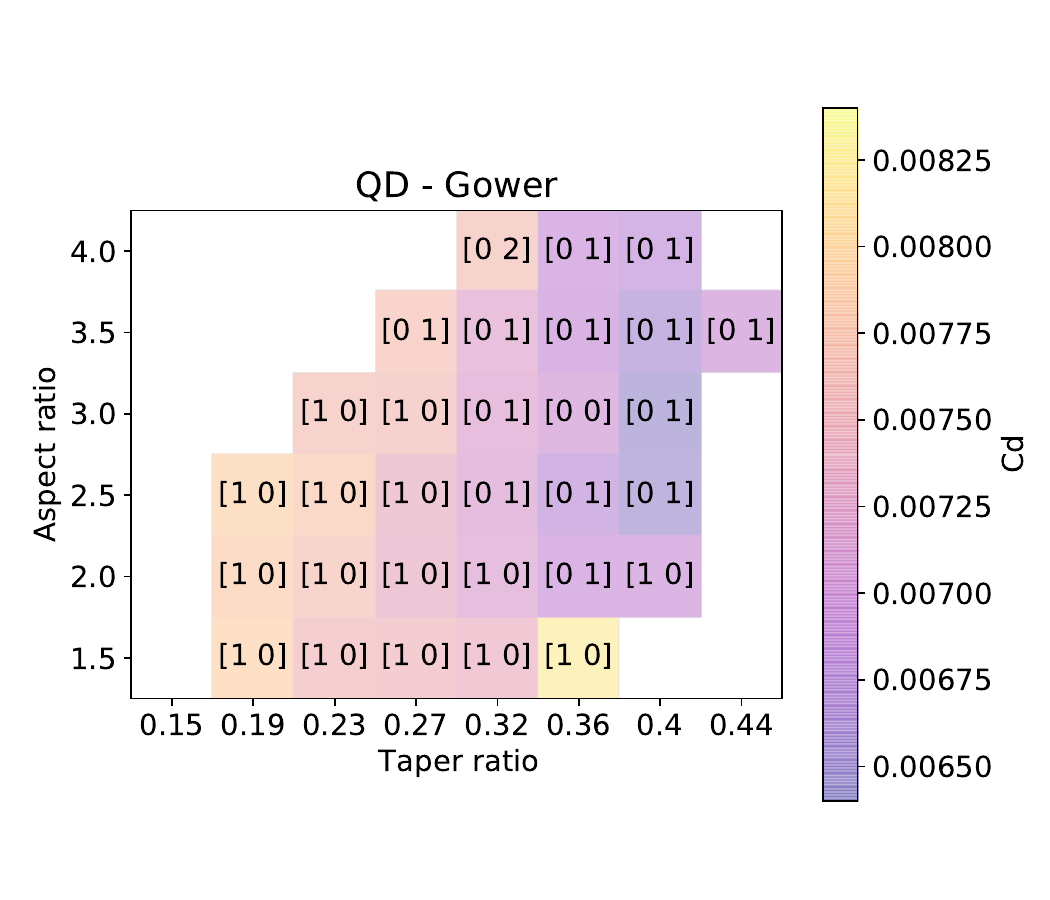}
\caption{Comparison of wing geometries extracted from the final archive obtained by one repetition with Bayesian QD with Gower kernel and with the found categories $\mathbf{x}^q=[W_g,A_i]\in\{0,1\} \times \{0,1,2\}$ in each niche. $W_g=0$ means absence of winglet and $W_g=1$ the presence of winglets. $A_i$ corresponds to the choice of airfoil with 0: NACA-0010, 1: NACA-63010 and 2: NACA-1210.}\label{Wing_QD_Gower_illustration}    
\end{center}
\end{figure}

Figure \ref{Wing_QD_Gower_illustration} illustrates the comparison of wing geometries extracted from the final archive obtained by one repetition with Bayesian QD with Gower kernel. Moreover, it presents the found categories $\mathbf{x}^q=[W_g, A_i]^T$ in each niche illustrating the convergence to different optimal solutions in terms of architecture choices (presence or absence of winglets and airfoil type). Due to the features combination (aspect ratio and taper ratio) different categorical choices are identified in the niches. The handling of mixed continuous, discrete and categorical variables offers a diversity of architectures and technology choices to the designers and the decision-makers. Moreover, the feature functions offer a diversity in terms of wing geometry. It is possible to identify thanks to the final archive the consequences of an increase of aspect ratio in terms of drag coefficient for the wing. \red{The final archive is a valuable asset for the design to make trade-offs in early design phases and to balance the consequences of these trade-offs. The discovery of new niches with a high quality (an optimal fitness value) allows to provide new choices for the decision-makers with respect to interesting characteristics that are described by the features. Moreover, the consequences in terms of fitness (objective function) for these new choices are directly available to the decision-markers via the optimal QD archive.}
\red{In the next section, the final aerospace test problem is presented increasing the complexity in terms of problem dimension and number of combination possibilities for the discrete and categorical variables.}

\subsection{Multidisciplinary design of a two-stage solid sounding rocket problem}\label{Rocket_section}

\red{Sounding rocket are very precious aerospace vehicles in order to carry out scientific experiments at a low cost \cite{seibert2006history}. A large number of sounding rockets are currently operational and the launch rate of such vehicles is quite large (several hundred per year \cite{NASA2022,mcdowell2023space}). These launch vehicles can be composed of multiple stages in order to reach high altitude with heavy payload. Numerous architectural and technological choices can be made when considering a solid propellant rocket. Some of these technological choices are related to the type of solid propellant, the geometry of the grain, the type of rocket material case or the number of stages.}

\begin{figure}[h!]
\begin{center}
    \includegraphics[angle=90,width=0.75\linewidth]{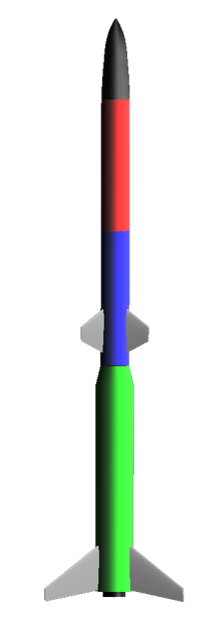}
        \caption{Illustrative view of sounding rocket with in green: the first stage, in blue: the second stage, in red: the scientific experiments (payload) and in black: the fairing. The lift surfaces are required for the guidance, control and navigation of the vehicle. }
        \label{fig:representation_sounding}
\end{center}
\end{figure}

\begin{figure}[h!]
\begin{center}
    \includegraphics[width=0.75\linewidth]{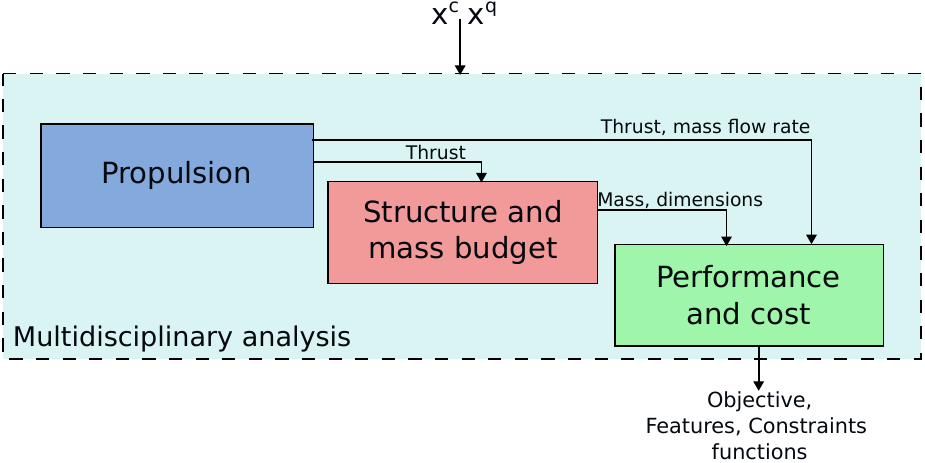}
        \caption{Design process of sounding rocket using a multidisciplinary analysis coupling different disciplines (propulsion, structure and mass budget, performance and cost)}
            \label{fig:design_process_sounding}

\end{center}
\end{figure}

\red{This test problem involves the design and the optimization of a two-stage sounding rocket. The two stages are powered using solid propellant. The objective function, features and constraints result from complex multidisciplinary coupled simulations involving propulsion, structure and performance evaluations (Figure \ref{fig:design_process_sounding}). The goal of the design process is to minimize the cost of the rocket (including development, production and launch costs) while considering two features: the payload mass and the altitude at the apogee of the trajectory. The design optimization problem of each sounding rocket stage involves 4 design variables (ratio between throat diameter and nozzle exit diameter, propellant mass, combustion pressure, nozzle exit diameter) and 3 discrete variables (type of engine, type of material of the casing, type of propellant). The payload mass is also considered as a design variable. Therefore, the optimization problem is in dimension 15. The dimension of the continuous design space is 9 and the dimension of the discrete design space is 6. Concerning the discrete variables, 4 types of solid propellant are considered (butalite and butalane - composed of Hydroxyl-terminated polybutadiene and ammonium perchlorate, and in addition aluminum for butalane ; nitramite - nitrocellulose/nitroglycerine ; and pAIM-120) \cite{north1988design}, 3 types of material (aluminium, steel and composite) for the rocket structure and 3 types of engines (Table \ref{combinatory_LV}). The engines differ by their geometry and efficiency. A representative view of the two-stage sounding rocket is displayed in Figure \ref{fig:representation_sounding}.}

\begin{table}[!h]
\begin{center}
\caption{Categorical design variables and their definition domains}\label{combinatory_LV}
\small
\begin{tabular}{||l c ||} 
 \hline
 Variable & Definition domain  \\ [0.5ex] 
 \hline\hline
\small
Propellant stage 1 & \{Butalane, Butalite, Nitramite, p-AIM120\}  \\ 
 \hline
\small
Casing material stage 1 & \{Steel, Aluminum, Composite\}  \\
 \hline
\small
Engine type stage 1 & \{Type 1, Type 2, Type 3\}  \\
\hline
  \hline
\small
\small
Propellant stage 2 & \{Butalane, Butalite, Nitramite, p-AIM120\} \\ 
 \hline
\small
Casing material stage 2 & \{Steel, Aluminum, Composite\}   \\
 \hline
\small
Engine type stage 2 & \{Type 1, Type 2, Type 3\}  \\
 \hline
 \hline
 Number of possible combinations&1296\\
 \hline
\end{tabular}
\end{center}
\end{table}

\red{The combination of the possible categorical architectures for this test-case is equal to 1296. The design process involves 3 disciplines: propulsion, structure and performance assessment. The design process is summarized in Figure \ref{fig:design_process_sounding}.}\\

The associated QD problem is the following:

\begin{eqnarray}\label{Eq_DQ_wing1}
\forall \tilde{\mathbf{f}} \in \mathcal{F}_t, \;\;\;  \min_{\mathbf{x}^c,\mathbf{x}^q} & & C(\mathbf{x}^c,\mathbf{x}^q) \\ 
\text{s.t.} & &  g_i(\mathbf{x}^c,\mathbf{x}^q) \leq 0\;\; i=1,\dots,8  \\ 
& & [f_{m_{CU}}(\mathbf{x}^c,\mathbf{x}^q),f_{alt}(\mathbf{x}^c,\mathbf{x}^q)]^T \in  \tilde{\mathbf{f}} \\ 
& & \mathbf{x}^c_\text{lb} \leq \mathbf{x}^c \leq \mathbf{x}^c_\text{ub} \\
& & \mathbf{x}^q \in \mathcal{X}^q
\label{Eq_RLV2}
\end{eqnarray}
\red{where $C(\cdot,\cdot)$ is the normalized cost of the vehicle, $\mathbf{x}^c=[R_1,Mp_1,Pc_1,Ne_1,$
$R_2,Mp_2,Pc_2,Ne_2,m_{CU}]^T$ with $R_i$ the ratio between the throat diameter and the nozzle exit diameter (for $i\in\{1,2\}$), $Mp_i$ the propellant mass, $Pc_i$ the combustion pressure, $Ne_i$ the nozzle exit diameter, $m_{CU}$ the payload mass and $i$ the stage number. The vector of categorical variables is $\mathbf{x}^q=[e_1,m_1,p_1,e_2,m_2,p_2]^T$ with $e_i$ the type of engine, $m_i$ the type of material of the casing and $p_i$ the type of propellant. Eight inequality constraints are also present in the optimization problem. They refer to specifications related to the thrust to weight ratios at the beginning of the flight phase (one constraint per stage), the pressure consistency in the engine (two constraints) and the geometrical layout of the solid propellant bloc (four constraints). The two features correspond to the payload mass $m_{CU}$ and the apogee altitude $alt$. These two feature functions are of interest for decision-makers as in early design phase the target mission might not be fully defined and therefore different choices in terms of payload mass or in terms of target culmination altitude could be made. The interest of QD approach is to provide  the decision-maker with the consequences in terms of cost (objective function) for different choices of mission target (payload mass and culmination altitude) and also in terms of sounding rocket architecture. }

\red{Considering the computational cost for this problem and the results obtained with the previous QD problems, only the MAP-Elites with a population of 10 individuals is considered. Five repetitions with different initial DoEs and populations are carried out.}\\

\begin{figure}[!h]
\begin{center}
\includegraphics[width=0.75\textwidth]{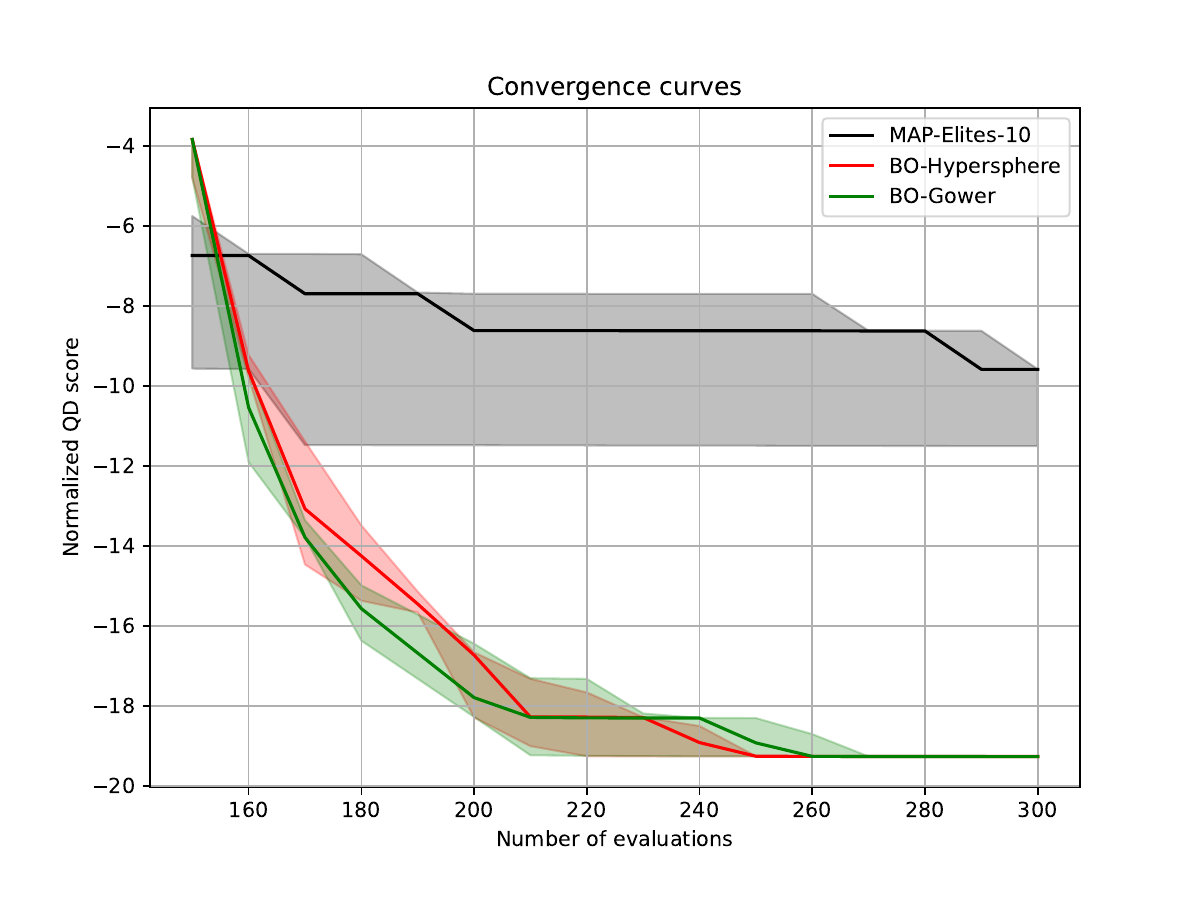}
\caption{Convergence curves (normalized QD score, the lower, the better) for the sounding rocket design problem with MAP-Elites and Bayesian QD algorithm with Gower and hypersphere kernels. For the five repetitions, the curves correspond to the median whereas the upper and lower limits of the shade area corresponds to the $75^\text{th}$ and $25^\text{th}$ quantiles.}\label{LV_Convergence}    
\end{center}
\end{figure}

\begin{figure}[!h]
\begin{center}
\includegraphics[width=0.75\textwidth]{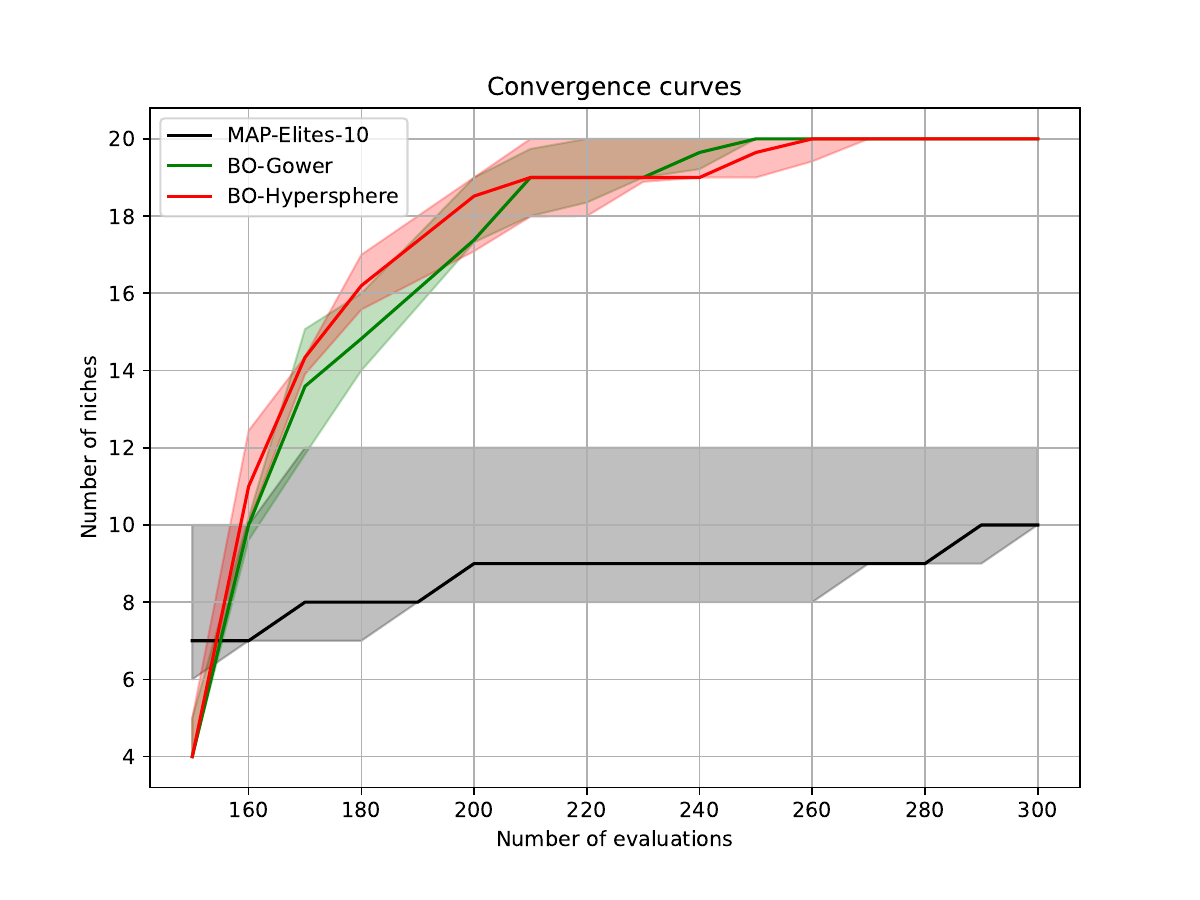}
\caption{Number of discovered niches (the higher, the better) for the sounding rocket design problem with MAP-Elites and Bayesian QD algorithm with Gower and hypersphere kernels. For the five repetitions, the curves correspond to the median whereas the upper and lower limits of the shade corresponds to the $75^\text{th}$ and $25^\text{th}$ quantiles. }\label{LV_Convergence_niche}  
\end{center}
\end{figure}

\red{The convergence curves of the three considered algorithms are displayed in Figure \ref{LV_Convergence} for the QD score and Figure \ref{LV_Convergence_niche} for the number of illuminated niches. The two versions of Bayesian algorithms overcome the MAP-Elites algorithm in terms convergence speed and final results for the QD-score and number of discovered niches. Indeed, both proposed algorithms tend to discover about 20 niches whereas the MAP-Elites algorithm illuminates only 10 niches. Furthermore, the convergence of the Bayesian algorithms is more robust to the initialization than with the MAP-Elites. Even if the possible combinations of discrete and categorical architectures is large on this test case (1296), the two versions of the proposed algorithms succeed to illuminate all the accessible niches.} \\

\begin{figure}[!h]
\begin{center}
\includegraphics[width=0.49\textwidth]{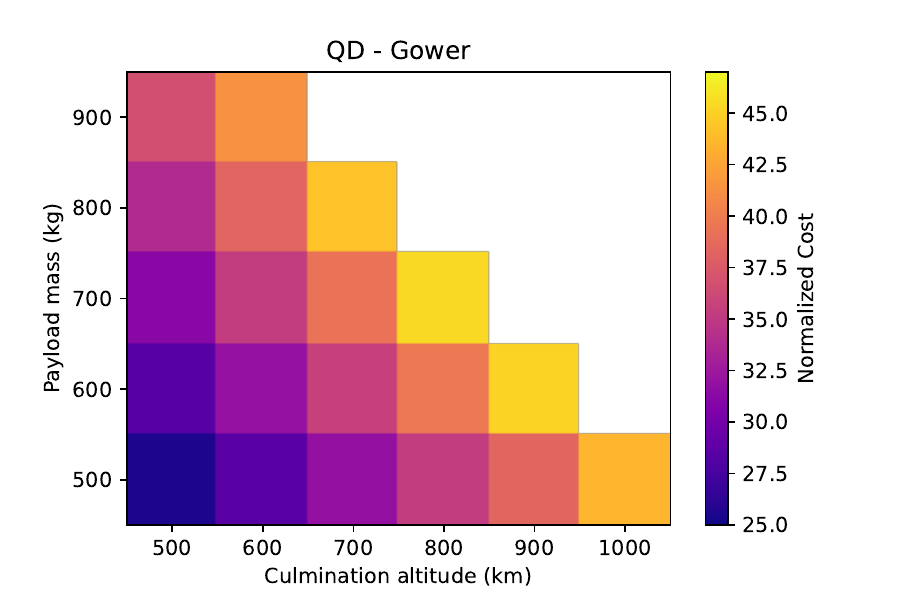}
\includegraphics[width=0.49\textwidth]{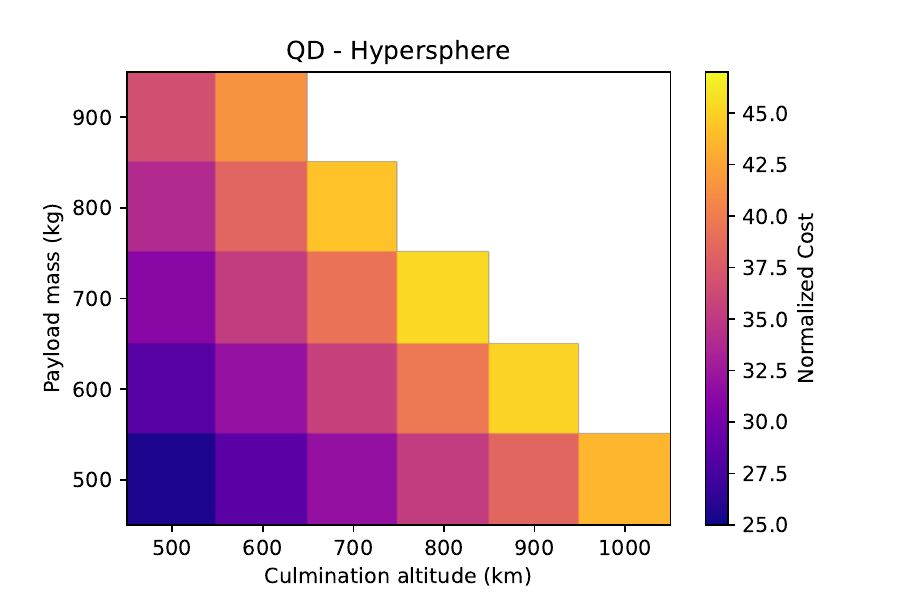}
\includegraphics[width=0.49\textwidth]{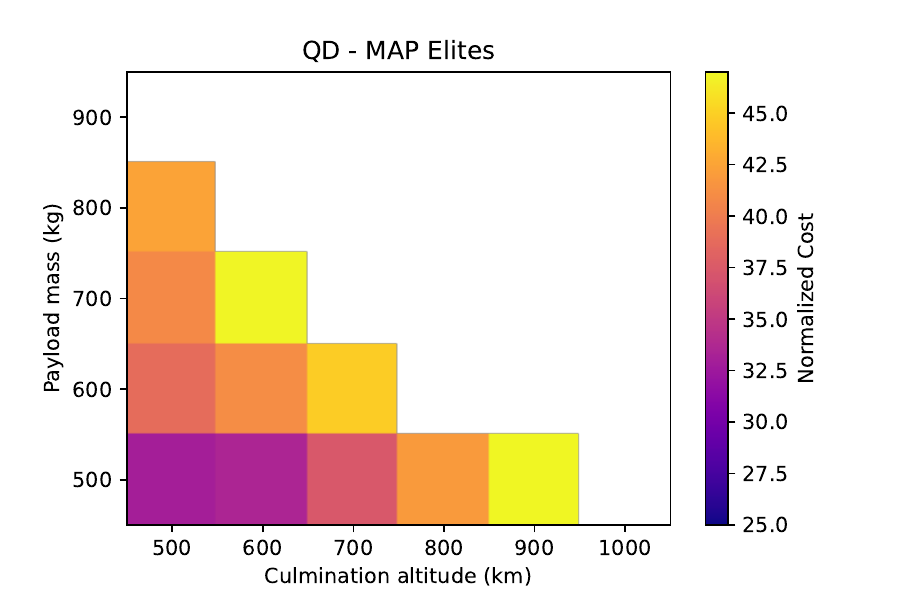}
\caption{Final archive for the sounding rocket design problem obtained by Bayesian QD with Gower kernel (top left), by Bayesian QD with hypersphere kernel (top right), by MAP-Elites (bottom)}\label{LV_final_archive}    
\end{center}
\end{figure}

\red{Figure \ref{LV_final_archive} represents the final illuminated maps for the three considered algorithms after 300 evaluations of the exact multidisciplinary process. As it can be seen, the maps provided by the two proposed algorithms are very similar on this test case. On the contrary, the classical MAP-Elites algorithm does not achieve to illuminate all the niches and the solutions in the illuminated niches are less efficient in terms of objective function.}

\begin{figure}[h!]
\begin{center}
    \includegraphics[width=1\linewidth]{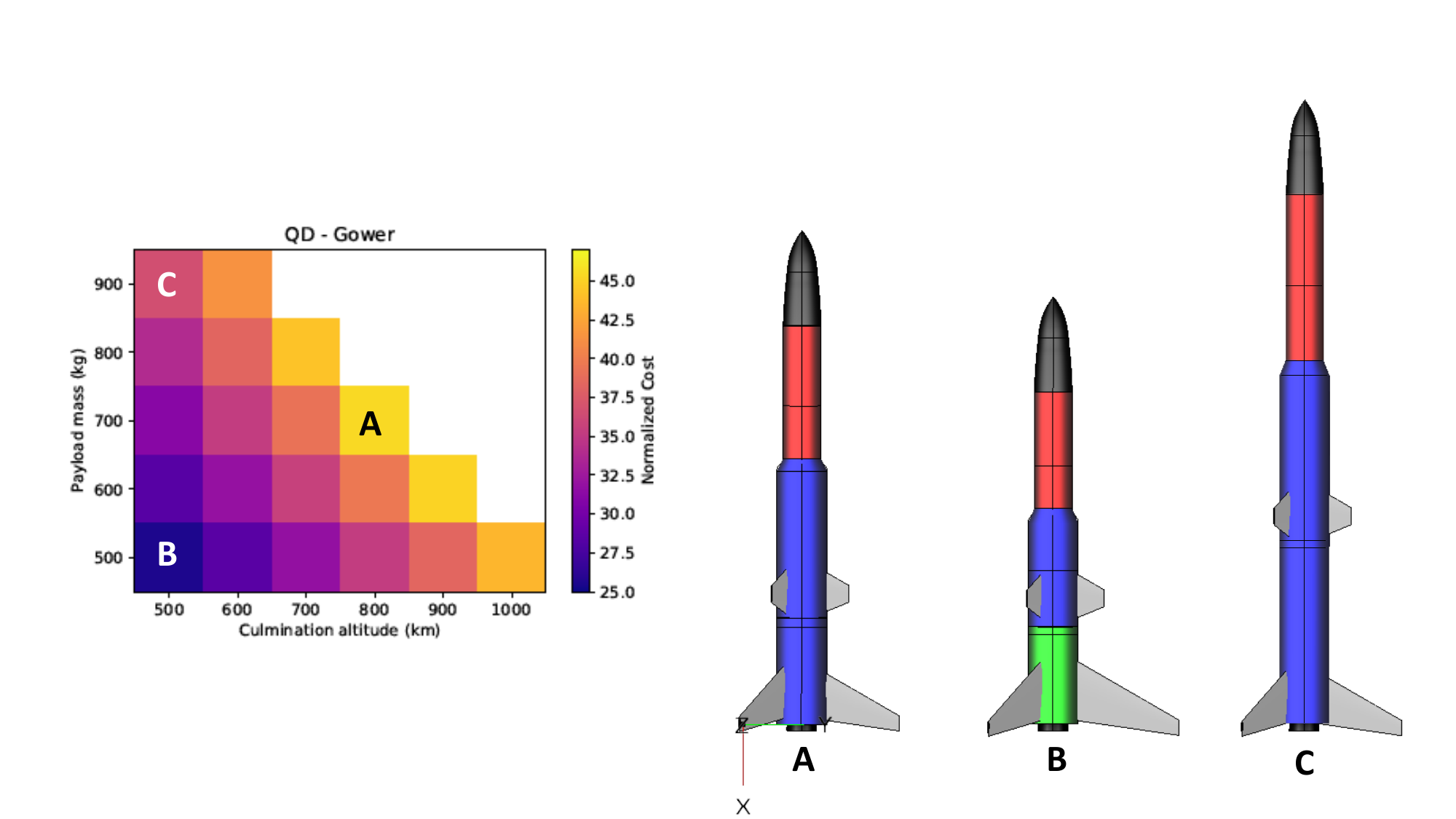}
        \caption{Sounding rocket architectures obtained for different niches of the BO-Gower map, with in black: the fairing, in red: the scientific experiment (payload), in blue: the stage powered using p-AIM120 propellant, and in green: the stage powered with butalane propellant.}
            \label{fig:rlv_map}
\end{center}
\end{figure}

\red{Figure \ref{fig:rlv_map} illustrates some different architectures obtained using the BO-Gower algorithm. It can be seen that the categorical variables, especially the propellant type, vary from one niche to another. This variation is particularly important on the first stage. With this map, one can choose between different payload masses and target culmination altitudes depending on the scientific experiment to carry out. As it can be expected, the mission cost increases with respect to both the mass of the payload and the culmination altitude to be reached. This final map summarizes a set of valuable information for the decision-makers as it is possible to get the direct consequences on the rocket cost of the choice of another target mission. The choice with respect to the features can therefore be made by the decision-makers in relation with the scientific experiment needs and the available budget. It offers a higher-level of exploration of the design space with respect to the features compared to classical optimization approaches. In addition, for this particular test problem, the objective function and the two features are not antagonistic as the sounding rocket with the lower cost corresponds to the sounding rocket that can reach the lower culmination altitude with the smaller payload mass. Therefore, multi-objective optimization algorithms could not offer the same diversity of solutions compared to QD approaches. }

\begin{figure}[h!]
\begin{center}
    \includegraphics[width=0.75\linewidth]{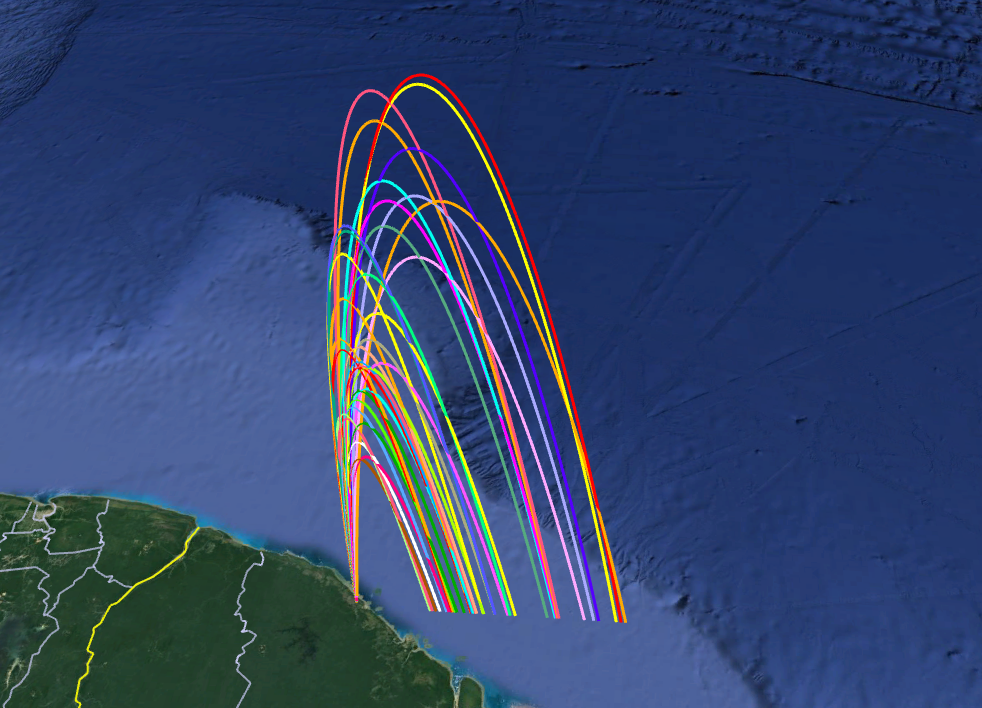}
        \caption{Obtained trajectories on the different niches obtained for one repetition of Bayesian QD with Gower kernel.}
            \label{fig:traj}

\end{center}
\end{figure}

\red{The obtained trajectories for the different niches are displayed in Figure \ref{fig:traj}. It can be seen that the culmination altitude is different between the niches which is of major interest depending on the type of scientific experiments that is carried out. }

\newpage
\section{Conclusions}
\label{sec:conclusion}

In this paper, a Bayesian Quality-Diversity (QD) optimization approach is proposed to solve constrained problem involving mixed continuous, discrete and categorical variables. \red{The proposed algorithm extends the existing Bayesian QD algorithms in order to offer the possibility to handle architectural and technological choices (through discrete and categorical variables) in engineering design problems while accounting for specification constraints.} Quality-Diversity approaches allow to identify a diversity of solutions with respect to some features with a high potential with respect to the objective. The proposed Bayesian QD algorithm relies on the use of Gaussian processes (GPs) to replace the exact objective, features and constraints functions. Moreover, in order to account for the presence of mixed continuous, discrete and categorical variables, two covariance models based on the Gower distance and the hypersphere decomposition are used for the GPs. The proposed Bayesian QD algorithm has been tested on a series of analytical QD problems of increasing complexity and on \red{two} engineering problems corresponding to the design of an aircraft wing \red{and the multidisciplinary design of a two-stage sounding rocket. These numerical experiments allowed to assess the efficiency of the proposed Bayesian QD algorithms in terms of convergence speed and robustness to the initialization compared to MAP-Elites algorithm. }

\red{Some limitations might be pointed out for the proposed approach. Firstly, Gaussian process may face difficulties for high dimensional  optimization problems due to the tuning of a large number of hyperparameters associated to the covariance model. Secondly, the proposed algorithm has been numerically evaluated on a limited number of feature functions and the full exploration of the Quality-Diversity map could be more difficult with a higher number of features.} 

\red{In future works, an adaptation of the proposed approach to high dimensional optimization problems could be investigated in order to overcome the dimension limitation of classical Gaussian process and Bayesian optimization.}

\section*{Acknowledgments}
This work is part of the PHOBOS project funded by ONERA - The French Aerospace. The authors want to thank Dr. R. Wuilbercq and M. G. Sire for fruitful discussions.

\newpage
\appendix

\section{Numerical settings}

\subsection{Numerical settings for the involved algorithms}\label{Appendix_num_settings}
\begin{itemize}
\item MAP-Elites:
\begin{itemize}
    \item MAP-Elites population size: 10 individuals or $10\times (d_c+d_q)$ individuals with $d_c$ the number of continuous variables and $d_q$ the number of categorical variables. 
    \item Mutation : probability of 0.4 and mutation according to a Normal distribution centered on 0 and with standard deviation of 0.3 for normalized continuous variables in $[0,1]^{d_c}$.
    \item \red{MAP-Elites generations: in the solving of the QD problem associated to the infill criterion, 4000 generations are carried out. The computational cost for solving the QD problem and to obtain the final archive is in the order of 20 seconds.}
\end{itemize}
\item Gaussian processes:
\begin{itemize}
    \item Initial LHS size for Gaussian process: $10\times (d_c+d_q)$ using pyDOE2 library \cite{sjogren2018pydoe2};
    \item Gaussian process library: SMT \cite{bouhlel2019python};
    \item \red{Gaussian process training algorithm: COBYLA with 20 multi-starts. The GP hyperparameter training cost is limited to few seconds (below 5s) as the number of training data is small in the context of limited simulation budget};
    \item Kernel type for continuous variables: squared exponential;
    \item Nugget for Gaussian process: $10^{-6}$.
\end{itemize}
\item Bayesian Optimization:
\begin{itemize}
    \item Exploration factor for infill criterion: $k=2$;
    \item Expected violation threshold: $t_i=0.0001 \;\;\; \mbox{for} \: i=1,\ldots,n_g$ with $n_g$ the number of constraints. 
\end{itemize}
\item Benchmark:
\begin{itemize}
    \item Number of repetitions of each analytical problem: 10, and 5 repetitions of the engineering wing design problem;
    \item Run on a cluster of 12 cores of Skylake Intel® Xeon® Gold 6152 CPU.
\end{itemize}
\end{itemize}

\subsection{Rosenbrock problem}\label{Rosen_pb_matrix}
The correspondence between the values of the categorical variables and the values of  $a_q, b_q, e_q, f_q, j_q, k_q, r_q, s_q, t_q, u_q$ and $v_q$ are given in the following matrix:
\small 
\begin{equation}
     \bordermatrix{ & x^q_{1} & x^q_{2}  & \vr a_q & b_q & e_q & f_q & j_q & k_q & r_q & s_q & t_q & u_q & v_q\cr
      & 0 & 0 & \VR 100 &  1 & 0.7 & 2000 &1 & 0 & 1 & -1.2 & 0 & 0 & -1\cr
      & 0 & 1 & \VR 103 & 1.6 & 0.2 & 1950 &-1  & 0 & 1 & -0.2 & 0 & 0 & 0.97\cr
      & 1 & 0 & \VR 98 & 2 & 0.3 & 2100  &1& 0 & 1 & -0.7 &  0 & 0& 0.95 \cr
      & 1 & 1 & \VR 100 & 1.7 & 0.5 & 2020 &  1 & 0 & 1 & 0.15 & 0 & 0 & 1.1 \cr
      & 2 & 0 & \VR 95 & 4.7 & 1.5 & 1970 & 1 & 0.15 & 2 & 0 & 0.5 & 0 & -0.8 \cr
      & 2 & 1 & \VR 97 & 2.4 & 1.2 & 2100 & 1 & -0.55 & 2 & 0.4 & 0 & -0.8 & 0.7 \cr
      & 3 & 0 & \VR 103 & 1.7 & 2.5 & 2070 & -1 & -1.15 & 2 & 0 & -1.5 & 0 & 1.8 \cr
      & 3 & 1 & \VR 100 & 0.2 & 1 & 1890 & 1 & -1.3 & 2 & 1.4 & 0 & 0.8& -1.7 \cr
      & 4 & 0 & \VR 96 & 1.1 & 0.5 & 2140 & -1 & 0.5 & 2 & 0 & -2.3 & 0 & -0.8 \cr
      & 4 & 1 & \VR 104 & 1.5 & 2 & 1930 & -1 & 1.4 & 2 & -2.4 & 0 & 1.8 & -0.8 \cr
      & 5 & 0 & \VR 99 & 1.1 & 0.5 & 2140 & 1 & -1.5 & 2 & 0& 2 & 0 & -0.9 \cr
      & 5 & 1 & \VR 104 & 1.5 & 2 & 2030 & 1 & 1.8 & 2 & 0.4 & 0 & 1 & -0.3 } \qquad
\end{equation}
\normalsize
\subsection{Trid problem}\label{Trid_pb_matrix}
The correspondence between the values of the categorical variables and the values of  $a_q, b_q, c_q, e_q, f_q, j_q, k_q, r_q, s_q, t_q$ and $u_q$ are given in the following matrix:
\small 
\begin{equation}
     \bordermatrix{ & x^q_{1} & x^q_{2}  & \vr a_q & b_q & c_q & e_q & f_q & j_q & k_q & r_q & s_q & t_q & u_q \cr
      & 0 & 0 & \VR 1. &  1 & 1 & 1 & 1 & 0.7 & 1 & 1 & 1.5 & 1 & 0.4\cr
      & 1 & 0 & \VR 0.95 & 1 & 1.1 & 0.8  &1 & 0.4 & 1.1 & 1 & 1.9 & 1 & 0.1 \cr
      & 2 & 0 & \VR 1 & 1.3 & 0.97 & 1.1 & 0.8  & 0.1 & 1 & 0.9 &1.5. & 1.1 & 0.4\cr
      & 0 & 1 & \VR 1.1 & 0.7 & 1 & 1 &  1 & 0.7 & 1 & 1 & 0.7 & 1 & 1.4 \cr
      & 1 & 1 & \VR 0.7 & 0.5 & 0.4 & 1.5 & 1 & 1.7 & 0.7& 0.7 & 0.5 & 1 & 0.9 \cr
      & 2 & 1 & \VR 0.7 & 1 & 1.5 & 1 & 1.3 & 0.91 & 1 & 1 & 1.5 & 0.7 & 0.1 } \qquad
\end{equation}

\subsection{Styblinski-Tang problem problem}\label{Sty_pb_matrix}
The correspondence between the values of the categorical variables and the values of  $a_q, b_q, c_q, e_q, f_q, j_q, $ and $k_q$ are given in the following matrix:
\small 
\begin{equation}
     \bordermatrix{ & x^q_{1} & x^q_{2} & x^q_{3} & \vr a_q & b_q & c_q & e_q & f_q & j_q & k_q  \cr
      & 0 & 0 & 0 & \VR 1 &  16 & 5 & 1.2 & 0.7 & 3.5 & 0.7 \cr
      & 1 & 0 & 0 & \VR 1.1 & 18 & 6.1 & 1.4  & 0.9 & 3.8 & 0.2 \cr
      & 1 & 1 & 0 & \VR 0.95 & 17 & 4.9 & 1.7 & 1.3  & 2.8 & 0.7 \cr
      & 0 & 1 & 0 & \VR 0.94 & 12 & 6.9 & 1.4 &  0.2 & 1.4 & 0.2  \cr
      & 0 & 0 & 1 & \VR 0.75 & 10 & 7 & 2.2 & 1.7 & 1.5 & 0.5 \cr
      & 1 & 0 & 1 & \VR 1.2 & 19 & 4.2 & 1.5  & 2.9 & 1.4 & 1.2 \cr
      & 1 & 1 & 1 & \VR 0.97 & 12 & 1.9 & 0.7 & 2.3  & 3.8 & 0.4 \cr
      & 0 & 1 & 1 & \VR 1.1 & 18 & 4.2 & 1.9 &  0.7 & 2.7 & 0.4 } \qquad
\end{equation}

\subsection{Modified MAP-Elites algorithm for the solving of the auxiliary optimization problem involved in Bayesian Optimization}\label{Appendix_MAP_Elites_algo}
\begin{algorithm}[h!]
\caption{MAP-Elites algorithm for constrained optimization problem with mixed continuous, discrete and categorical variables. Objective, constraints and features functions correspond to the functions of the infill problem Eqs.(\ref{probaux_QD_1})-(\ref{probaux_QD_2})} \label{alg:MAP_Elites}
\begin{algorithmic}
\State 1) Initialize the archives $\mathcal{A}$ (empty maps for the objective function $\mathcal{A}_{Y}$, the feature functions $\mathcal{A}_{F_{t_{j=1,\dots,n}}}$, and the constraint functions $\mathcal{A}_{G_{i=1,\dots,n_g}}$ and the genome $\mathcal{A}_{X}$). 
\State 2) Generate the initial population of size $M$:
\State \hspace{0.5cm}  $\{\tilde{\mathbf{x}}_1= \left[\mathbf{x}_1^c,\mathbf{x}_1^d,\mathbf{x}_1^q\right]^T, \dots, \tilde{\mathbf{x}}_M=\left[\mathbf{x}_M^c,\mathbf{x}_M^d,\mathbf{x}_M^q\right]^T \}$ using LHS
\State 3) Evaluate  the initial population: 
\State \hspace{0.5cm}  $\{y_k = f(\tilde{\mathbf{x}}_k)  \}_{k\in\{1,\dots,M\}}$, $\{f_{t_j}(\tilde{\mathbf{x}}_k) \}_{k\in\{1,\dots,M\}; \; j\in\{1,\dots,n\}} \;\;$ and $\{g_{i}(\tilde{\mathbf{x}}_k)\}_{k\in\{1,\dots,M\}; \; i\in\{1,\dots,n_g\}} $ %\Comment{Evaluation of the objective function}
\State 4) Update of the archive 
\hspace{0.5cm}     \For{$k \gets 1$ to $M$}
    \hspace{1.0cm}     \If{$\forall i\in \{1,\dots,n_g\}, \;\; g_i(\tilde{\mathbf{x}}_k) \leq 0 $}
            \State \hspace{0.5cm} $\mathcal{A}_{Y} \gets  y_k$, \Comment{Update of objective function map}  
            \State \hspace{0.5cm} $\mathcal{A}_{{F_{t_j}}} \gets {f_{t_j}(\tilde{\mathbf{x}}_k)} \text{ for } j\in\{1,\dots,n\}$, \Comment{Update of features maps} 
            \State \hspace{0.5cm} $\mathcal{A}_{G_i} \gets g_{i}(\tilde{\mathbf{x}}_k) \text{ for } i\in\{1,\dots,n_g\}$ \Comment{Update of constraints maps}
            \State \hspace{0.5cm} $\mathcal{A}_{X} \gets \tilde{\mathbf{x}}_k$ \Comment{Update of genome map}
        \EndIf 
    \EndFor
\For{$\text{iter} \gets 1$ to \text{itermax}}
\hspace{0.5cm}     \State 5) $\mathbf{x'}_{\text{Mselec}} \gets \text{RandomSelection} \left( \mathcal{A}_{X} \right)$ of size $M'$
\hspace{0.5cm}     \State 6) $\mathcal{\tilde{X}}' \gets \text{RandomVariation} \left(\mathbf{x'}_{\text{Mselec}}\right)$
\hspace{0.5cm}     \State 7) Evaluate the new population and store results $\mathcal{\tilde{X}}'$ of size $M'$ in temporary archives: 
 \State \hspace{0.cm}  $\mathbf{Y}_{M'}=\left[f(\tilde{\mathbf{x}}'_k)  \right]^T_{k\in\{1,\dots,M'\}}$, \Comment{Objective temporary archive}  
  \State \hspace{0.cm}  $\mathbf{F_t}_{M'}[j]=\left[f_{t_j}(\tilde{\mathbf{x}}'_k) \right]_{k\in\{1,\dots,M'\}; \; j\in\{1,\dots,n\}} $ \hspace{-0.2cm}\Comment{Features temporary archives}  
   \State \hspace{0.cm}  $\mathbf{G}_{M'}[i]=\left[g_{i}(\tilde{\mathbf{x}}'_k)\right]_{k\in\{1,\dots,M'\}; \; i\in\{1,\dots,n_g\}} $ \hspace{-0.2cm}\Comment{Constraints temporary archives} 
\hspace{0.5cm}     \State 8) Update of the  archives 
\hspace{1.5cm}     \For{$k \gets 1$ to $M'$}
    \hspace{1.0cm}    \If{$ \mathbf{G}_{M'}[i][k] \leq 0 $ for $i\in\{1,\dots,n_g\}$ and $\left( \mathcal{A}_{F_{t}}(\mathbf{F_t}_{M'}[k]) = \emptyset \right.$ or $\left. \mathcal{A}_{Y}(\mathbf{F_t}_{M'}[k]) \geq \mathbf{Y}_{M'}[k]\right)$} \Comment{Check if valuable solution}  
            \State \hspace{-0.5cm} $\mathcal{A}_{Y} \gets  \mathbf{Y}_{M'}[k]$, \Comment{Update of objective function map}  
            \State \hspace{-0.5cm} $\mathcal{A}_{{F_{t_j}}} \gets \mathbf{F_t}_{M'}[j][k] \text{ for } j\in\{1,\dots,n\}$, \hspace{-0.5cm}\Comment{Update of features maps} 
            \State \hspace{-0.5cm} $\mathcal{A}_{G_i} \gets \mathbf{G}_{M'}[i][k] \text{ for } i\in\{1,\dots,n_g\}$, \hspace{-0.2cm}\Comment{Update of constraints map}
            \State \hspace{-0.5cm} $\mathcal{A}_{X} \gets \tilde{\mathbf{x}}'_k$ \Comment{Update of genome map}
    \EndIf 
    \EndFor
 \EndFor
 \State \textbf{Return} final maps: $\mathcal{A}_{Y}, \mathcal{A}_{F_{t_{j\in\{1,\dots,n\}}}}, \mathcal{A}_{G_{i\in\{1,\dots,n_g\}}}, \mathcal{A}_{X}$
\end{algorithmic}
\end{algorithm} 

 \bibliographystyle{elsarticle-num} 
 \bibliography{cas-refs}





% \end{thebibliography}
\end{document}